\documentclass{imanum}
\usepackage{amssymb,amsmath,amsfonts}
\usepackage[utf8]{inputenc}
\usepackage[T1]{fontenc}
\usepackage{mathrsfs}
\usepackage{enumerate}
\usepackage{tabularx}
\usepackage{graphicx}
\usepackage{color}
\usepackage[automark]{scrpage2}

\renewcommand{\a}{\alpha}
\renewcommand{\b}{\beta}
\newcommand{\g}{\gamma}
\renewcommand{\d}{\delta}
\renewcommand{\r}{\varrho}
\newcommand{\s}{\sigma}
\renewcommand{\t}{\tau}
\newcommand{\ps}{\psi}
\newcommand{\D}{\Delta}
\newcommand{\Ps}{\Psi}
\newcommand{\Om}{\Omega}
\newcommand{\sets}[1]{\mathbb{#1}}
\newcommand{\NN}{\sets{N}}
\newcommand{\RR}{\sets{R}}
\newcommand{\nice}[1]{\mathscr{#1}}
\newcommand{\nC}{\nice{C}}
\newcommand{\nE}{\nice{E}}

\newcommand{\nO}{\nice{O}}
\newcommand{\nEA}{\nice{E}_A}
\newcommand{\nEB}{\nice{E}_B}
\newcommand{\nEF}{\nice{E}_F}
\newcommand{\nEG}{\nice{E}_G}
\newcommand{\nSF}{\nice{S}_F}
\newcommand{\nLF}{\nice{L}_F}
\newcommand{\dd}{{\rm d}}
\newcommand{\ddt}{\tfrac{\dd}{\dd t}}
\newcommand{\ddtwo}{\tfrac{\dd^2}{\dd t^2}}
\newcommand{\ee}{{\rm e}}
\newcommand{\pa}{\partial}
\newcommand{\pat}{\pa_{t}}
\newcommand{\patt}{\pa_{tt}}
\newcommand{\pattt}{\pa_{ttt}}
\newcommand{\pax}{\pa_x}
\newcommand{\ve}[2]{\begin{pmatrix} #1 \\ #2 \end{pmatrix}}
\newcommand{\ma}[4]{\begin{pmatrix} #1 & \; #2 \\ #3 & \; #4 \end{pmatrix}}
\newcommand{\absbig}[1]{\big|#1\big|}
\newcommand{\norm}[2]{\|#2\|_{#1}}
\newcommand{\normbig}[2]{\big\|#2\big\|_{#1}}
\newcommand{\normBig}[2]{\Big\|#2\Big\|_{#1}}
\newcommand{\atil}{\widetilde{\a}}
\newcommand{\btil}{\widetilde{\b}}
\newcommand{\dtil}{\widetilde{\d}}
\newcommand{\pstil}{\widetilde{\ps}}
\newcommand{\ftil}{\widetilde{f}}
\newcommand{\bbtil}{\widetilde{b}}
\newcommand{\vtil}{\widetilde{v}}

\newcommand{\Dtil}{\widetilde{D}}
\newcommand{\Xtil}{\widetilde{X}}
\newcommand{\Ytil}{\widetilde{Y}}
\newcommand{\ahat}{\widehat{\a}}
\newcommand{\bhat}{\widehat{\b}}
\newcommand{\dhat}{\widehat{\d}}
\numberwithin{equation}{section}
\begin{document} 
\title{Efficient time integration methods based on operator splitting and application to the Westervelt equation} 
\shorttitle{Time-splitting methods for the Westervelt equation}
\author{\textsc{Barbara Kaltenbacher} \\
Alpen-Adria-Universität Klagenfurt, Institut für Mathematik, \\ Universitätsstraße~65--67, 9020 Klagenfurt, Austria. \\ E-mail: Barbara.Kaltenbacher@uni-klu.ac.at \\[2mm]
\textsc{Vanja Nikolic} \\ 
Alpen-Adria-Universität Klagenfurt, Institut für Mathematik, \\ Universitätsstraße 65--67, 9020~Klagenfurt, Austria. \\ E-mail: Vanja.Nikolic@aau.at \\[2mm]
\textsc{Mechthild Thalhammer} \\
Leopold--Franzens Universität Innsbruck, Institut für Mathematik, Technikerstra\ss{}e 13/VII, 6020 Innsbruck, Austria. \\ 
E-mail: Mechthild.Thalhammer@uibk.ac.at}
\shortauthorlist{B. Kaltenbacher, V. Nikolic, M. Thalhammer}
\maketitle
$\,$
\begin{abstract}%
{%
Efficient time integration methods based on operator splitting are introduced for the Westervelt equation, a nonlinear damped wave equation that arises in nonlinear acoustics as mathematical model for the propagation of sound waves in high intensity ultrasound applications.
For the first-order Lie--Trotter splitting method a global error estimate is deduced, confirming that the splitting method remains stable and that the nonstiff convergence order is retained in situations where the problem data are sufficiently regular. 
Fundamental ingredients in the stability and error analysis are regularity results for the Westervelt equation and related linear evolution equations of hyperbolic and parabolic type.
Numerical examples illustrate and complement the theoretical investigations.
} 
{Nonlinear evolution equations \and Westervelt equation \and Regularity of solutions \and Time-splitting methods \and Stability \and Local error expansion \and Convergence}
\end{abstract}
\addtocontents{toc}{\protect\setcounter{tocdepth}{1}}
\tableofcontents
\section{Introduction}
\paragraph{Scope of applications.}
High intensity ultrasound plays a crucial role in numerous practical settings ranging from medical treatment like lithotripsy or thermotherapy to industrial applications like ultrasound cleaning or welding and sonochemistry;
as a small selection, we mention the contributions~\cite{DreyerKrausBauerRiedlinger2000,Kaltenbacher2007} and refer to the literature cited therein. 
Numerical simulation of high intensity ultrasound propagation is a valuable tool for the design and improvement of high intensity ultrasound devices but poses major challenges due to the nonlinearity of the underlying partial differential equations.
Time integration methods for the equations of nonlinear acoustics have been investigated in~\cite{DreyerKrausBauerRiedlinger2000,KaltenbacherLandesHoffelnerSimkovics2002}; 
however, the use of transient simulations within the mathematical optimisation of high intensity ultrasound devices still seems to be beyond the scope of these approaches. 
On the other hand, operator splitting methods have proven to be efficient time integration methods in the context of other classes of nonlinear partial differential equations,  
see for instance~\cite{AuzingerHofstaetterKochThalhammer2013Nonlinear,DescombesThalhammer2012} and the references given therein.
This motivates the investigation of time-splitting methods for the equations of nonlinear acoustics.

\paragraph{Westervelt equation.}
In this work, we study time integration methods based on operator splitting for one of the classical models of nonlinear acoustics, the non-degenerate Westervelt equation
\begin{subequations}
\label{eq:Equation}
\begin{equation}
\label{eq:Equation1}
\begin{split}
&\patt \ps(x,t) - b \, \D \, \pat \ps(x,t) - c^2 \D \, \ps(x,t) \\
&\quad= \tfrac{\beta_a}{c^2} \, \pat \big(\pat \ps(x,t)\big)^2\,, \qquad (x,t) \in \Om \times (0,T] \subset \RR^d \times \RR_{> 0}\,,
\end{split}
\end{equation}
describing the propagation of the acoustic velocity potential $\ps: \overline{\Om} \times [0,T] \to \RR$, where $ c>0$ denotes the speed of sound, $b > 0$ the diffusivity of sound, and $\beta_a > 1$ the parameter of nonlinearity; 
for details on the physical background as well as the derivation of the Westervelt equation we refer to the original work~\cite{Westervelt1963}, see also~\cite{Kaltenbacher2007}.
For our purposes, it is advantageous to employ a reformulation of~\eqref{eq:Equation1} as nonlinear evolutionary system
\begin{equation}
\label{eq:Equation2}
\ddt \, u(t) = F\big(u(t)\big)\,, \qquad t \in (0,T]\,,
\end{equation}
\end{subequations}
for the abstract function $u: [0,T] \to X = X_1 \times X_2: t \mapsto u(t) = \Ps(\cdot,t) = (\ps(\cdot,t),\pat \ps(\cdot,t))$ with values in the underlying Banach space~$(X,\norm{X}{\cdot})$. 

\paragraph{Time-splitting methods.}
Our main concern is to introduce and analyse operator splitting methods for the time integration of the non-degenerate Westervelt equation~\eqref{eq:Equation}. 
Time-splitting methods for nonlinear evolution equations of the form~\eqref{eq:Equation2} utilise a natural decomposition of the defining operator
\begin{subequations}
\label{eq:EAB}
\begin{equation}
F = A + B
\end{equation}
and rely on the presumption that efficient numerical solvers are available for the associated subproblems 
\begin{equation}
\ddt \, v(t) = A\big(v(t)\big)\,, \quad \ddt \, w(t) = B\big(w(t)\big)\,, \qquad t \in (0,T]\,.
\end{equation}
\end{subequations}
We propose different decompositions for~\eqref{eq:Equation} which have in common that they require the numerical solution of a nonlinear diffusion equation in each time step. 

\paragraph{Convergence analysis.}
For the decomposition with the best performance in numerical tests we provide a detailed convergence analysis, adapting the approach exploited in~\cite{DescombesThalhammer2012} for the first-order Lie--Trotter splitting methods applied to nonlinear Schrödinger equations, see also~\cite{AuzingerHofstaetterKochThalhammer2013Nonlinear} for an extension to the second-order Strang splitting method.
A rigorous stability and error analysis of time-splitting methods for the Westervelt equation~\eqref{eq:Equation} is a complex task;
on the one hand, two unbounded nonlinear operators~$F = A + B$ are present, contrary to nonlinear Schrödinger equations comprising a linear differential operator and a nonlinear multiplication operator, and, on the other hand, the underlying Banach space $X = X_1 \times X_2$ is composed of two different function spaces reflecting the regularity of the solution components~$(\ps, \pat \ps)$.
As this considerably simplifies the local error analyis, we focus on the first-order Lie--Trotter splitting method, given by 
\begin{equation*}
u_1 = \nSF(h,u_0) = \nEB\big(h,\nEA(h,u_0)\big) \; \approx \; u(h) = \nEF\big(h,u(0)\big)\,, 
\end{equation*}
and indicate the extension to higher-order splitting methods; 
here, we denote by $\nEF, \nEA, \nEB$ the evolution operators associated with the Westervelt equation~\eqref{eq:Equation} and the subproblems~\eqref{eq:EAB}.
In situations where the problem data are sufficiently regular the obtained convergence result ensures that the Lie--Trotter splitting method remains stable and retains the nonstiff order of convergence. 
Fundamental ingredients in the stability and error analyis are regularity results for the Westervelt equation and related linear evolution equations of hyperbolic  and parabolic type. 
Numerical examples illustrate and complement the theoretical investigations. 
\section{Westervelt equation}
In this section, we introduce the initial-boundary value problem for the Westervelt equation and its abstract formulation as Cauchy problem.
A regularity result, needed as a basic ingredient in the convergence analysis of time-splitting methods for the Westervelt equation, is deduced in Section~\ref{sec:Regularity}. 

\paragraph{Westervelt equation.} 
We study the following initial-boundary value problem for a function $\ps: \overline{\Om} \times [0,T] \subset \RR^{d} \times \RR_{\geq 0} \to \RR: (x,t) \mapsto \ps(x,t)$ 
\begin{subequations}
\label{eq:WesterveltReformulations}
\begin{equation}
\label{eq:WesterveltOriginal}
\begin{cases}
&\patt \ps(x,t) - \a \, \D \, \pat \ps(x,t) - \b \, \D \, \ps(x,t) \\
&\quad= \g \, \pat \big(\pat \ps(x,t)\big)^2 = \d \, \pat \ps(x,t) \, \patt \ps(x,t)\,, \quad (x,t) \in \Om \times (0,T]\,, \\
&\ps(x,t) \; \text{ given}\,, \quad (x,t) \in \pa \Om \times [0,T]\,, \\
&\ps(x,0) \; \text{ given}\,, \quad \pat \ps(x,0) \; \text{ given}\,, \quad x \in \Om\,,
\end{cases} 
\end{equation}
involving positive constants $\a,\b,\g > 0$ and $\d = 2 \g > 0$.
With regard to the time integration by first- and second-order splitting methods, we restrict ourselves to situations where a sufficiently regular solution to~\eqref{eq:WesterveltOriginal} exists such that pointwise evaluations of the arising time and space derivatives of~$\ps$ are justified; 
in particular, we suppose the spatial domain and the prescribed initial data to be sufficiently regular. 

\paragraph{Reformulation (Non-degenerate case).} 
A regularity result for~\eqref{eq:WesterveltOriginal}, given in Section~\ref{sec:Regularity}, ensures non-degeneracy of the Westervelt equation;
that is, for any initial state of sufficiently small norm the relation $0 < 1 - \d \, \pat \ps(x,t) < \infty$ holds for all $(x,t) \in \overline{\Om} \times [0,T]$. 
As a consequence, the equivalent formulation 
\begin{equation}
\label{eq:WesterveltNondegenerate}
\begin{cases}
&\patt \ps(x,t) = \a \, \big(1- \d \, \pat \ps(x,t)\big)^{-1} \D \, \pat \ps(x,t) \\ 
&\qquad\qquad\qquad+\, \b \, \big(1- \d \, \pat \ps(x,t)\big)^{-1} \D \, \ps(x,t)\,, \quad (x,t) \in \Om \times (0,T]\,, \\
&\ps(x,t) \; \text{ given}\,, \quad (x,t) \in \pa \Om \times [0,T]\,, \\
&\ps(x,0) \; \text{ given}\,, \quad \pat \ps(x,0) \; \text{ given}\,, \quad x \in \Om\,,
\end{cases} 
\end{equation}
is obtained. 

\paragraph{Reformulation as first-order system.} 
We next rewrite~\eqref{eq:WesterveltNondegenerate} as a first-order system for the function $\Ps = (\Ps_1,\Ps_2) = (\ps, \pat \ps): \overline{\Om} \times [0,T] \to \RR^2$
\begin{equation}
\label{eq:WesterveltFirstOrder}
\begin{cases}
&\pat \Ps_1(x,t) = \Ps_2(x,t)\,, \quad (x,t) \in \Om \times (0,T]\,, \\
&\pat \Ps_2(x,t) = \a \, \big(1 - \d \, \Ps_2(x,t)\big)^{-1} \D \, \Ps_2(x,t) \\ 
&\qquad\qquad\qquad+\, \b \, \big(1- \d \, \Ps_2(x,t)\big)^{-1} \D \Ps_1(x,t)\big)\,, \quad (x,t) \in \Om \times (0,T]\,, \\
&\Ps_1(x,t) \; \text{ given}\,, \quad (x,t) \in \pa \Om \times [0,T]\,, \\ 
&\Ps_1(x,0) \; \text{ given}\,, \quad \Ps_2(x,0) \; \text{ given}\,, \quad x \in \Om\,.
\end{cases} 
\end{equation}
\end{subequations}

\paragraph{Abstract Cauchy problem.} 
In regard to the introduction and analysis of time-splitting methods for the Westervelt equation it is convenient to formulate~\eqref{eq:WesterveltReformulations} as abstract Cauchy problem for the function $u: [0,T] \to X: t \mapsto u(t) = \Ps(\cdot,t)$
\begin{subequations}
\label{eq:Problem}
\begin{equation}
\begin{cases}
&\ddt \, u(t) = F\big(u(t)\big)\,, \quad t \in (0,T]\,, \\ 
&u(0) \text{ given}\,,
\end{cases} 
\end{equation}
where $(X,\norm{X}{\cdot})$ denotes the underlying Banach space and the nonlinear operator $F: D(F) \to X$ is given by 
\begin{equation}
\begin{gathered}
F(v) = \ve{v_2}{\atil(v_2) \, \D v_2 + \btil(v_2) \, \D v_1}\,, \qquad v = (v_1,v_2) \in D(F)\,, \\
\atil(v_2) = \a \, \big(1 - \d v_2\big)^{-1}\,, \qquad \btil(v_2) = \b \, \big(1 - \d v_2\big)^{-1}\,.
\end{gathered}
\end{equation}
\end{subequations}
According to the situation under consideration, the domain $D(F) \subset X$ is chosen such that it reflects the regularity requirements on the solution to the Westervelt equation as well as the imposed boundary conditions. 
It is notable that~\eqref{eq:Problem} can be cast into the form of a quasilinear problem, since 
\begin{equation*}
F(v) = \ma{0}{I}{\btil(v_2) \, \D}{\;\atil(v_2) \, \D} \ve{v_1}{v_2}\,,
\end{equation*}
where~$I$ denotes the identity operator.
\section{Time-splitting methods}
In this section, we introduce the general format of exponential operator splitting methods for the time integration of nonlinear evolution equations and specify the different decompositions considered for the Westervelt equation. 
Numerical examples comparing the accuracy and efficiency of the resulting time dicretisations are found in Section~\ref{sec:Numerics}.
\subsection{Time-splitting methods for nonlinear evolution equations}
\paragraph{Nonlinear evolution equation.}
In accordance with the reformulation of the initial-boundary value problem for the Westervelt equation~\eqref{eq:WesterveltReformulations} as abstract Cauchy problem~\eqref{eq:Problem}, we consider the initial value problem 
\begin{subequations}
\label{eq:ProblemABEF}
\begin{equation}
\label{eq:ProblemAB}
\begin{cases}
&\ddt \, u(t) = F\big(u(t)\big) = A\big(u(t)\big) + B\big(u(t)\big)\,, \quad t \in (0,T]\,, \\
&u(0) \text{ given}\,, 
\end{cases}
\end{equation}
assuming that the nonlinear operators $A: D(A) \to X$ and $B: D(B) \to X$ are chosen such that the intersection $D(A) \cap D(B)$ coincides with the domain of the defining operator $F: D(F) \to X$.
For the following investigations it is useful to introduce the evolution operator associated with~\eqref{eq:ProblemAB} 
\begin{equation}
\label{eq:EF}
u(t) = \nEF\big(t,u(0)\big)\,, \qquad t \in [0,T]\,.
\end{equation}
\end{subequations}

\paragraph{Subproblems.}
Exponential operator splitting methods for the time integration of the initial value problem~\eqref{eq:ProblemAB} rely on the numerical solution of the subproblems 
\begin{subequations}
\label{eq:Subproblems}
\begin{equation}
\begin{split}
&\text{Subproblem A:} \qquad
\begin{cases}
&\ddt \, v(t) = A\big(v(t)\big)\,, \quad t \in (0,T]\,, \\ 
&v(0) \text{ given}\,, 
\end{cases} \\
&\text{Subproblem B:} \qquad
\begin{cases}
&\ddt \, w(t) = B\big(w(t)\big)\,, \quad t \in (0,T]\,, \\
&w(0) \text{ given}\,.
\end{cases} 
\end{split}
\end{equation}
In accordance with~\eqref{eq:EF} the associated evolution operators are given by 
\begin{equation}
v(t) = \nEA\big(t,v(0)\big)\,, \quad w(t) = \nEB\big(t,w(0)\big)\,, \qquad t \in [0,T]\,.
\end{equation}
\end{subequations}

\paragraph{Time-discrete solution.}
For an initial approximation $u_0 \approx u(0)$ and a sequence of time grid points $0 = t_0 < t_1 < \cdots < t_N \leq T$ with corresponding time stepsizes $h_n = t_{n+1} - t_n$ for $0 \leq n \leq N-1$ the time-discrete solution to~\eqref{eq:ProblemABEF} is determined through a recurrence relation of the form 
\begin{subequations}
\label{eq:Splitting}
\begin{equation}
u_n = \nSF(h_{n-1}, u_{n-1}) \; \approx \; u(t_n) = \nEF\big(h_{n-1}, u(t_{n-1})\big)\,, \qquad 1 \leq n \leq N\,. 
\end{equation}
The numerical evolution operator associated with a general splitting method of (nonstiff) order $p \geq 1$ can be cast into the format 
\begin{equation}
\begin{gathered}
u_n = \nSF(h_{n-1}, u_{n-1}) = U_{n-1,s}\,, \\
U_{n-1,j} = \nE_{B_j}\big(h_{n-1},\nE_{A_j}(h_{n-1},U_{n-1,j-1})\big)\,, \quad 1 \leq j \leq s\,, \quad U_{n-1,0} = u_{n-1}\,,
\end{gathered}
\end{equation}
involving certain coefficients $(a_j,b_j)_{j=1}^{s}$. 
Here, we employ the abbreviations 
\begin{equation}
A_j = a_j A\,, \quad B_j = b_j B\,, \qquad 1 \leq j \leq s\,;
\end{equation}
\end{subequations}
due to the fact that the considered evolution equation is autonomous, a scaling of the operators~$A, B$ corresponds to a scaling in time. 

\paragraph{Lie--Trotter and Strang splitting methods.}
The Lie--Trotter splitting method of (non\-stiff) order $p = 1$, given by 
\begin{equation}
\label{eq:Lie}
\nSF(t,v) = \nEB\big(t,\nEA(t,v)\big) \quad \text{or} \quad \nSF(t,v) = \nEA\big(t,\nEB(t,v)\big)\,, 
\end{equation}
can be cast into the format~\eqref{eq:Splitting} for the choices 
\begin{equation*}
\begin{gathered}
s = 1: \qquad a_1 = 1\,, \quad b_1 = 1\,, \qquad \text{or} \\
s = 2: \qquad a_1 = 0\,, \quad b_1 = 1\,, \quad a_2 = 1\,, \quad b_2 = 0\,,
\end{gathered}
\end{equation*}
respectively. 
The widely used Strang splitting method of (nonstiff) order $p = 2$, given~by 
\begin{equation}
\label{eq:Strang}
\begin{gathered}
\nSF(t,v) = \nE_{\frac{1}{2} A}\big(t,\nEB\big(t,\nE_{\frac{1}{2} A}(t,v)\big)\big) \quad \text{or} \\
\nSF(t,v) = \nE_{\frac{1}{2} B}\big(t,\nEA\big(t,\nE_{\frac{1}{2} B}(t,v)\big)\big)
\end{gathered}
\end{equation}
results for the choices 
\begin{equation*}
\begin{gathered}
s = 2: \qquad a_1 = \tfrac{1}{2}\,, \quad b_1 = 1\,, \quad a_2 = \tfrac{1}{2}\,, \quad b_2 = 0\,, \qquad \text{or} \\
s = 2: \qquad a_1 = 0\,, \quad b_1 = \tfrac{1}{2}\,, \quad a_2 = 1\,, \quad b_2 = \tfrac{1}{2}\,, 
\end{gathered}
\end{equation*}
respectively. 
\subsection{Time-splitting methods for the Westervelt equation}
\label{sec:Decompositions}
In the following, we propose different decompositions 
\begin{equation*}
F = A + B
\end{equation*}
of the nonlinear operator defining the Westervelt equation~\eqref{eq:WesterveltReformulations}--\eqref{eq:Problem} and discuss the computational effort for the numerical solution of the associated subproblems~\eqref{eq:Subproblems}.  
\subsubsection{Decomposition~I}
\paragraph{Decomposition.}
A first decomposition involves the nonlinear operators 
\begin{equation*}
\begin{gathered}
A(v) = \ve{v_2}{\atil(v_2) \, \D v_2}\,, \qquad B(v) = \ve{0}{\btil(v_2) \, \D v_1}\,, \\
\atil(v_2) = \a \, \big(1 - \d v_2\big)^{-1}\,, \qquad \btil(v_2) = \b \, \big(1 - \d v_2\big)^{-1}\,.
\end{gathered}
\end{equation*}

\paragraph{Subproblems.}
The resolution of the subproblem associated with~$A$
\begin{equation*}
\begin{cases}
&\pat \Ps_1(x,t) = \Ps_2(x,t)\,, \\
&\pat \Ps_2(x,t) = \a \, \big(1 - \d \, \Ps_2(x,t)\big)^{-1} \D \Ps_2(x,t)\,,
\end{cases}
\end{equation*}
amounts to the numerical solution of a nonlinear diffusion equation for the second component~$\Ps_2 = \pat \ps$
\begin{equation*}
\pat \Ps_2(x,t) = \a \, \big(1 - \d \, \Ps_2(x,t)\big)^{-1} \D \Ps_2(x,t)\,;
\end{equation*}
the first component $\Ps_1 = \ps$ is then retained by integration
\begin{equation*}
\Ps_1(x,t) = \Ps_1(x,0) + \int_{0}^{t} \Ps_2(x,\t) \; \dd\t\,.
\end{equation*}
For the subproblem associated with~$B$
\begin{equation*}
\begin{cases}
&\pat \Ps_1(x,t) = 0\,, \\
&\pat \Ps_2(x,t) = \b \, \big(1 - \d \, \Ps_2(x,t)\big)^{-1} \D \Ps_1(x,t)\,,
\end{cases} 
\end{equation*}
the first component remains constant on the considered time interval 
\begin{equation*}
\Ps_1(x,t) = \Ps_1(x,0)\,.
\end{equation*}
Consequently, the second component is a solution to 
\begin{equation*}
\pat \Ps_2(x,t) = \b \, \big(1 - \d \, \Ps_2(x,t)\big)^{-1} \D \Ps_1(x,0) 
\end{equation*}
with explicit representation 
\begin{equation*}
\Ps_2(x,t) = \tfrac{1}{\d} \, \Big(1 - \sqrt{\big(1 - \d \, \Ps_2(x,0)\big)^2 - 2 \, \b \d \, t \, \D \Ps_1(x,0)}\;\Big)\,. 
\end{equation*}
We note that in the non-degenerate case the relation $\d \, \Ps_2(x,0) < 1$ holds\,; 
thus, the other admissable choice for a solution to the subproblem leads to a contradiction when evaluating at $t = 0$.  
Provided that the time increment $t > 0$ is chosen sufficiently small, it is ensured that $(1 - \d \, \Ps_2(x,0))^2 - 2 \, \b \d \, t \, \D \Ps_1(x,0) > 0$ and hence $\Ps_2(x,t) \in \RR$. 
\subsubsection{Decomposition~II}
\paragraph{Decomposition.}
A second decomposition involves the nonlinear operators 
\begin{equation*}
A(v) = \ve{\frac{1}{2} \, v_2}{\atil(v_2) \, \D v_2}\,, \qquad B(v) = \ve{\frac{1}{2} \, v_2}{\btil(v_2) \, \D v_1}\,.
\end{equation*}

\paragraph{Subproblems.}
The resolution of the subproblem associated with~$A$
\begin{equation*}
\begin{cases}
&\pat \Ps_1(x,t) = \frac{1}{2} \, \Ps_2(x,t)\,, \\
&\pat \Ps_2(x,t) = \a \, \big(1 - \d \, \Ps_2(x,t)\big)^{-1} \D \Ps_2(x,t)\,,
\end{cases}
\end{equation*}
requires the numerical solution of a nonlinear diffusion equation for the second component
\begin{equation*}
\pat \Ps_2(x,t) = \a \, \big(1 - \d \, \Ps_2(x,t)\big)^{-1} \D \Ps_2(x,t)\,;
\end{equation*}
the first component is then obtained by integration 
\begin{equation*}
\Ps_1(x,t) = \Ps_1(x,0) + \tfrac{1}{2} \int_{0}^{t} \Ps_2(x,\t) \; \dd\t\,.
\end{equation*}
The resolution of the subproblem associated with~$B$
\begin{equation*}
\begin{cases}
&\pat \Ps_1(x,t) = \frac{1}{2} \, \Ps_2(x,t)\,, \\
&\pat \Ps_2(x,t) = \b \, \big(1 - \d \, \Ps_2(x,t)\big)^{-1} \D \Ps_1(x,t)\,, 
\end{cases} 
\end{equation*}
amounts to the numerical solution of a nonlinear wave equation for the first component 
\begin{equation*}
\patt \Ps_1(t) = \tfrac{1}{2} \, \b \, \big(1 - 2 \, \d \, \pat \Ps_1(x,t)\big)^{-1} \D \Ps_1(x,t)\,;
\end{equation*}
the second component is then given by 
\begin{equation*}
\Ps_2(x,t) = 2 \, \pat \Ps_1(x,t)\,.
\end{equation*}
\subsubsection{Decomposition~III}
\paragraph{Decomposition.}
A third decomposition involves the nonlinear operators 
\begin{equation*}
A(v) = \ve{0}{\atil(v_2) \, \D v_2}\,, \qquad B(v) = \ve{v_2}{\btil(v_2) \, \D v_1}\,.
\end{equation*}

\paragraph{Subproblems.}
The resolution of the subproblem associated with~$A$ 
\begin{equation*}
\begin{cases}
&\pat \Ps_1(x,t) = 0\,, \\
&\pat \Ps_2(x,t) = \a \, \big(1 - \d \, \Ps_2(x,t)\big)^{-1} \D \Ps_2(x,t)\,,
\end{cases}
\end{equation*}
amounts to the numerical solution of a nonlinear diffusion equation for the second component 
\begin{equation*}
\pat \Ps_2(x,t) = \a \, \big(1 - \d \, \Ps_2(x,t)\big)^{-1} \D \Ps_2(x,t)\,,
\end{equation*}
whereas the first component remains constant in time
\begin{equation*}
\Ps_1(x,t) = \Ps_1(x,0)\,.
\end{equation*}  
The resolution of the subproblem associated with~$B$
\begin{equation*}
\begin{cases}
&\pat \Ps_1(x,t) = \Ps_2(x,t)\,, \\
&\pat \Ps_2(x,t) = \b \, \big(1 - \d \, \Ps_2(x,t)\big)^{-1} \D \Ps_1(x,t)\,, 
\end{cases} 
\end{equation*}
requires the numerical solution of a nonlinear wave equation for the first component 
\begin{equation*}
\patt \Ps_1(t) = \b \, \big(1 - \d \, \pat \Ps_1(x,t)\big)^{-1} \D \Ps_1(x,t)\,; 
\end{equation*}
the second component is then given by 
\begin{equation*}
\Ps_2(x,t) = \pat \Ps_1(x,t)\,.
\end{equation*}
\subsubsection{Decomposition~IV}
\paragraph{Decomposition.}
A fourth decomposition involves the nonlinear operators 
\begin{equation*}
A(v) = \ve{0}{\atil(v_2) \, \D v_2 + \d \, v_2 \, \btil(v_2) \, \D v_1}\,, \qquad 
B(v) = \ve{v_2}{\b \, \D v_1}\,;  
\end{equation*}
we note that the identity $\b + \d \, v_2 \, \btil(v_2) = \btil(v_2)$ holds. 

\paragraph{Subproblems.}
For the subproblem associated with~$A$ 
\begin{equation*}
\begin{cases}
&\pat \Ps_1(x,t) = 0\,, \\
&\pat \Ps_2(x,t) = \a \, \big(1 - \d \, \Ps_2(x,t)\big)^{-1} \D \Ps_2(x,t) \\
&\qquad\qquad\qquad+ \; \b \, \d \, \Ps_2(x,t) \, \big(1 - \d \, \Ps_2(x,t)\big)^{-1} \D \Ps_1(x,t)\,,
\end{cases}
\end{equation*}
the first component remains constant in time
\begin{equation*}
\Ps_1(x,t) = \Ps_1(x,0)\,, 
\end{equation*}  
whereas the computation of the second component requires the numerical solution of a nonlinear reaction-diffusion equation 
\begin{equation*}
\pat \Ps_2(x,t) = \a \, \big(1 - \d \, \Ps_2(x,t)\big)^{-1} \D \Ps_2(x,t) + \b \, \d \, \Ps_2(x,t) \, \big(1 - \d \, \Ps_2(x,t)\big)^{-1} \D \Ps_1(x,0)\,.
\end{equation*}
The resolution of the subproblem associated with~$B$ 
\begin{equation*}
\begin{cases}
&\pat \Ps_1(x,t) = \Ps_2(x,t)\,, \\
&\pat \Ps_2(x,t) = \b \, \D \Ps_1(x,t)\,,
\end{cases} 
\end{equation*}
amounts to the numerical solution of a linear wave equation for the first component 
\begin{equation*}
\patt \Ps_1(x,t) = \b \, \D \Ps_1(x,t)\,;
\end{equation*}
the second component is then given by 
\begin{equation*}
\Ps_2(x,t) = \pat \Ps_1(x,t)\,.
\end{equation*}
\subsubsection{Computational effort}
\paragraph{Subproblem~A.}
In all decompositions the realisation of the subproblem associated with the nonlinear operator~$A$ requires the numerical solution of a nonlinear diffusion equation for the second component $\Ps_2 = \pat \ps$, where in case of Decomposition~IV the nonlinear diffusion equation involves an additional zero-order term; 
the first component $\Ps_1 = \ps$ is obtained by integration or remains constant in time. 

\paragraph{Subproblem~B.}
In Decomposition~I the realisation of the subproblem associated with the nonlinear operator~$B$ reduces to the application of the Laplace operator to the first component of the imposed initial state; 
the second component is then computed through an explicit representation. 
Contrary, in Decomposition~II-III the realisation of the subproblem amounts to the numerical solution of a nonlinear wave equation for the first component; 
the second component is then given by the time derivative~$\pat \Ps_1$.
We point out that the solvability of the arising nonlinear wave equations is an open question and that even local solvability in time cannot be guaranteed.
As alternative we thus introduce Decomposition~IV involving a linear wave equation.
\section{Stability and error analysis}
\label{sec:Convergence}
In this section, we describe the general approach for a convergence analysis of time-splitting methods, employing a framework of abstract Cauchy problems;
details on the specialisation to the Westervelt equation are given in Section~\ref{sec:GlobalErrorEstimate}. 

\paragraph{Uniform time grid.}
For notational simplicity, we henceforth restrict ourselves to constant time increments $h > 0$ with associated equidistant time grid points given by $t_n = n h$ for $0 \leq n \leq N$.
However, it is straightforward to extend the arguments to a non-uniform time grid with stepsize ratios bounded from above and below; 
in the global error estimate~\eqref{eq:GlobalErrorEstimateuN} the constant time increment $h > 0$ is then replaced with the maximum stepsize~$h_{\max} = \max\{h_{n}: 0 \leq n \leq N-1\}$. 

\paragraph{Global and local error.}
In order to deduce a convergence result for time-splitting methods~\eqref{eq:Splitting} applied to nonlinear evolution equations~\eqref{eq:ProblemABEF}, we utilise the  telescopic identity 
\begin{subequations}
\label{eq:GlobalError}
\begin{equation}
\begin{split}
u_N - u(t_N) 
&= \nSF^N\big(h,u_0\big) - \nSF^N\big(h,u(0)\big) \\
&\qquad+ \sum_{n=0}^{N-1} \bigg(\nSF^{N-n-1}\Big(h,\nLF\big(h,u(t_n)\big) + \nEF\big(h,u(t_n)\big)\Big) \\
&\qquad\qquad\qquad- \nSF^{N-n-1}\Big(h,\nEF\big(h,u(t_n)\big)\Big)\bigg)
\end{split}
\end{equation}
relating the global error to compositions of the splitting operator and local errors 
\begin{equation}
\begin{gathered}
\nSF^{n+1}(h,v) = \nSF\big(h,\nSF^{n}(h,v)\big)\,, \quad 0 \leq n \leq N-1\,, \qquad \nSF^0(h,v) = v\,, \\
\nLF\big(h,u(t_n)\big) = \nSF\big(h,u(t_n)\big) - \nEF\big(h,u(t_n)\big)\,, \quad 0 \leq n \leq N-1\,. 
\end{gathered}
\end{equation}
\end{subequations}

\paragraph{Approach.}
Our aim is to deduce an estimate for the global error~\eqref{eq:GlobalError} with respect to the norm of the underlying Banach space~$X$ or a suitable subspace $\Xtil \subset X$, respectively, such that the expected dependence on the time stepsize is retained under appropriate regularity requirements on the problem data.  
In the context of the Westervelt equation, we shall make use of the fact that the regularity result stated in Section~\ref{sec:Regularity} implies 
\begin{subequations}
\label{eq:GlobalErrorEstimate}
\begin{equation}
u(0) \in D\,, \; \normbig{D}{u(0)} \leq C_0 \quad \Longrightarrow \quad u(t) \in D\,, \; \normbig{D}{u(t)} \leq C\,, \quad t \in [0,T]\,, 
\end{equation}
for certain constants $C_0, C > 0$; 
the subspace $D \subset \Xtil \subset X$, chosen accordingly to the (nonstiff) order of the splitting method, captures additional regularity and compatibility requirements.    
Regularity results for the arising subproblems and associated variational equations shall ensure stability of the splitting procedure, that is, boundedness of compositions of the splitting operator 
\begin{equation}
\label{eq:GlobalErrorStability}
\normbig{\Xtil}{\nSF^n(h,v) - \nSF^n(h,\vtil)} \leq \ee^{C t_n} \normbig{\Xtil}{v - \vtil}\,, \qquad 1 \leq n \leq N\,, 
\end{equation}
provided that the arguments $v, \vtil \in \Dtil$ remain bounded in a certain subspace $D \subset \Dtil \subset \Xtil \subset X$. 
Concerning a local error analysis, we restrict ourselves to the least technical case, the first-order Lie--Trotter splitting method given by 
\begin{equation*}
\nSF(h,v) = \nEB\big(h,\nEA(h,v)\big)\,, 
\end{equation*}
see also~\eqref{eq:Lie};
here, a suitable expansion yields the representation  
\begin{equation}
\label{eq:L}
\begin{split}
\nLF(h,v) 
&= \int_{0}^{h} \int_{0}^{\t_1} \pa_{2} \nEF\big(h-\t_1,\nSF(\t_1,v)\big) \, \pa_{2} \nEB(\t_1,w) \, \big(\pa_{2} \nEB(\t_2,w)\big)^{-1} \\
&\qquad\qquad\quad \times \big[B,A\big]\big(\nEB(\t_2,w)\big) \Big\vert_{w = \nEA(\t_1,v)} \; \dd\t_2 \, \dd\t_1\,,
\end{split}
\end{equation}
see~\cite{AuzingerHofstaetterKochThalhammer2013Nonlinear,DescombesThalhammer2012} and formula~\eqref{eq:LLie} in Section~\ref{sec:LocalError}.
Bounds for the evolution operators~$\nEF, \nEA, \nEB$, their Fréchet derivatives with respect to the second argument~$\pa_2 \nEF, \pa_2 \nEA, \pa_2 \nEB$, and the first Lie-commutator~$[A,B]$ shall lead to the local error estimate
\begin{equation}
p = 1: \qquad 
\normbig{\Xtil}{\nLF(h,v)} \leq C \, h^{p+1}\,. 
\end{equation}
Altogether, these ingredients allow to establish the global error estimate 
\begin{equation}
\label{eq:GlobalErrorEstimateuN}
\begin{split}
\normbig{\Xtil}{u_N - u(t_N)} 
&\leq \normbig{\Xtil}{\nSF^N\big(h,u_0\big) - \nSF^N\big(h,u(0)\big)} \\
&\qquad+ \sum_{n=0}^{N-1} \Big\|\nSF^{N-n-1}\Big(h,\nLF\big(h,u(t_n)\big) + \nEF\big(h,u(t_n)\big)\Big)\\ 
&\qquad\qquad\qquad- \nSF^{N-n-1}\Big(h,\nEF\big(h,u(t_n)\big)\Big)\Big\|_{\Xtil} \\
&\leq C \, \bigg(\normbig{\Xtil}{u_0 - u(0)} + \sum_{n=1}^{N} \normbig{\Xtil}{\nLF\big(h,u(t_{n-1})\big)}\bigg) \\
&\leq C \, \Big(\normbig{\Xtil}{u_0 - u(0)} + h^{p}\Big)
\end{split}
\end{equation}
\end{subequations}
with constant depending in particular on the bound for the exact solution values with respect to the norm in~$D$, the bound for the initial approximation~$u_0$ in $\Dtil$, and the final time; 
this shows that the nonstiff order of concergence is retained, provided that the prescribed initial state is sufficiently regular.
\section{Global error estimate}
\label{sec:GlobalErrorEstimate}
In this section, we establish a global error estimate for the Lie--Trotter splitting method applied to the Westervelt equation.
We include a detailed stability and error analysis for Decomposition~I showing the best performance in numerical tests, see Section~\ref{sec:Decompositions}\,. 
In accordance with Section~\ref{sec:Regularity} we study the Westervelt equation in three space dimensions subject to homogeneous Dirichlet boundary conditions on a regular domain; 
we consider this to be the least technical case, as far as the boundary conditions are concerned, and the most relevant case, as far as the space dimension is concerned. 
Fundamental auxiliary results for the proof of Theorem~\ref{thm:MainResult} are deduced in Sections~\ref{sec:DerivativesCommutator}--\ref{sec:EvolutionOperators}.  
\subsection{Fréchet derivative and Lie-commutator}
\label{sec:DerivativesCommutator}
In the following, we specify the Fréchet derivative of the operator defining~\eqref{eq:WesterveltReformulations}--\eqref{eq:Problem} and determine the first Lie-commutator arising in the local error representation~\eqref{eq:L}.

\paragraph{Notation and auxiliary estimates.}
The product space $X = X_1 \times X_2$ of two Banach spaces is endowed with the norm $\norm{X}{(x_1,x_2)} = \norm{X_1}{x_1} + \norm{X_2}{x_2}$.
We employ standard abbreviations for Sobolev spaces such as $(W^{k,p}(\Om),\norm{W^{k,p}}{\cdot})$ and $H^{k}(\Om) = W^{k,2}(\Om)$ for $k \in \NN_{\geq 0}$ and $p \in \NN_{\geq 1}$, see for instance~\cite{AdamsFournier}; 
for notational simplicity, in the norms we do not indicate the dependence on the domain. 
For convenience, we recall Young's inequality 
\begin{subequations}
\label{eq:AuxiliaryEstimates}
\begin{equation}
a b \leq \tfrac{a^p}{p} + \tfrac{b^q}{q}\,, \qquad a, b \geq 0\,, \quad 1 < p,q < \infty\,, \quad \tfrac{1}{p} + \tfrac{1}{q} = 1\,.
\end{equation}
By means of Hölder's inequality 
\begin{equation}
\normbig{L^1}{f_1 \, f_2} \leq \normbig{L^p}{f_1} \, \normbig{L^q}{f_2}\,, \qquad 1 \leq p,q \leq \infty\,, \quad \tfrac{1}{p} + \tfrac{1}{q} = 1\,, 
\end{equation}
and the continuous embeddings $H^{2}(\Om) \hookrightarrow \nC(\Om)$ as well as $H^{1}(\Om) \hookrightarrow L^4(\Om), L^6(\Om)$, valid for any bounded regular domain $\Om \subset \RR^3$, the relations
\begin{equation}
\begin{gathered}
\normbig{L^2}{f_1 \, f_2} \leq \normbig{L^{\infty}}{f_1} \normbig{L^2}{f_2} \leq C \, \normbig{H^2}{f_1} \normbig{L^2}{f_2}\,, \\
\normbig{L^2}{f_1 \, f_2} \leq \normbig{L^4}{f_1} \normbig{L^4}{f_2} \leq C \, \normbig{H^1}{f_1} \normbig{H^1}{f_2}\,, \\
\normbig{L^2}{f_1 \, f_2} \leq \normbig{L^3}{f_1} \normbig{L^6}{f_2} \leq C \, \normbig{H^1}{f_1} \normbig{H^1}{f_2}\,, \\
\normbig{L^2}{f_1 \, f_2 \, f_3} \leq \normbig{L^6}{f_1} \normbig{L^6}{f_2} \normbig{L^6}{f_3} \leq C \, \normbig{H^1}{f_1} \normbig{H^1}{f_2} \normbig{H^1}{f_3}\,, \end{gathered}
\end{equation}
follow. 
Moreover, we make use of the fact that~$H^{k+2}(\Om)$ forms an algebra for arbitrary exponents $k \in \NN_{\geq 0}$, that is, the relation   
\begin{equation}
\normbig{H^{k+2}}{f_1 \, f_2} \leq C \, \norm{H^{k+2}}{f_1} \, \norm{H^{k+2}}{f_2}\,, \qquad k \in \NN_{\geq 0}\,,
\end{equation}
holds. 
This in particular implies  
\begin{equation}
\normbig{H^{k+2}}{\nabla f_1 \cdot \nabla f_2} \leq C \, \norm{H^{k+3}}{f_1} \, \norm{H^{k+3}}{f_2}\,, \qquad k \in \NN_{\geq 0}\,.
\end{equation}
\end{subequations}

\paragraph{Defining operators.}
For the convenience of the reader, we recall the definitions of the nonlinear operators associated with Decomposition~I
\begin{equation*}
\begin{gathered}
F(v) = A(v) + B(v)\,, \\
A(v) = \ve{v_2}{\atil(v_2) \, \D v_2}\,, \qquad B(v) = \ve{0}{\btil(v_2) \, \D v_1}\,, \\ 
\atil(v_2) = \a \, \big(1 - \d v_2\big)^{-1}\,, \qquad \btil(v_2) = \b \, \big(1 - \d v_2\big)^{-1}\,. 
\end{gathered}
\end{equation*}
Under the non-degeneracy condition $0 < \underline{\nu} \leq 1 - \d \, v_2(x) \leq \overline{\nu} < \infty$ for all $x \in \overline{\Om}$, justified by Theorem~\ref{thm:TheoremRegularity} in Section \ref{sec:Regularity}, the auxiliary results in~\eqref{eq:AuxiliaryEstimates} imply 
\begin{equation*}
\normbig{H^{k+2}}{\atil(v_2)} + \normbig{H^{k+2}}{\btil(v_2)} \leq C \big(\norm{H^{k+2}}{v_2}\big)\,, \qquad k \in \NN_{\geq 0}\,, 
\end{equation*}
and further lead to the estimates 
\begin{equation*}
\begin{gathered}
\normbig{L^2}{\atil(v_2) \, \D v_2} \leq C\big(\norm{H^{2}}{v_2}\big)\,, \\
\normbig{L^2}{\btil(v_2) \, \D v_1} \leq C\big(\norm{H^2}{v_1}, \norm{H^{2}}{v_2}\big)\,, \\
\normbig{H^{k+2}}{\atil(v_2) \, \D v_2} \leq C\big(\norm{H^{k+4}}{v_2}\big)\,, \qquad k \in \NN_{\geq 0}\,, \\
\normbig{H^{k+2}}{\btil(v_2) \, \D v_1} \leq C\big(\norm{H^{k+4}}{v_1}, \norm{H^{k+2}}{v_2}\big)\,, \qquad k \in \NN_{\geq 0}\,, 
\end{gathered}
\end{equation*}
with constants~$C(\cdot)$ depending on bounds for the arising norms of the solution components~$v_1, v_2$.

\paragraph{Fréchet derivatives.}
The Fréchet derivatives of the nonlinear operators associated with Decomposition~I are given by 
\begin{equation*}
\begin{gathered}
F'(v) = A'(v) + B'(v)\,, \\
A'(v) = \ma{0}{I}{0}{\atil(v_2) \, \D + \atil'(v_2) \, \D v_2}\,, \qquad B'(v) = \ma{0}{0}{\btil(v_2) \, \D}{\btil'(v_2) \, \D v_1}\,, \\
\atil'(v_2) = \a \, \d \, \big(1 - \d v_2\big)^{-2} = \tfrac{\d}{\a} \, \big(\atil(v_2)\big)^2\,, \qquad
\btil'(v_2) = \b \, \d \, \big(1 - \d v_2\big)^{-2} = \tfrac{\d}{\b} \, \big(\btil(v_2)\big)^2\,;
\end{gathered}
\end{equation*}
more precisely, application to an element $w = (w_1,w_2)$ yields 
\begin{equation*}
\begin{gathered}
A'(v) \, w = \ve{w_2}{\atil(v_2) \, \D w_2 + \atil'(v_2) \, \D v_2 \, w_2}\,, \\
B'(v) \, w = \ve{0}{\btil(v_2) \, \D w_1 + \btil'(v_2) \, \D v_1 \, w_2}\,.
\end{gathered}
\end{equation*}

\paragraph{Lie-commutator.}
In order to determine the first Lie-commutator 
\begin{equation*}
\big[A,B\big](v) = A'(v) \, B(v) - B'(v) \, A(v)\,,
\end{equation*}
see also~\eqref{eq:Commutator}, we employ the auxiliary relations
\begin{equation*}
\begin{gathered}
\nabla \btil(v_2) = \btil'(v_2) \, \nabla v_2\,, \\
\D \, \btil(v_2) = \btil''(v_2) \, \nabla v_2 \cdot \nabla v_2 + \btil'(v_2) \, \D v_2\,, \\
\D \, \big(\btil(v_2) \, \D v_1\big) = \D \, \btil(v_2) \, \D v_1 + 2 \, \nabla \btil(v_2) \cdot \nabla \D v_1 + \btil(v_2) \, \D^2 v_1\,, \\
\btil''(v_2) = 2 \, \b \, \d^2 \big(1 - \d v_2\big)^{-3} = 2 \, \big(\tfrac{\d}{\b}\big)^2 \, \big(\btil(v_2)\big)^3\,.
\end{gathered}
\end{equation*}
Noting that the identity $\atil'(v_2) \, \btil(v_2) - \atil(v_2) \, \btil'(v_2) = 0$ holds, a brief calculation yields  
\begin{equation*}
\begin{split}
&\big[A,B\big](v) \\
&\;\;= \ma{0}{I}{0}{\atil'(v_2) \, \D v_2 + \atil(v_2) \, \D} \, \ve{0}{\btil(v_2) \, \D v_1} 
- \ma{0}{0}{\btil(v_2) \, \D}{\btil'(v_2) \, \D v_1} \ve{v_2}{\atil(v_2) \, \D v_2} \\
&\;\;= \ve{\btil(v_2) \, \D v_1}{\big(\atil'(v_2) \, \btil(v_2) - \atil(v_2) \, \btil'(v_2)\big) \, \D v_1 \, \D v_2 + \atil(v_2) \, \D \, \big(\btil(v_2) \, \D v_1\big) - \btil(v_2) \, \D v_2} \\
&\;\;= \ve{\btil(v_2) \, \D v_1}{\atil(v_2) \, \D \, \big(\btil(v_2) \, \D v_1\big) - \btil(v_2) \, \D v_2} \\
&\;\;= \ve{\btil(v_2) \, \D v_1}{\atil(v_2) \, \big(\D \, \btil(v_2) \, \D v_1 + 2 \, \nabla \btil(v_2) \cdot \nabla \D v_1 + \btil(v_2) \, \D^2 v_1\big) - \btil(v_2) \, \D v_2}\,,
\end{split}
\end{equation*}
that is, we have 
\begin{equation*}
\begin{gathered}
\big[A,B\big](v) = \ve{\zeta_1}{\zeta_2}\,, \qquad \zeta_1 = \btil(v_2) \, \D v_1\,, \\
\zeta_2 = \atil(v_2) \, \big(\D \, \btil(v_2) \, \D v_1 + 2 \, \nabla \btil(v_2) \cdot \nabla \D v_1 + \btil(v_2) \, \D^2 v_1\big) - \btil(v_2) \, \D v_2\,.
\end{gathered}
\end{equation*}
As a consequence, by means of the relations 
\begin{equation*}
\begin{gathered}
\normbig{L^2}{\zeta_1} \leq C\big(\norm{H^2}{v_1}, \norm{H^2}{v_2}\big)\,, \\
\normbig{L^2}{\zeta_2} \leq C\big(\norm{H^4}{v_1}, \norm{H^2}{v_2}\big)\,, \\
\normbig{H^{k+2}}{\zeta_1} \leq C\big(\norm{H^{k+4}}{v_1}, \norm{H^{k+2}}{v_2}\big)\,, \qquad k \in \NN_{\geq 0}\,, \\
\normbig{H^{k+2}}{\zeta_2} \leq C\big(\norm{H^{k+6}}{v_1}, \norm{H^{k+4}}{v_2}\big)\,, \qquad k \in \NN_{\geq 0}\,,
\end{gathered}
\end{equation*}
the following estimate is obtained
\begin{equation*}
\normbig{H^{k} \times H^{k}}{\big[A,B\big](v)} + \normbig{H^{k+2} \times H^k}{\big[A,B\big](v)} \leq C\big(\norm{H^{k+4} \times H^{k+2}}{v}\big)\,, \qquad k \in \NN_{\geq 0}\,,
\end{equation*}
see also~\eqref{eq:AuxiliaryEstimates}; 
the special case $k = 1$ is included by interpolation, see~\cite[Thm~7.23]{AdamsFournier}. 
\subsection{Estimates for evolution operators}
\label{sec:EvolutionOperators}
In the following, we state estimates for the evolution operators associated with the Westervelt equation and the subproblems corresponding to Decomposition~I as well as their derivatives, see~\eqref{eq:WesterveltReformulations}--\eqref{eq:Problem} and~\eqref{eq:Subproblems}.  
For this purpose, we make use of the fact that the Fréchet derivative of the evolution operator~$\nEG(\cdot,v)$ with respect to the initial state~$v$, denoted by~$\pa_2 \nEG(\cdot,v)$, is a solution to the variational equation  
\begin{equation*}
\begin{cases}
&\tfrac{\dd}{\dd t} \, \pa_{2} \nE_G(t,v) = G'\big(\nE_G(t,v)\big) \, \pa_{2} \nE_G(t,v)\,, \quad t \in (0,T]\,, \\
&\pa_{2} \nE_G(0,v) = I\,,
\end{cases}
\end{equation*}
see also~\eqref{eq:dEG}.
For the convenience of the reader, we collect the relations 
\begin{equation*}
\begin{gathered}
F(v) = A(v) + B(v)\,, \qquad F'(v) = A'(v) + B'(v)\,, \\
A(v) = \ve{v_2}{\atil(v_2) \, \D v_2}\,, \qquad A'(v) = \ma{0}{I}{0}{\atil(v_2) \, \D + \atil'(v_2) \, \D v_2}\,, \\
B(v) = \ve{0}{\btil(v_2) \, \D v_1}\,, \qquad B'(v) = \ma{0}{0}{\btil(v_2) \, \D}{\btil'(v_2) \, \D v_1}\,, \\
\atil(v_2) = \a \, \big(1 - \d v_2\big)^{-1}\,, \qquad \atil'(v_2) = \a \, \d \, \big(1 - \d v_2\big)^{-2}\,, \\
\btil(v_2) = \b \, \big(1 - \d v_2\big)^{-1}\,, \qquad \btil'(v_2) = \b \, \d \, \big(1 - \d v_2\big)^{-2}\,, 
\end{gathered}
\end{equation*}
see also Section~\ref{sec:DerivativesCommutator}.
We point out that suitable applications of the linear variation-of-constants formula and Gronwall's lemma imply that the constants arising in the estimates for the evolution operators are of the special form~$\ee^{Ct}$, ensuring that repeated applications of the evolution operators remain bounded; 
this is crucial in view of the stability estimate~\eqref{eq:GlobalErrorStability}.  
For simplicity, we do not distinguish an abstract function $v = (v_1,v_2)$ and the associated function $v: \overline{\Om} \times [0,T] \to \RR^2: (x,t) \mapsto v(x,t) = v(t)(x)$ in notation. 

\paragraph{Estimates for evolution operators.}
Let $t \in [0,T]$.
\begin{enumerate}[(i)]
\item
\begin{enumerate}[(a)]
\item 
In the situation of Theorem~\ref{thm:TheoremRegularity} in Section \ref{sec:Regularity} the evolution operator~$\nEF$ satisfies 
\begin{equation*}
\normbig{H^{k+6} \times H^{k+5}}{\nEF(t,v)} \leq \ee^{Ct} \, \norm{H^{k+6} \times H^{k+5}}{v}\,, \qquad k \in \NN_{\geq 0}\,.
\end{equation*}
\item
For the associated derivative with respect to the initial state the estimate 
\begin{equation*}
\normbig{H^{\ell+1} \times H^{\ell}}{\pa_{2} \nEF(t,v) \, w} \leq \ee^{C(\norm{H^{4} \times H^{4}}{v}) t} \, \norm{H^{\ell+1} \times H^{\ell}}{w}\,, \qquad \ell = 0,1,2,3\,,
\end{equation*}
is valid. 
\end{enumerate}
\paragraph{Proof.}
(a) \, 
The stated estimate follows from Theorem~\ref{thm:TheoremRegularity} and interpolation~\cite[Thm~7.23]{AdamsFournier}. \\
(b) \, 
We denote by $v(t) = \nEF(t,v(0))$ the evolution operator associated with~$F$ and by $V(t) = \pa_{2} \nEF(t,v(0))$ the corresponding Fréchet dervative with respect to the initial state, which satisfies the initial value problem  
\begin{equation*}
\begin{cases}
&\ddt \, V(t) = F'\big(v(t)\big) \, V(t) \,, \quad t \in (0,T]\,, \\
&V(0) = I\,,
\end{cases} 
\end{equation*}
see~\eqref{eq:dEG};
due to linearity, the solution to the variational equation is given by a relation of the form 
\begin{equation*}
\pa_{2} \nEF\big(t,v(0)\big) \, w = V(t) \, w = \ee^{D(t,0)} \, w\,, \qquad t \in (0,T]\,.
\end{equation*}
More precisely, using that 
\begin{equation*}
\begin{gathered}
V = \ma{V_{11}}{V_{12}}{V_{21}}{V_{22}}\,, \qquad F' = \ma{0}{I}{F'_{21}}{F'_{22}}\,, \\
F'_{21}(v) = \btil(v_2) \, \D\,, \qquad F'_{22}(v) = \atil(v_2) \, \D + \atil'(v_2) \, \D v_2 + \btil'(v_2) \, \D v_1\,, \\
\end{gathered}
\end{equation*}
we obtain the two decoupled systems  
\begin{equation*}
\begin{split}
&\begin{cases}
&\ddt \, V_{11}(t) = V_{21}(t)\,, \quad t \in (0,T]\,, \\
&\ddt \, V_{21}(t) = F'_{21}\big(v(t)\big) \, V_{11}(t) + F'_{22}\big(v(t)\big) \, V_{21}(t)\,, \quad t \in (0,T]\,, \\
&V_{11}(0) = I\,, \qquad V_{21}(0) = 0\,, 
\end{cases} \\
&\begin{cases}
&\ddt \, V_{12}(t) = V_{22}(t)\,, \quad t \in (0,T]\,, \\
&\ddt \, V_{22}(t) = F'_{21}\big(v(t)\big) \, V_{12}(t) + F'_{22}\big(v(t)\big) \, V_{22}(t)\,, \quad t \in (0,T]\,, \\
&V_{12}(0) = 0\,, \qquad V_{22}(0) = I\,,
\end{cases}
\end{split}
\end{equation*}
which correspond to the strongly damped linear wave equations 
\begin{equation*}
\begin{split}
&\begin{cases}
&\ddtwo \, V_{11}(t) = F'_{21}\big(v(t)\big) \, V_{11}(t) + F'_{22}\big(v(t)\big) \, \ddt \, V_{11}(t)\,, \quad t \in (0,T]\,, \\
&V_{11}(0) = I\,, \quad \ddt \, V_{11}(0) = V_{21}(0) = 0\,, \\
\end{cases} \\
&\begin{cases}
&\ddtwo \, V_{12}(t) = F'_{21}\big(v(t)\big) \, V_{12}(t) + F'_{22}\big(v(t)\big) \, \ddt \, V_{12}(t)\,, \quad t \in (0,T]\,, \\
&V_{12}(0) = 0\,, \quad \ddt \, V_{12}(0) = V_{22}(0) = I\,.
\end{cases}
\end{split}
\end{equation*}
Both equations can be cast into the form~\eqref{eq:LinearHyperbolic} with coefficient functions 
\begin{equation*}
a = \atil(v_2)\,, \qquad b = \atil'(v_2) \, \D v_2 + \btil'(v_2) \, \D v_1\,, \qquad c = \btil(v_2)\,, \qquad f = 0\,.
\end{equation*}
We note that the assumptions of Proposition~\ref{prop:PropositionHyperbolic} are verified as follows:
From Theorem~\ref{thm:TheoremRegularity} (i) we have
\begin{equation*}
\begin{gathered}
0 < \underline{\nu} \leq 1 - \d \, v_2(t) \leq \overline{\nu} < \infty\,, \quad t \in [0,T]\,, \\
(v_1,v_2) \in \nC\big([0,T], H^4(\Om) \times H^4(\Om)\big)\,, \quad \pat v_2\in \nC\big([0,T], H^2(\Om)\big)\,, \\
\normbig{H^4 \times H^4\times H^2}{\big(v_1(t),v_2(t),\pat v_2(t)\big)} \leq \ee^{C t} \, \normbig{H^4 \times H^4}{\big(v_1(0),v_2(0)\big)}\,, \quad t \in [0,T]\,, 
\end{gathered}
\end{equation*}
where the first relation is understod pointwise for $x \in \overline{\Om}$; 
these estimates can be made use of in 
\begin{equation*} 
\begin{gathered}
a(t) \geq \a \, \underline{\nu} > 0\,, \quad c(t) \geq \b \, \underline{\nu} > 0\,, \qquad t \in [0,T]\,, \\
\normbig{L^2((0,T),L^3)}{b} \leq \sqrt{T} \, \normbig{\nC((0,T),H^1)}{b} \leq \sqrt{T} \, C\big(\normbig{\nC((0,T),H^3 \times H^3)}{(v_1,v_2)}\big)\,, \\
\normbig{L^2((0,T),L^\infty)}{\pat c} \leq \sqrt{T} \, \normbig{\nC([0,T],L^\infty)}{\pat c} \leq \sqrt{T} \, C\big(\normbig{\nC([0,T],H^2)}{\pat v_2}\big)\,.
\end{gathered}
\end{equation*}
Thus, Proposition~\ref{prop:PropositionHyperbolic} implies 
\begin{equation*}
\begin{split}
\normbig{H^{\ell+1} \times H^{\ell}}{V(t) \, w} 
&\leq \ee^{C(\norm{\nC([0,t],H^{4} \times H^{4})}{v} + \norm{\nC([0,t],H^{2})}{\pat v_2})} \, \norm{H^{\ell+1} \times H^{\ell}}{w} \\
&\leq \ee^{C(\norm{H^{4} \times H^{4}}{v(0)}) t} \, \norm{H^{\ell+1} \times H^{\ell}}{w}\,, \qquad \ell = 0,1,2,3\,, 
\end{split}
\end{equation*}
which proves the stated result. 
$\hfill \diamondsuit$
\item 
\begin{enumerate}[(a)]
\item 
Arguments based on the regularity result~\cite[Thm~6, p.\,386]{Evans} for linear evolution equations of parabolic type imply 
\begin{equation*}
\normbig{H^{k+4} \times H^{k+2}}{\nEA(t,v)} \leq \ee^{Ct} \, \norm{H^{k+4} \times H^{k+2}}{v}\,, \qquad k \in \NN_{\geq 0}\,. 
\end{equation*}
\item 
For any $k \in \NN_{\geq 0}$ the estimate 
\begin{equation*}
\normbig{H^k \times H^{\ell}}{\pa_2 \nEA(t,v) \, w} 
\leq \begin{cases}
\ee^{C(\norm{H^5 \times H^3}{v}) t} \, \norm{H^k \times H^{\ell}}{w}\,, \quad \ell = 0,1,2\,, \\
\ee^{C(\norm{H^7 \times H^5}{v}) t} \, \norm{H^k \times H^{\ell}}{w}\,, \quad \ell \in \NN_{\geq 3}\,, 
\end{cases}
\end{equation*}
is valid. 
\end{enumerate}
\paragraph{Proof.}
(a) \, 
For the following, we set $v(t) = \nEA(t,v(0))$.
The second component of the solution, governed by the nonlinear diffusion equation
\begin{equation*}
\pat v_2(t) = \a \, \big(1 - \d v_2(t)\big)^{-1} \, \D v_2(t)\,, 
\end{equation*}
inherites the regularity of the prescribed initial state;
more precisely, an application of Proposition~\ref{prop:PropositionParabolic} with 
\begin{equation*}
a = \a \, \big(1 - \d v_2\big)^{-1}\,, \qquad b = 0\,, \qquad f=0\,, 
\end{equation*}
see also~\cite[Thm~6, p.\,386]{Evans}, a fixed point argument, and Gronwall's lemma imply 
\begin{equation*}
\normbig{H^{k+2}}{v_2(t)} \leq \ee^{Ct} \, \normbig{H^{k+2}}{v_2(0)}\,, \qquad k \in \NN_{\geq 0}\,.
\end{equation*}
The regularity of the second component carries over to the first component 
\begin{equation*}
v_1(t) = v_1(0) + \int_{0}^{t} v_2(\t) \; \dd\t\,, 
\end{equation*}
and we thus obtain the estimate 
\begin{equation*}
\normbig{H^{k+2} \times H^{k+2}}{v(t)} \leq \ee^{Ct} \, \normbig{H^{k+2} \times H^{k+2}}{v(0)}\,, \qquad k \in \NN_{\geq 0}\,.
\end{equation*}
Furthermore, employing the identity 
\begin{equation*}
\pat \, \Big(- \, \tfrac{1}{2 \d} \, \big(1 - \d v_2(t)\big)^2\Big) = \big(1 - \d v_2(t)\big) \, \pat v_2(t) = \a \, \D v_2(t) 
\end{equation*}
and performing integration 
\begin{equation*}
\begin{split}
\tfrac{1}{2 \d} \, \Big(\big(1 - \d v_2(0)\big)^2 - \, \big(1 - \d v_2(t)\big)^2\Big) &= \a \, \D \int_{0}^{t} v_2(\t) \; \dd\t \\
&= \a \, \D \, \big(v_1(t) - v_1(0)\big)\,, 
\end{split}
\end{equation*}
leads to the relation 
\begin{equation*}
\D \, v_1(t) = \D \, v_1(0) + \tfrac{1}{2 \a \d} \, \Big(\big(1 - \d v_2(0)\big)^2 - \, \big(1 - \d v_2(t)\big)^2\Big)\,.
\end{equation*}
Making use of the fact that~$H^{k+2}(\Om)$ forms an algebra for any $k \in \NN_{\geq 0}$ leads to 
\begin{equation*}
\normbig{H^{k+4} \times H^{k+2}}{v(t)} \leq \ee^{Ct} \, \normbig{H^{k+4} \times H^{k+2}}{v(0)}\,, \qquad k \in \NN_{\geq 0}\,, 
\end{equation*}
and thus proves the statement. \\
(b) \, 
Analogously to before, we denote $v(t) = \nEA(t,v(0))$ and consider the initial value problem 
\begin{equation*}
\begin{cases}
&\ddt \, V(t) = A'\big(v(t)\big) \, V(t) \,, \quad t \in (0,T]\,, \\
&V(0) = I\,,
\end{cases} 
\end{equation*}
for $V(t) = \pa_{2} \nEA(t,v(0))$. 
A reformulation in components 
\begin{equation*}
V = \ma{V_{11}}{V_{12}}{V_{21}}{V_{22}}\,, \qquad A' = \ma{0}{I}{0}{A'_{22}}\,, 
\qquad A'_{22}(u) = \atil(v_2) \, \D + \atil'(v_2) \, \D v_2\,,
\end{equation*}
leads to two decoupled systems involving a linear diffusion-reaction equation with time-dependent coefficients 
\begin{equation*}
\begin{split}
&\begin{cases}
&\ddt \, V_{11}(t) = V_{21}(t)\,, \quad t \in (0,T]\,, \qquad V_{11}(0) = I\,, \\
&\ddt \, V_{21}(t) = A'_{22}\big(u(t)\big) \, V_{21}(t)\,, \quad t \in (0,T]\,, \qquad V_{21}(0) = 0\,, \\
\end{cases} \\
&\begin{cases}
&\ddt \, V_{12}(t) = V_{22}(t)\,, \quad t \in (0,T]\,, \qquad V_{12}(0) = 0\,, \\
&\ddt \, V_{22}(t) = A'_{22}\big(u(t)\big) \, V_{22}(t)\,, \quad t \in (0,T]\,, \qquad V_{22}(0) = I\,.
\end{cases}
\end{split}
\end{equation*}
An application of~Proposition \ref{prop:PropositionParabolic} with 
\begin{equation*}
a = \atil(v_2)\,, \qquad b = \atil'(v_2) \, \D v_2\,, \qquad f = 0\,, 
\end{equation*}
leads to the following estimates 
\begin{equation*}
\begin{gathered}
\norm{H^k \times H^{\ell}}{V(t) \, w} \leq \ee^{C(\norm{\nC([0,t],H^3)}{v_2})} \, \norm{H^k \times H^{\ell}}{w}\,, \qquad k \in \NN_{\geq 0}\,, \quad \ell = 0,1,2\,, \\
\norm{H^k \times H^{\ell}}{V(t) \, w} \leq \ee^{C(\norm{\nC([0,t],H^5)}{v_2}} \, \norm{H^k \times H^{\ell}}{w}\,, \qquad k \in \NN_{\geq 0}\,, \quad \ell \in \NN_{\geq 3}\,,
\end{gathered}
\end{equation*}
which imply the stated result. 
$\hfill \diamondsuit$
\item 
\begin{enumerate}[(a)]
\item 
The evolution operator associated with~$B$ fulfills 
\begin{equation*}
\normbig{H^{k+2} \times H^{k}}{\nEB(t,v)} \leq \ee^{Ct} \, \norm{H^{k+2} \times H^{k}}{v}\,, \qquad k \in \NN_{\geq 0}\,.
\end{equation*}
\item 
The estimate 
\begin{equation*}
\normbig{H^{k+2} \times H^k}{\pa_2 \nEB(t,v) \, w} \leq \ee^{C(\norm{H^{k+4} \times H^{k+2}}{v}) t} \, \norm{H^{k+2} \times H^{k}}{w}\,, \qquad k \in \NN_{\geq 0}\,,
\end{equation*}
is valid.
\end{enumerate}
\paragraph{Proof.}
(a) \, 
Let $v(t) = \nEB(t,v(0))$. 
The first component of the solution remains constant in time 
\begin{equation*}
v_1(t) = v_1(0)\,;
\end{equation*}
using that the second component is given by the explicit representation 
\begin{equation*}
v_2(t) = \tfrac{1}{\d} \, \Big(1 - \sqrt{\big(1 - \d v_2(0)\big)^2 - 2 \, \b \d \, t \, \D v_1(0)}\;\Big)\,,
\end{equation*}
the stated result follows. \\
(b) \, 
Differentiation of the explicit solution representation 
\begin{equation*}
v(t) = \ve{v_{1}(0)}{\tfrac{1}{\d} \, \Big(1 - \sqrt{(1 - \d v_{2}(0))^2 - 2 \, \b \d \, t \, \D v_{1}(0)}\;\Big)} 
\end{equation*}
with respect to the initial state yields 
\begin{equation*}
\begin{split}
V(t) \, w 
&= \ma{I}{0}{\frac{\b \, t}{\sqrt{(1 - \d v_{2}(0))^2 - 2 \, \b \d \, t \, \D v_{1}(0)}} \, \D}{\frac{1 - \d \, v_{2}(0)}{
\sqrt{(1 - \d v_{2}(0))^2 - 2 \, \b \d \, t \, \D v_{1}(0)}}} \, \ve{w_1}{w_2} \\
&= \ve{w_1}{\frac{1}{\sqrt{(1 - \d v_{2}(0))^2 - 2 \, \b \d \, t \, \D v_{1}(0)}} \, \Big(\b \, t \, \D w_1 + \big(1 - \d \, v_{2}(0)\big) \, w_2\Big)}\,, 
\end{split}
\end{equation*}
where $V(t) = \pa_2 \nEB(t,v(0))$.
Suitable estimation proves the statement. 
$\hfill \diamondsuit$
\end{enumerate}
\subsection{Main result}
In the following, we deduce a global error estimate for the Lie--Trotter splitting method applied to the Westervelt equation. 
Under the assumption that the prescribed initial state fulfills the regularity requirements 
\begin{equation}
\begin{gathered}
\label{eq:ConditionPs0}
u(0) = \big(\ps(\cdot,0), \pat \ps(\cdot,0)\big) \in H^6(\Om) \times H^5(\Om)\,, \\
\normbig{H^6 \times H^5}{u(0)} = \normbig{H^6}{\ps(\cdot,0)} + \normbig{H^5}{\pat \ps(\cdot,0)} \leq C_0\,,  
\end{gathered}
\end{equation}
with suitably chosen constant $C_0 > 0$ as well as certain compatibility conditions, Theorem~\ref{thm:TheoremRegularity} given in Section~\ref{sec:Regularity} guarantees that the solution to the Westervelt equation~\eqref{eq:WesterveltReformulations}--\eqref{eq:Problem} satisfies 
\begin{equation*}
u(t) = \big(\ps(\cdot,t), \pat \ps(\cdot,t)\big) \in H^6(\Om) \times H^5(\Om)\,, \qquad t \in [0,T]\,, 
\end{equation*}
and remains bounded. 
As discussed in Section~\ref{sec:Convergence}, the extension of the global error estimate to variable time stepsizes is straightforward. 

\begin{theorem}[Lie--Trotter splitting method, Decomposition~I]
\label{thm:MainResult}
Assume that the initial state fulfills the regularity requirement~\eqref{eq:ConditionPs0} and that the intial approximation $u_0 \approx u(0)$ remains bounded in $H^{5}(\Om) \times H^3(\Om)$.  
Then the Lie--Trotter splitting method applied to the Westervelt equation~\eqref{eq:WesterveltReformulations}--\eqref{eq:Problem} satisfies the global error estimate 
\begin{equation*}
\normbig{H^3 \times H^1}{u_N - u(t_N)} \leq C \, \Big(\normbig{H^3 \times H^1}{u_0 - u(0)} + h\Big)\,, \qquad 0 \leq t_N = N h \leq T\,, 
\end{equation*}
with constant depending on bounds for $\norm{\nC([0,t_N],H^6 \times H^5)}{u}$ as well as $\norm{H^5 \times H^3}{u_0}$ and the final time~$t_N$.
\end{theorem}
\paragraph{Proof.} 
Let $t \in [0,T]$ and $0 \leq t_N = N h \leq T$.
In order to deduce the stated global error estimate for the Lie--Trotter splitting method based on Decomposition~I, we follow the approach described in Section~\ref{sec:Convergence}.
For convenience, we collect the fundamental auxiliary results deduced in Sections~\ref{sec:DerivativesCommutator}--\ref{sec:EvolutionOperators}, adapted to the present situation 
\begin{subequations}
\begin{gather}
\label{eq:EstimateEF}
\normbig{H^{k+6} \times H^{k+5}}{\nEF(t,v)} \leq \ee^{Ct} \, \norm{H^{k+6} \times H^{k+5}}{v}\,, \qquad k \in \NN_{\geq 0}\,, \\
\label{eq:EstimateEA}
\normbig{H^{k+4} \times H^{k+2}}{\nEA(t,v)} \leq \ee^{Ct} \, \norm{H^{k+4} \times H^{k+2}}{v}\,, \qquad k \in \NN_{\geq 0}\,, \\
\label{eq:EstimateEB}
\normbig{H^{k+2} \times H^{k}}{\nEB(t,v)} \leq \ee^{Ct} \, \norm{H^{k+2} \times H^{k}}{v}\,, \qquad k \in \NN_{\geq 0}\,, \\
\label{eq:EstimatedEF}
\normbig{H^{\ell+1} \times H^{\ell}}{\pa_{2} \nEF(t,v) \, w} \leq \ee^{C(\norm{H^{4} \times H^{4}}{v}) t} \, \norm{H^{\ell+1} \times H^{\ell}}{w}\,, \qquad \ell = 0,1,2,3\,, \\
\label{eq:EstimatedEA}
\normbig{H^{k+2} \times H^k}{\pa_2 \nEA(t,v) \, w} \leq 
\begin{cases}
\ee^{C(\norm{H^5 \times H^3}{v}) t} \, \norm{H^{k+2} \times H^k}{w}\,, \qquad k = 0,1,2\,, \\
\ee^{C(\norm{H^7 \times H^5}{v}) t} \, \norm{H^{k+2} \times H^k}{w}\,, \qquad k \in \NN_{\geq 3}\,, \\
\end{cases} \\
\label{eq:EstimatedEB}
\normbig{H^{k+2} \times H^k}{\pa_2 \nEB(t,v) \, w} \leq \ee^{C(\norm{H^{k+4} \times H^{k+2}}{v}) t} \, \norm{H^{k+2} \times H^{k}}{w}\,, \qquad k \in \NN_{\geq 0}\,, \\
\label{eq:EstimateAB}
\normbig{H^{k+2} \times H^k}{\big[A,B\big](v)} \leq C\big(\norm{H^{k+4} \times H^{k+2}}{v}\big)\,, \qquad k \in \NN_{\geq 0}\,.
\end{gather}
\end{subequations}
\begin{enumerate}[(i)]
\item 
The basic regularity assumptions on the initial state and the initial numerical approximation
\begin{equation*}
\begin{gathered}
u(0) \in H^6(\Om) \times H^5(\Om)\,, \qquad u_0 \in H^5(\Om) \times H^3(\Om)\,,
\end{gathered}
\end{equation*}
ensure that the exact and numerical solution values, given by 
\begin{equation*}
\begin{gathered}
u(t) = \nEF\big(t,u(0)\big)\,, \\ 
u_N = \nSF^N(h,u_0)\,, \qquad \nSF(h,v) = \nEB\big(h,\nEA(h,v)\big)\,, 
\end{gathered}
\end{equation*}
remain in the underlying product spaces; 
indeed, due to~\eqref{eq:EstimateEF} with $k = 0$ as well as~\eqref{eq:EstimateEA}--\eqref{eq:EstimateEB} with $k = 1$ and $k = 3$, respectively, we obtain 
\begin{equation}
\label{eq:EstimateEFSFn}
\begin{gathered}
\normbig{H^6 \times H^5}{\nEF\big(t,u(0)\big)} \leq \ee^{Ct} \, \normbig{H^6 \times H^5}{u(0)}\,, \\
\normbig{H^5 \times H^3}{\nSF^n(h,u_0)} \leq \ee^{C t_n} \, \norm{H^5 \times H^3}{u_0}\,, 
\end{gathered}
\end{equation}
for integer $1 \leq n \leq N$.
We note that the constant arising in the estimation of the splitting operator  
\begin{equation*}
\begin{split}
\normbig{H^5 \times H^3}{\nSF(h,v)} 
&= \normbig{H^5 \times H^3}{\nEB\big(h,\nEA(h,v)\big)} \leq \ee^{Ch} \, \normbig{H^5 \times H^3}{\nEA(h,v)} \\
&\leq \ee^{Ch} \, \norm{H^5 \times H^3}{v}
\end{split}
\end{equation*}
is of the special form~$\ee^{Ch}$ and thus repeated applications of the splitting method lead to the stated bound. 
\item 
In order to establish a stability estimate~\eqref{eq:GlobalErrorStability}, we employ the expansion 
\begin{equation*}
\begin{split}
\nEG(h,v) - \nEG(h,\vtil) 
&= \nEG\big(h,\s v + (1-\s) \vtil\big)\Big\vert_{\s = 0}^{1} \\
&= \int_{0}^{1} \pa_2 \nEG\big(h,\s v + (1-\s) \vtil\big) \, (v - \vtil) \; \dd\s\,, 
\end{split}
\end{equation*}
which implies the identity 
\begin{equation*}
\begin{split}
\nSF(h,v) - \nSF(h,\vtil)
&= \int_{0}^{1} \pa_2 \nEB\big(h,\s_2 \, \nEA(t,v) + (1-\s_2) \, \nEA(t,\vtil)\big) \; \dd\s_2 \\
&\qquad \times \int_{0}^{1} \pa_2 \nEA\big(h,\s_1 v + (1-\s_1) \vtil\big) \; \dd\s_1 \; (v - \vtil)\,. 
\end{split}
\end{equation*}
As a consequence, by~\eqref{eq:EstimatedEB}, a further application of~\eqref{eq:EstimateEA}, and~\eqref{eq:EstimatedEA} the estimates 
\begin{equation*}
\begin{split}
k = 0,1,2: \qquad 
&\normbig{H^{k+2} \times H^k}{\nSF(h,v) - \nSF(h,\vtil)} \\
&\quad\leq \ee^{C(\norm{H^5 \times H^3}{v}, \norm{H^{k+4} \times H^{k+2}}{v}, \norm{H^5 \times H^3}{\vtil}, \norm{H^{k+4} \times H^{k+2}}{\vtil}) h} \\ &\quad\qquad\times \norm{H^{k+2} \times H^{k}}{v - \vtil}\,, \\
k \geq 3: \qquad 
&\normbig{H^{k+2} \times H^k}{\nSF(h,v) - \nSF(h,\vtil)} \\
&\quad\leq \ee^{C(\norm{H^7 \times H^5}{v}, \norm{H^{k+4} \times H^{k+2}}{v}, \norm{H^7 \times H^5}{\vtil}, \norm{H^{k+4} \times H^{k+2}}{\vtil}) h} \\ 
&\quad\qquad\times \norm{H^{k+2} \times H^{k}}{v - \vtil}\,,
\end{split}
\end{equation*}
follow.
This suggest the choice $k = 1$ and implies the stability bound
\begin{equation*}
\begin{split}
&\normbig{H^3 \times H^1}{\nSF^n(h,v) - \nSF^n(h,\vtil)} 
\leq \ee^{C(\norm{H^5 \times H^3}{v}, \norm{H^5 \times H^3}{\vtil}) t_n} \, \norm{H^3 \times H^1}{v - \vtil}\,.
\end{split}
\end{equation*}
\item 
Applying the above considerations to the global error representation~\eqref{eq:GlobalError} yields 
\begin{equation*}
\begin{split}
\normbig{H^3 \times H^1}{u_N - u(t_N)} 
&\leq \normbig{H^3 \times H^1}{\nSF^N\big(h,u_0\big) - \nSF^N\big(h,u(0)\big)} \\
&\qquad+ \sum_{n=0}^{N-1} \Big\|\nSF^{N-n-1}\Big(h,\nLF\big(h,u(t_n)\big) + \nEF\big(h,u(t_n)\big)\Big) \\
&\qquad\qquad\qquad- \nSF^{N-n-1}\Big(h,\nEF\big(h,u(t_n)\big)\Big)\Big\|_{H^3 \times H^1} \\
&\leq C \, \bigg(\normbig{H^3 \times H^1}{u_0 - u(0)} + \sum_{n=0}^{N-1} \normbig{H^3 \times H^1}{\nLF\big(h,u(t_n)\big)}\bigg)
\end{split}
\end{equation*}
with constant depending in particular on bounds for $\norm{H^5 \times H^3}{u(0)}$ and $\norm{H^5 \times H^3}{u_0}$.  
\item 
It remains to estimate the local error
\begin{equation*}
\begin{gathered}
\nLF\big(h,u(t_n)\big) 
= \int_{0}^{h} \int_{0}^{\t_1} \pa_{2} \nEF\big(h-\t_1,\nSF\big(\t_1,u(t_n)\big)\big) \, Z_1(\t_1,\t_2) \; \dd\t_2 \, \dd\t_1\,, \\ 
Z_1(\t_1,\t_2) = \pa_{2} \nEB(\t_1,w) \, \big(\pa_{2} \nEB(\t_2,w)\big)^{-1} \, Z_2(\t_1,\t_2)\,, \\
Z_2(\t_1,\t_2) = \big[B,A\big]\big(\nEB(\t_2,w)\big)\,, \qquad w = \nEA\big(\t_1,u(t_n)\big)\,,
\end{gathered}
\end{equation*}
see~\eqref{eq:L}. 
For this purpose, we first apply~\eqref{eq:EstimatedEF} with $\ell = 2$
\begin{equation*}
\begin{split}
&\normbig{H^3 \times H^1}{\nLF\big(h,u(t_n)\big)} \leq \normbig{H^3 \times H^2}{\nLF\big(h,u(t_n)\big)} \\
&\qquad\leq C \, \int_{0}^{h} \int_{0}^{\t_1} \normbig{H^3 \times H^2}{\pa_{2} \nEF\big(h-\t_1,\nSF\big(\t_1,u(t_n)\big)\big) \, Z_1(\t_1,\t_2)} \\
&\qquad\leq C \, \int_{0}^{h} \int_{0}^{\t_1} \normbig{H^3 \times H^2}{Z_1(\t_1,\t_2)} \\
&\qquad\leq C \, \int_{0}^{h} \int_{0}^{\t_1} \normbig{H^4 \times H^2}{Z_1(\t_1,\t_2)}\,, 
\end{split}
\end{equation*}
with constant depending in particular on a bound for 
\begin{equation*}
\normbig{H^4 \times H^4}{\nSF\big(h,u(t_n)\big)} \leq \normbig{H^6 \times H^4}{\nSF\big(h,u(t_n)\big)}\,, 
\end{equation*}
which by~\eqref{eq:EstimateEA}--\eqref{eq:EstimateEB} reduces to a bound for 
\begin{equation*}
\normbig{H^6 \times H^4}{u(t_n)} \leq \normbig{H^6 \times H^5}{u(t_n)}
\end{equation*}
and further by~\eqref{eq:EstimateEF} to a bound for~$\norm{H^6 \times H^5}{u(0)}$; 
recall also the dependence on a bound for $\norm{H^5 \times H^3}{u_0}$, cf.~\eqref{eq:EstimateEFSFn}.  
Noting that the identity 
\begin{equation*}
\big(\pa_{2} \nEB(\t_2,w)\big)^{-1} = \pa_{2} \nEB\big(- \t_2,\nEB(\t_2,w)\big)
\end{equation*}
holds, see~\eqref{eq:Inverse}, we next apply~\eqref{eq:EstimatedEB} with $k = 2$ twice 
\begin{equation*}
\begin{split}
&\normbig{H^4 \times H^2}{\nLF\big(h,u(t_n)\big)} \\
&\qquad\leq C \, \int_{0}^{h} \int_{0}^{\t_1} \normbig{H^4 \times H^2}{\pa_{2} \nEB(\t_1,w) \, \pa_{2} \nEB\big(- \t_2,\nEB(\t_2,w)\big) \, Z_2(\t_1,\t_2)} \\
&\qquad\leq C \, \int_{0}^{h} \int_{0}^{\t_1} \normbig{H^4 \times H^2}{Z_2(\t_1,\t_2)}
\end{split}
\end{equation*}
and finally~\eqref{eq:EstimateAB} with $k = 2$ 
\begin{equation*}
\normbig{H^4 \times H^2}{\nLF\big(h,u(t_n)\big)} 
\leq C \, \int_{0}^{h} \int_{0}^{\t_1} \normbig{H^4 \times H^2}{\big[B,A\big]\big(\nEB(\t_2,w)\big)} \leq C \, h^2\,; 
\end{equation*}
the arising constant in addition depends on
\begin{equation*}
\begin{gathered}
\normbig{H^6 \times H^4}{w} \leq C \, \normbig{H^6 \times H^4}{u(t_n)} \leq C \, \normbig{H^6 \times H^5}{u(0)}\,, \\ 
\normbig{H^6 \times H^4}{\nEB(\t,w)} \leq C \, \normbig{H^6 \times H^4}{w} \leq C \, \normbig{H^6 \times H^5}{u(0)}\,, \qquad \t \in [0,h]\,.
\end{gathered}
\end{equation*}
Altogether this proves the stated global error estimate 
\begin{equation*}
\normbig{H^3 \times H^1}{u_N - u(t_N)} \leq C \, \Big(\normbig{H^3 \times H^1}{u_0 - u(0)} + h\Big)
\end{equation*}
with constant depending on bounds for~$\norm{H^6 \times H^5}{u(0)}$, $\norm{H^5 \times H^3}{u_0}$, and the final time.  
We note that exchanging the roles of~$A$ and~$B$, which corresponds to the application of the second scheme in~\eqref{eq:Lie}, leads to the same estimate. 
$\hfill \diamondsuit$
\end{enumerate}
\section{Numerical examples}
\label{sec:Numerics}
In this section, we include a numerical example comparing the accuracy of the Lie--Trotter and Strang splitting methods for the time integration of the Westervelt equation, based on the four decompositions that were introduced in Section~\ref{sec:Decompositions};  
furthermore, we illustrate the numerical solution obtained for a problem with more realistic parameter values. 
As our focus is on the time integration and in order to facilitate the numerical computations, we restrict ourselves to the Westervelt equation in a single space dimension;  
in particular, the spatial grid width is chosen sufficiently fine such that the global error is dominated by the time discretisation error. 
For the numerical solution of the subproblems we apply explicit and implicit time integration methods of the same order as the underlying splitting method; 
in combination with the first-order Lie--Trotter splitting method we use the explicit and implicit Euler methods, and in combination with the second-order Strang splitting method we use a second-order explicit Runge--Kutta method and the Crank--Nicolson scheme.
As expected, if an explicit solver is used for the numerical solution of the subproblems, sufficiently small time stepsizes are required to avoid instabilities; 
for a higher number of space grid points this unstable behaviour will change for the worse.   
We point out that for the chosen problem data a regular solution to the Westervelt equation exists such that no order reductions and thus no loss of accuracy is encountered when applying the Lie--Trotter and Strang splitting methods. 
We report that the application of the linearly implicit and semi-implicit Euler methods leads to numerical results of essentially the same accuracy as the implicit Euler method, with less computational effort.
For the space discretisation of the first model problem we use the Finite Difference Method with equidistant grid points, which is simple to implement, such that it is rendered possible for the reader to reproduce the numerical results; 
the same qualitative behaviour is expected for a space discretisation by the Finite Element Method, which will be the method of choice for practically relevant problems in two and three space dimensions. 

\begin{figure}[t!]
\begin{center}
\includegraphics[width=0.49\textwidth]{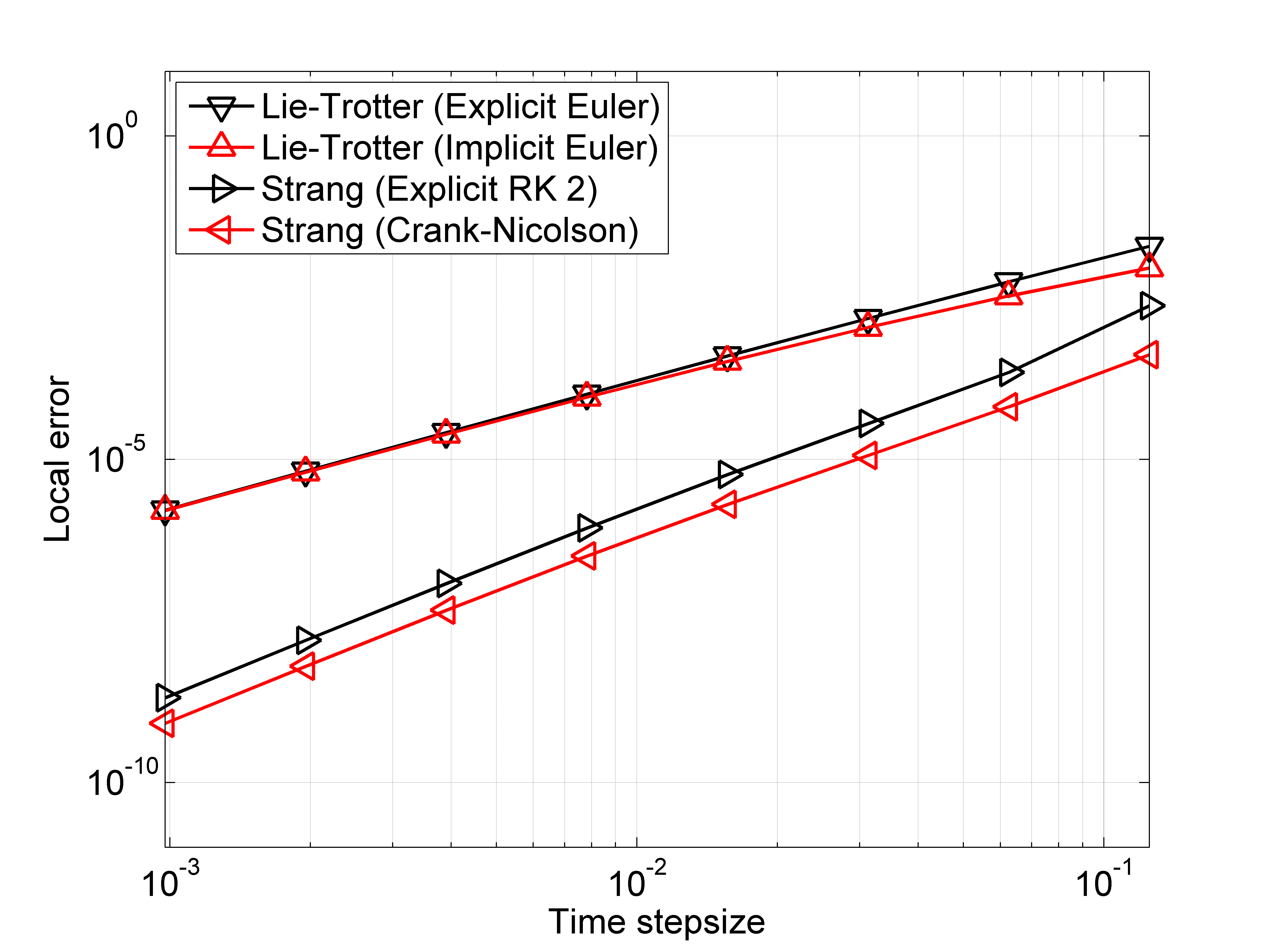} \includegraphics[width=0.49\textwidth]{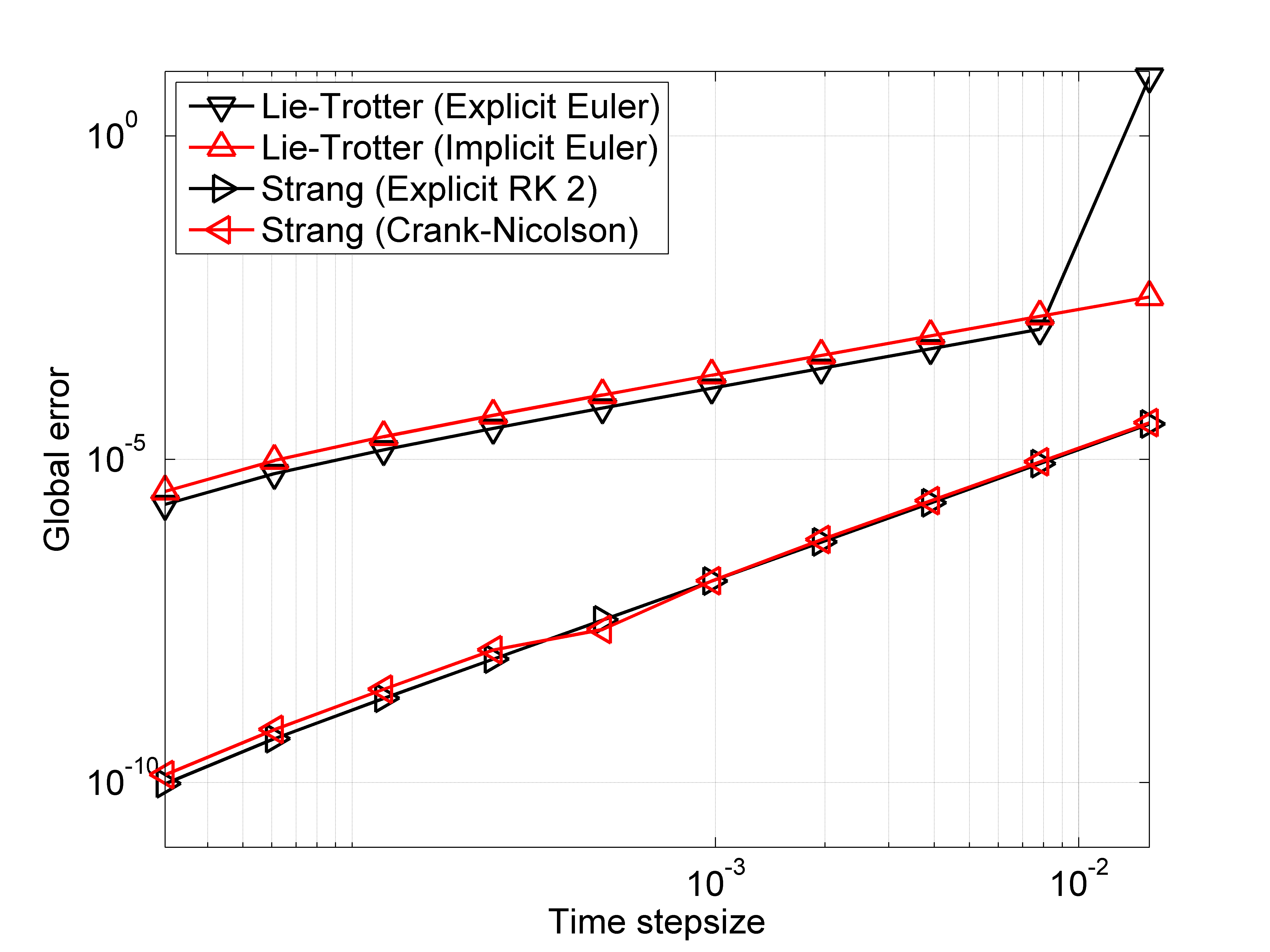} \\
\end{center}
\caption{Local (left) and global (right) errors for the Lie--Trotter and Strang splitting methods with respect to the $L^2 \times L^2$-norm obtained for Decomposition~I. 
Comparison of different time integration methods for the numerical solution of the subproblems.}
\label{fig:Figure1a}
\end{figure}

\begin{figure}[t!]
\begin{center}
\includegraphics[width=0.49\textwidth]{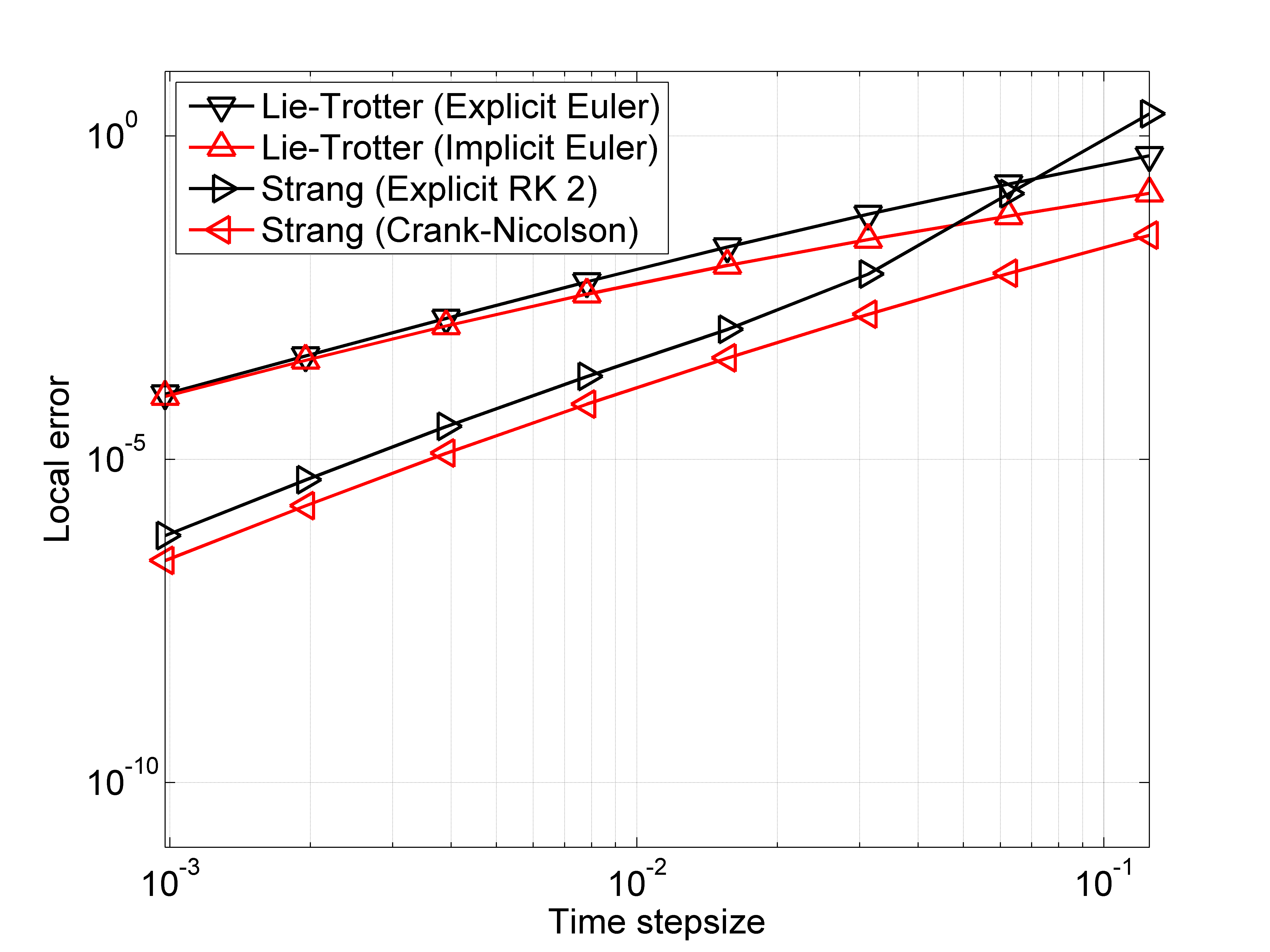} \includegraphics[width=0.49\textwidth]{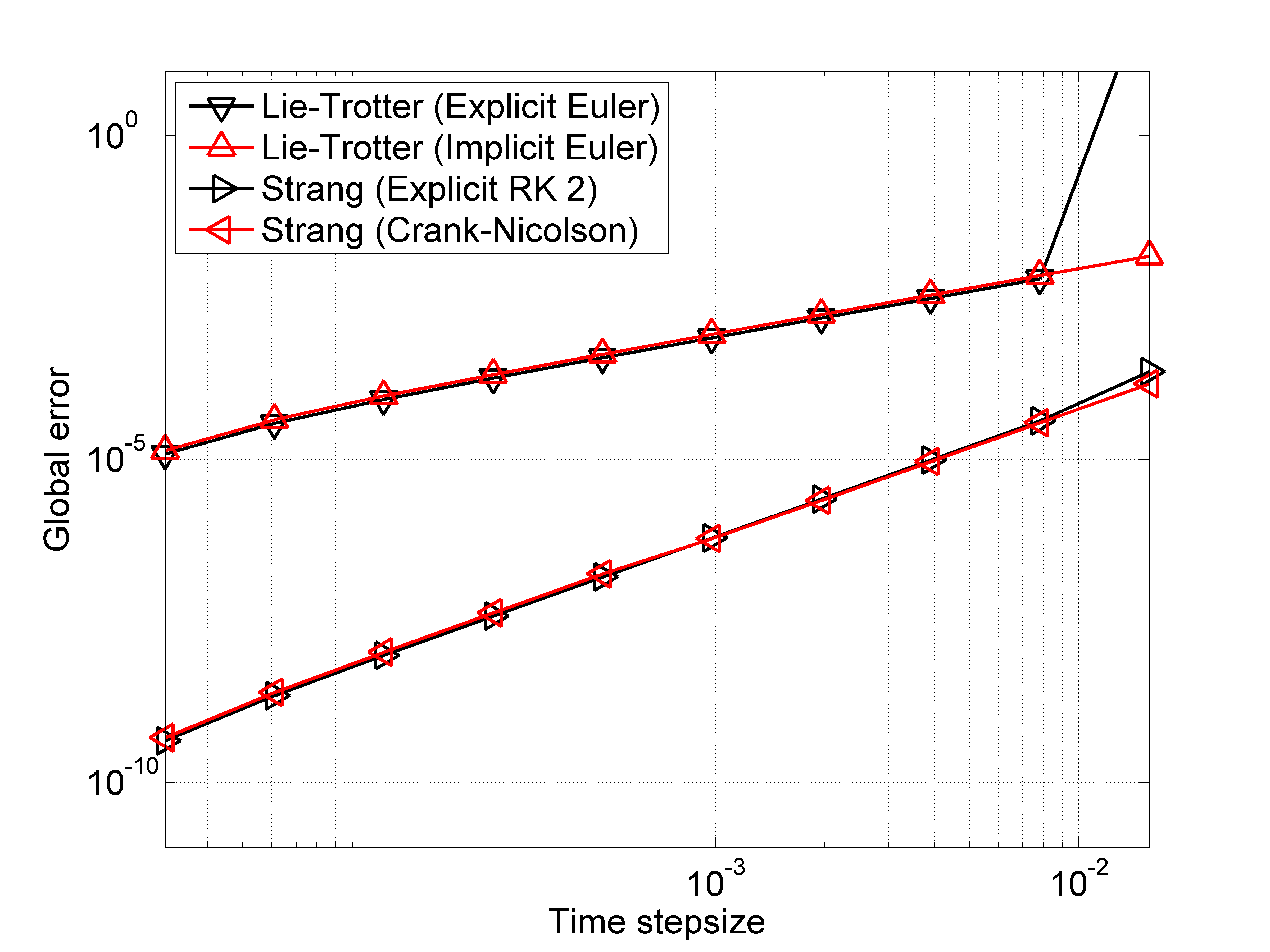} \\[-3.5cm]
\end{center}
\caption{Local (left) and global (right) errors for the Lie--Trotter and Strang splitting methods with respect to the $H^3 \times H^1$-norm obtained for Decomposition~I. 
Comparison of different time integration methods for the numerical solution of the subproblems.}
\label{fig:Figure2}
\end{figure}

\begin{figure}[t!]
\begin{center}
\includegraphics[width=0.49\textwidth]{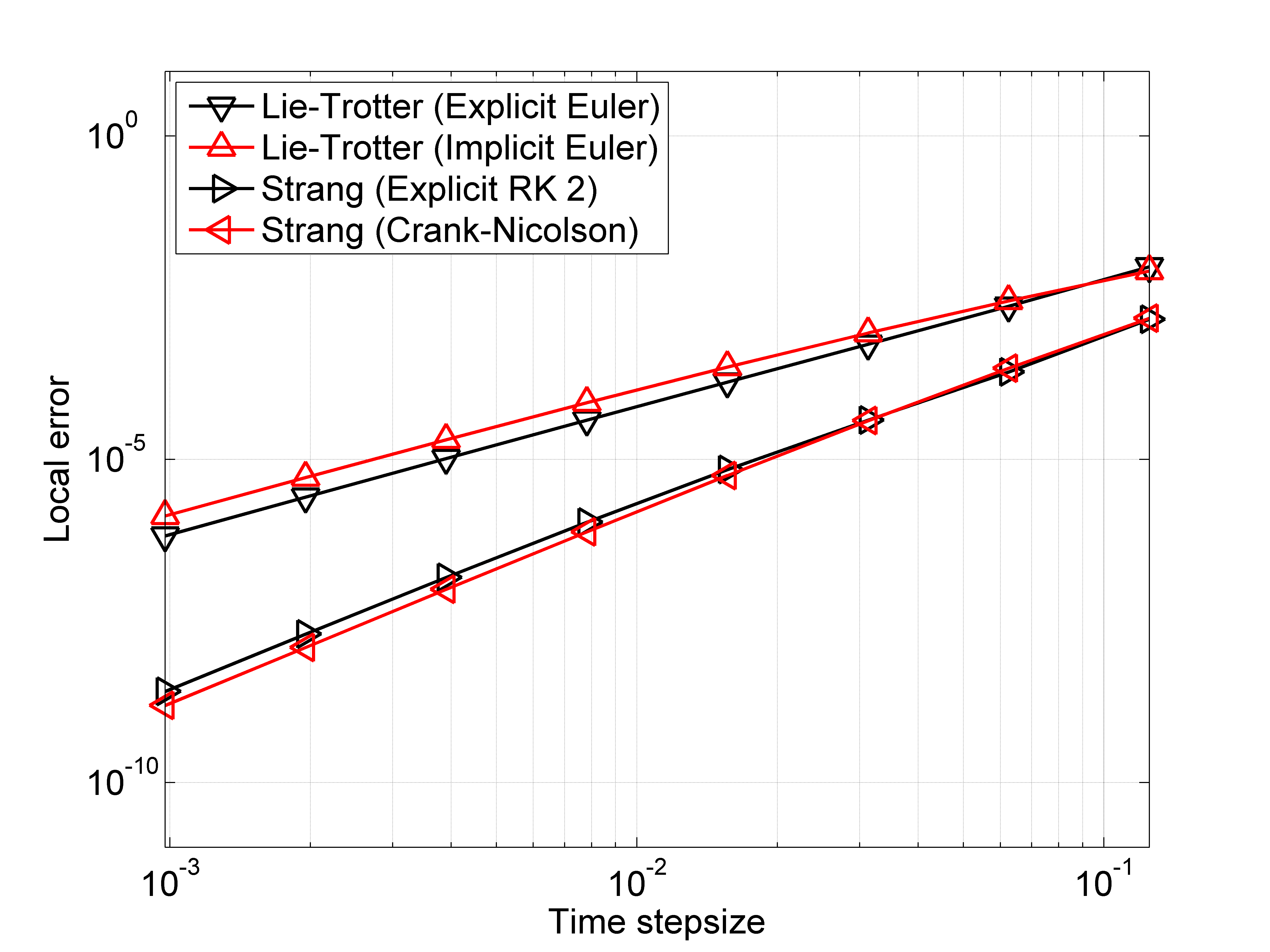} \includegraphics[width=0.49\textwidth]{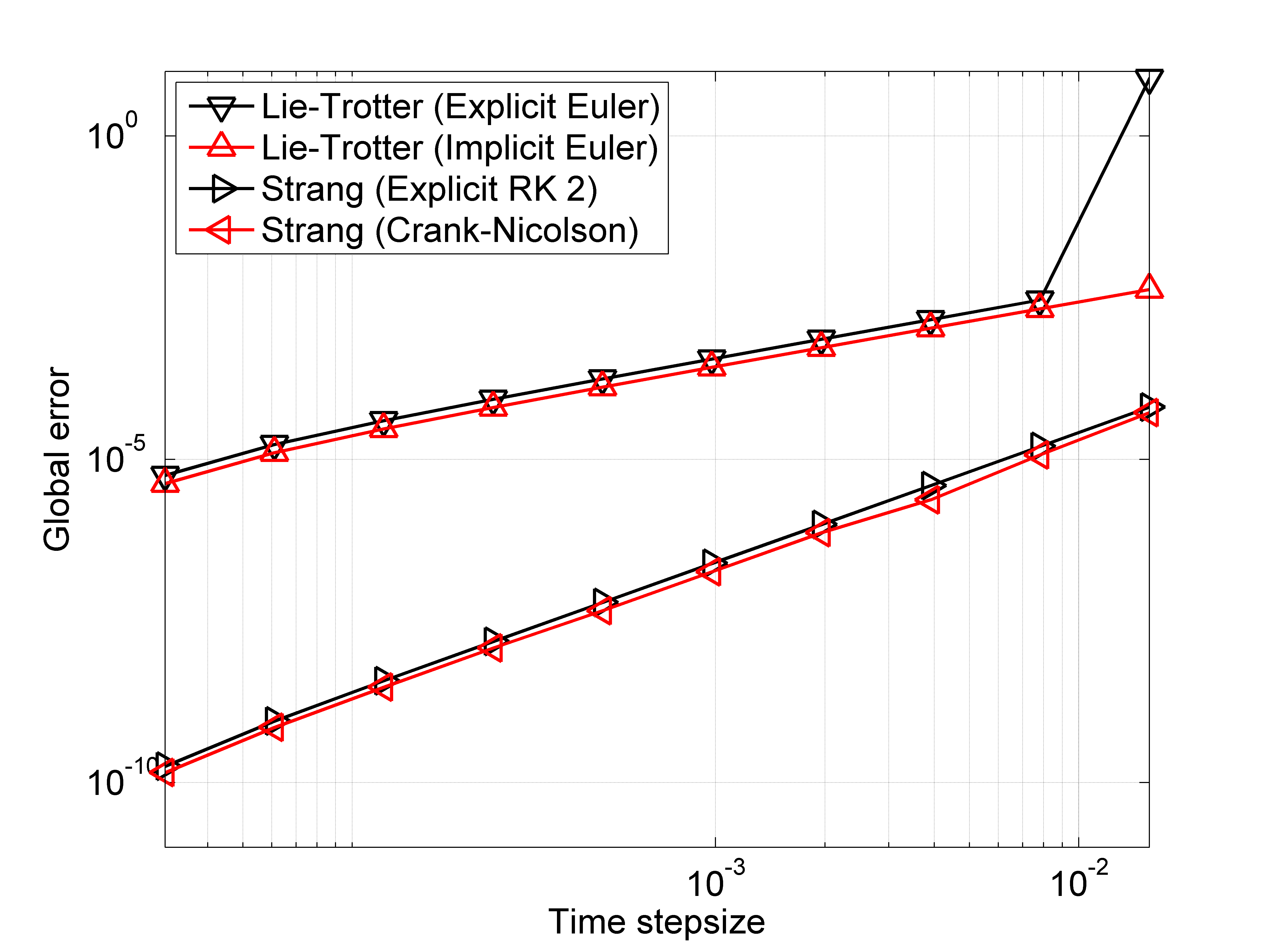} \\
\includegraphics[width=0.49\textwidth]{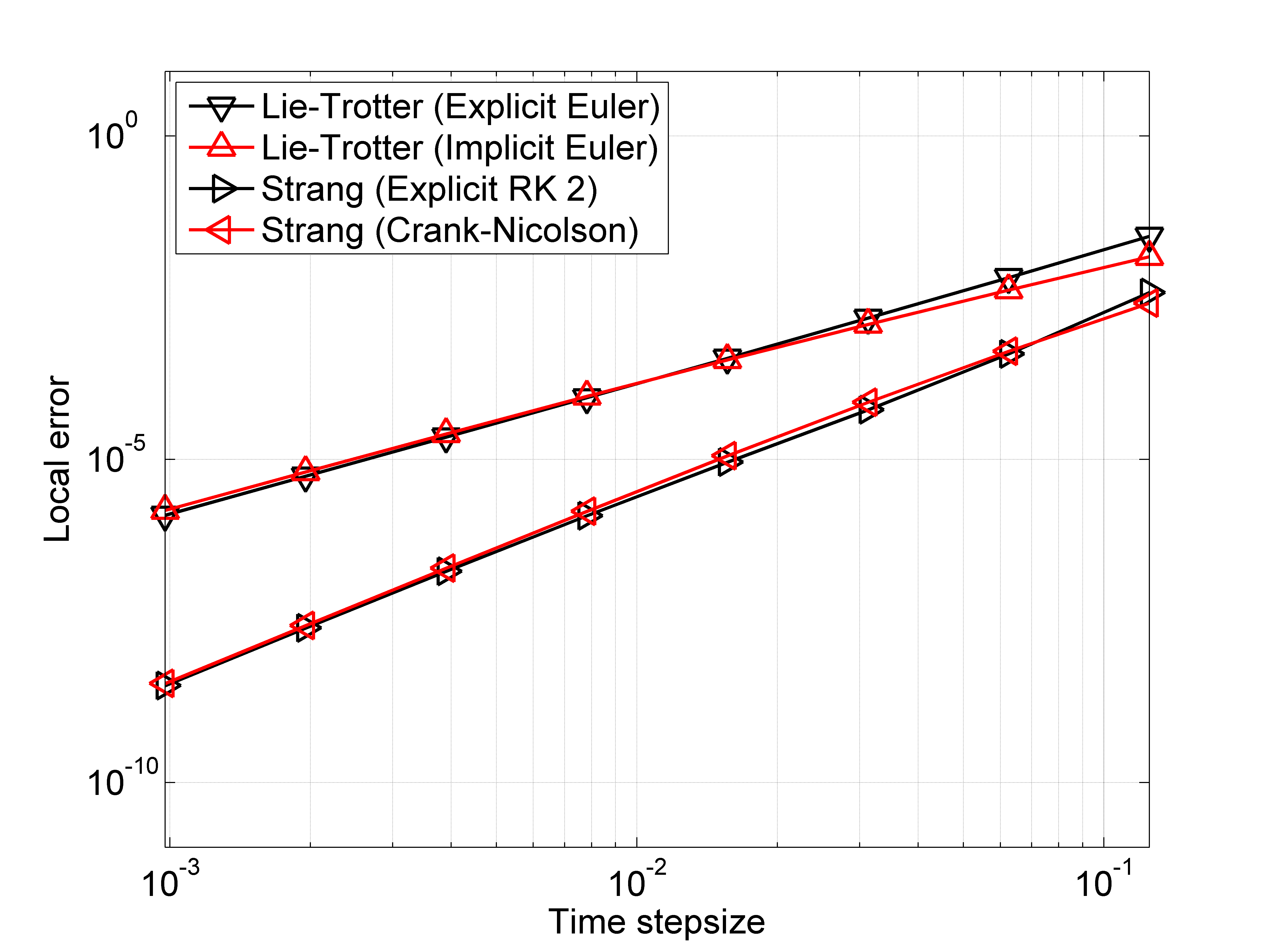} \includegraphics[width=0.49\textwidth]{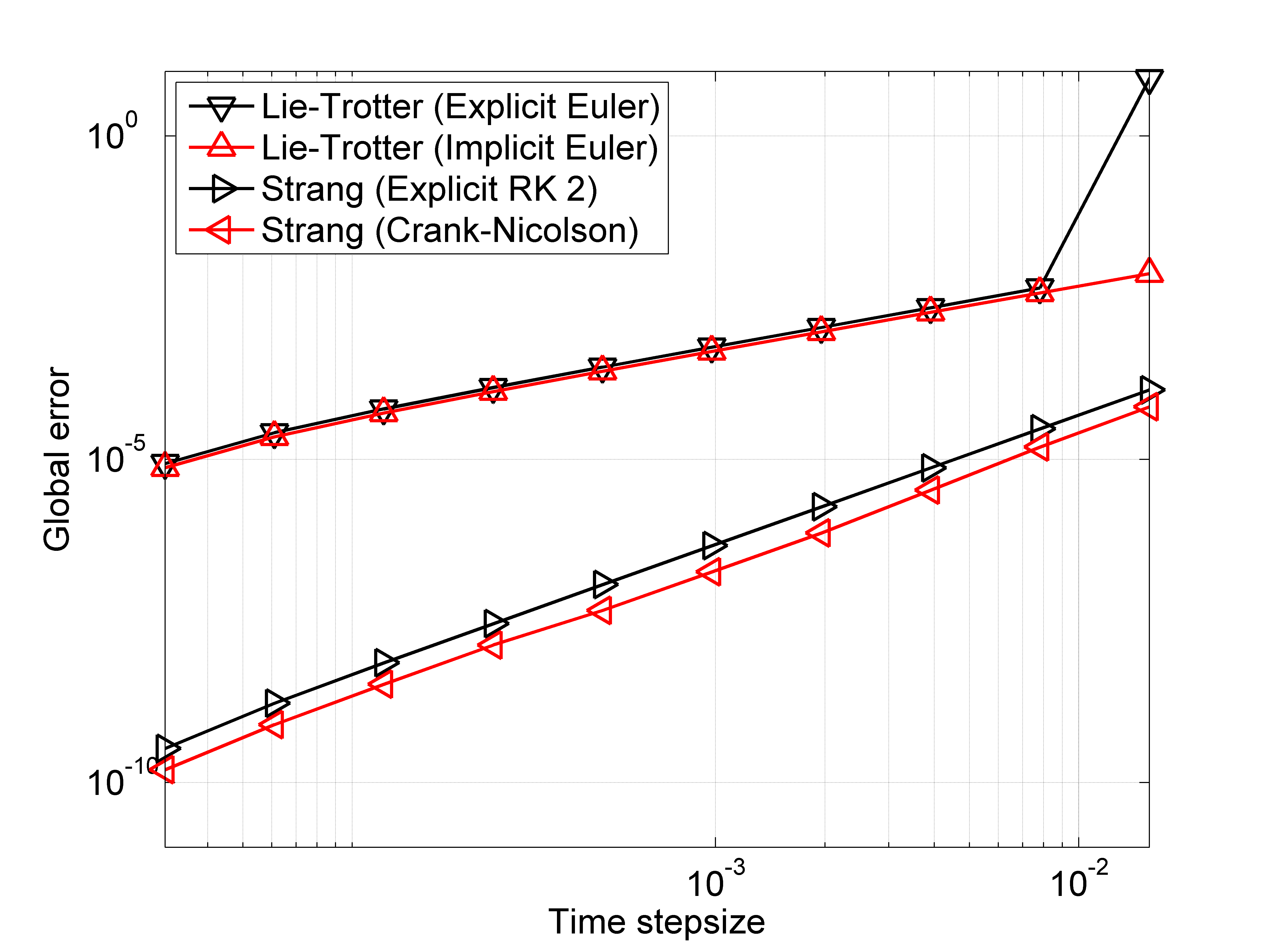} \\
\includegraphics[width=0.49\textwidth]{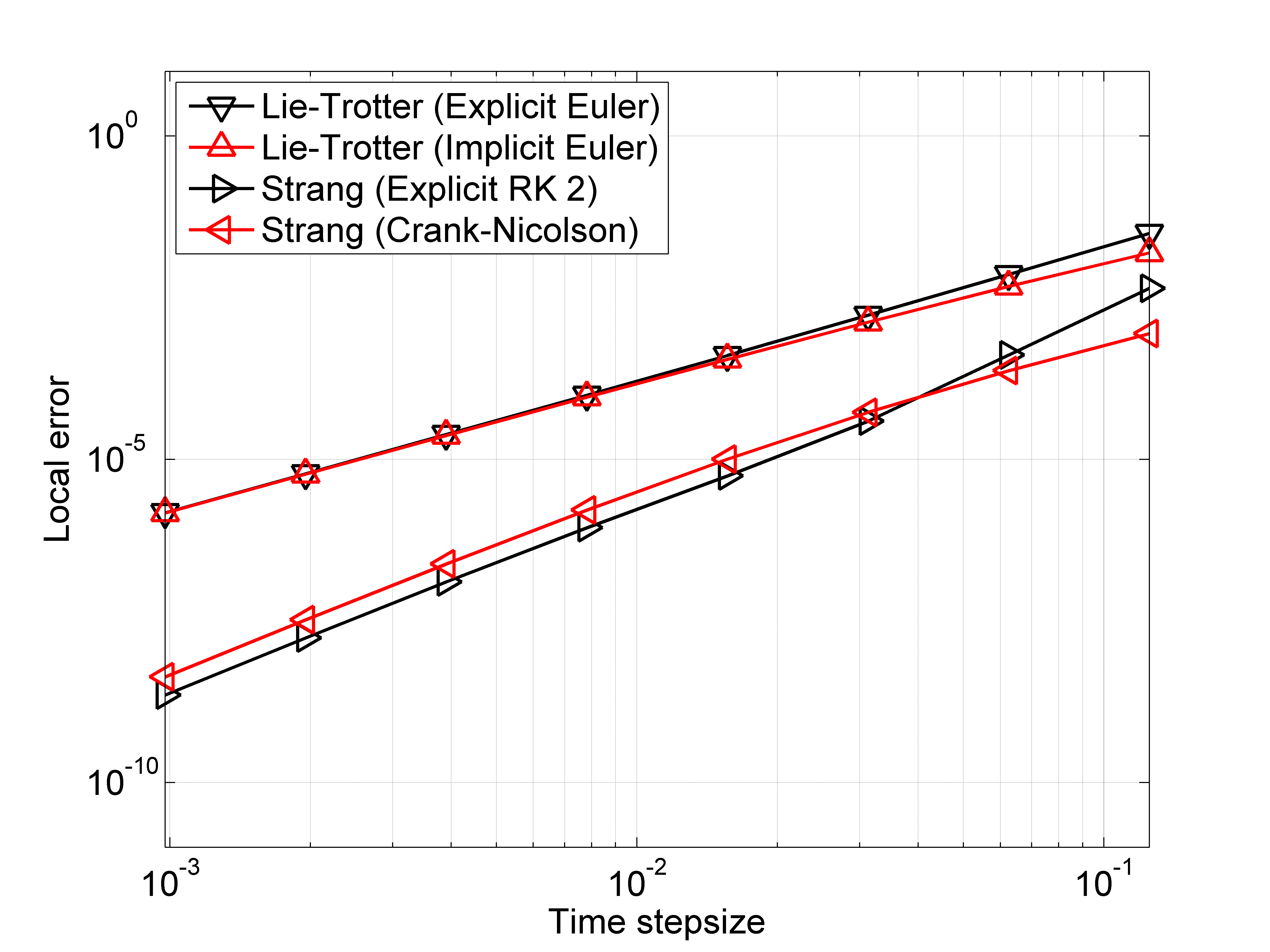} \includegraphics[width=0.49\textwidth]{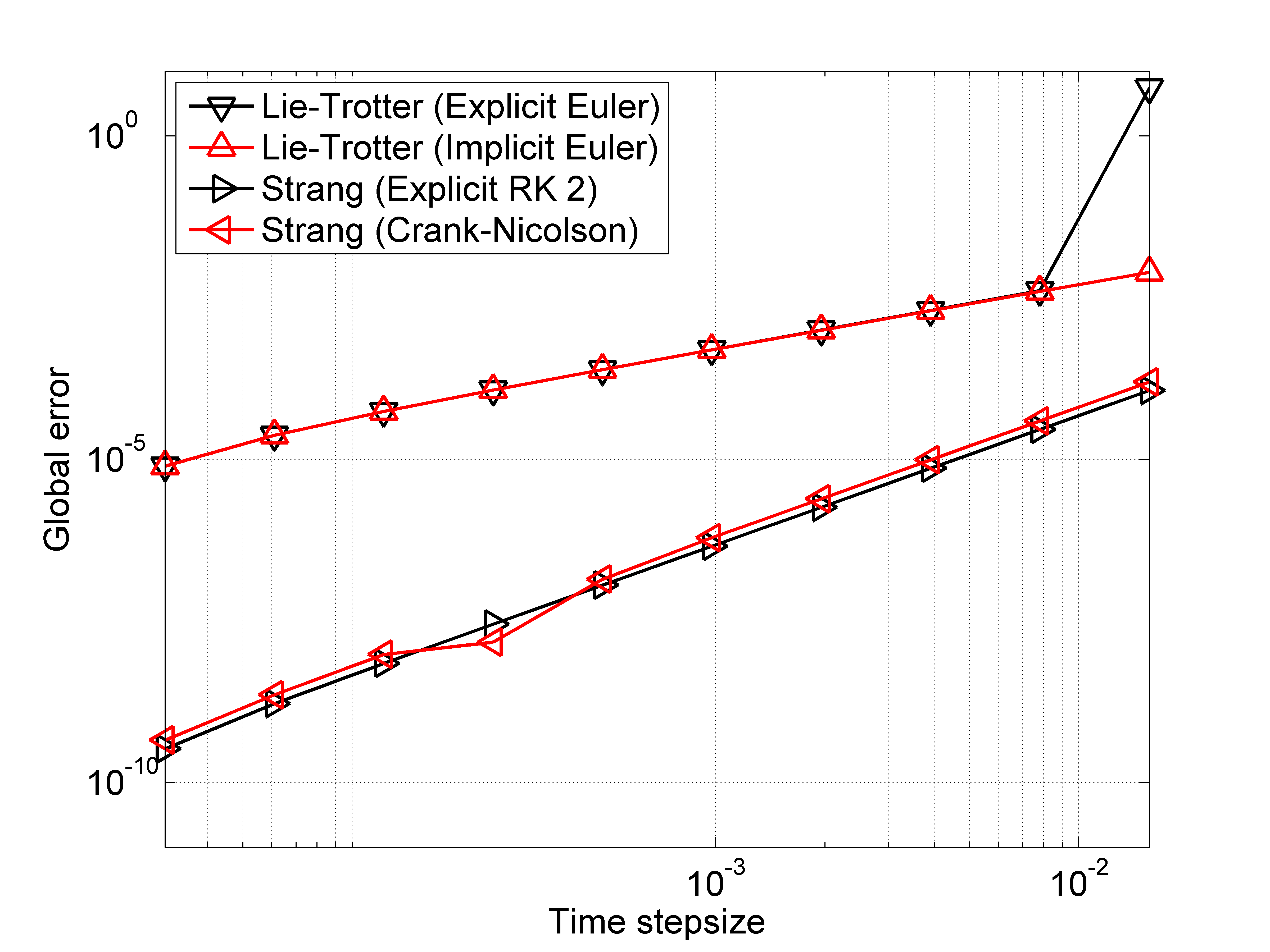} 
\end{center}
\caption{Local (left) and global (right) errors for the Lie--Trotter and Strang splitting methods with respect to the $L^2 \times L^2$-norm obtained for Decomposition~II (first row), Decomposition~III (second row) and Decomposition~IV (third row). 
Comparison of different time integration methods for the numerical solution of the subproblems.}
\label{fig:Figure1b}
\end{figure}

\begin{figure}[t!]
\begin{center}
\includegraphics[width=0.7\textwidth]{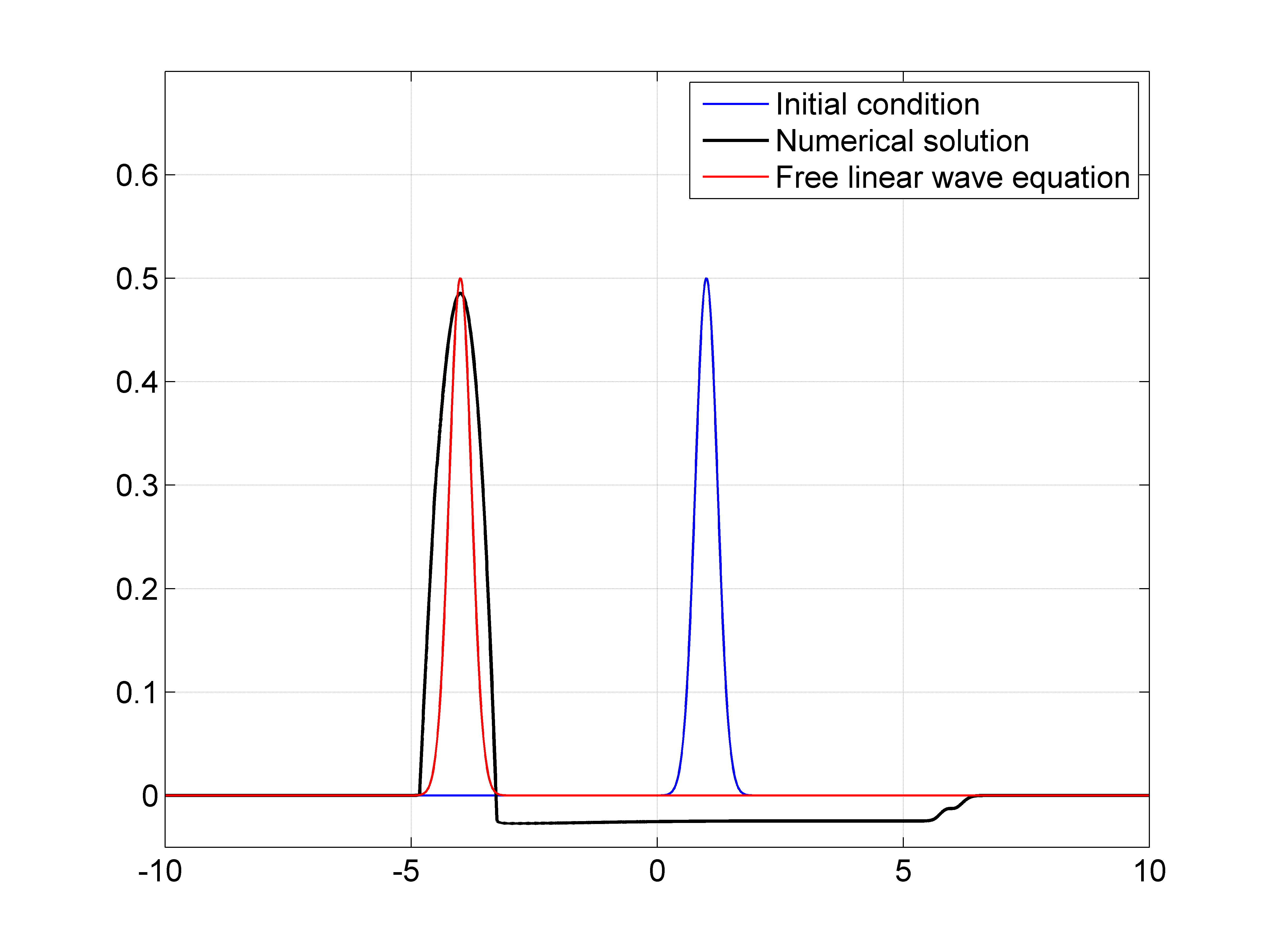}
\end{center}
\caption{Solution to the Westervelt equation with parameter values~\eqref{realisticparameters}, computed by the Lie--Trotter splitting method based on Decomposition~I. 
Comparison of the acoustic velocity potential with the solution to the linear wave equation at time $5 \cdot 10^{-3}$.}
\label{fig:Figure3}
\end{figure}

\paragraph{Model problem.} 
We consider the one-dimensional Westervelt equation 
\begin{subequations}
\label{eq:ModelProblem}
\begin{equation}
\begin{split}
&\patt \ps(x,t) - \a \, \pa_{xxt} \ps(x,t) - \b \, \pa_{xx} \ps(x,t) \\
&\quad= \d \, \pat \ps(x,t) \, \patt \ps(x,t)\,, \quad (x,t) \in (-a,a) \times (0,T]\,, \\
\end{split}
\end{equation}
with parameter values 
\begin{equation}
\a = 1\,, \qquad \b = 1\,, \qquad \g = \tfrac{1}{2}\,, \qquad \d = 2 \g = 1\,, 
\end{equation}
subject to homogeneous Dirichlet conditions 
\begin{equation}
\ps(-a,t) = 0 = \ps(a,t)\,, \quad t \in [0,T]\,, 
\end{equation}
and regular initial conditions 
\begin{equation}
\ps(x,0) = \ee^{- x^2}\,, \quad \pat \ps(x,0) = - \, x \, \ee^{- x^2}\,, \qquad x \in [-a,a]\,.
\end{equation}
\end{subequations}
For the numerical computations we set $a = 8$, $T = 1$, and choose $M = 100$ equidistant space grid points. 
In Figures~\ref{fig:Figure1a} and~\ref{fig:Figure1b} we display the local and global errors with respect to the $L^2 \times L^2$-norm, obtained for the Lie-Trotter and the Stang splitting methods based on Decompositions~I-IV. 
Accordingly to Theorem~\ref{thm:MainResult}, for Decomposition~I we also include the errors with respect to a stronger norm, the $H^3 \times H^1$-norm, see Figure~\ref{fig:Figure2}.
In all cases, the slopes of the lines reflect the convergence orders $p = 1$ for the Lie--Trotter splitting method and $p = 2$ for the Strang splitting method, respectively.  
From the numerical results we conclude that the application of Decompositions~II-IV does not improve the size of the errors; 
for this reason, we favour Decomposition~I with the least computational effort. 

\paragraph{Realistic parameter values.}
As a further illustration, we consider the one-dimensional Westervelt equation with more realistic parameter, geometry, and excitation values
\begin{equation}
\label{realisticparameters}
\begin{gathered}
a = 15\,, \quad \a = 10^{-2}\,, \quad c = 10^3\,, \quad \b = c^2 = 10^6\,, \, \quad \d = 2 \cdot 10^{-4}\,,\\ 
\ps(x,0) = 0.5 \, \ee^{-10 (x-1)^2}\,, \quad x \in [-a,a]\,, \\ 
\pat \ps(x,0) = c \, \pax \ps(x,0) = - 10 \, c \, (x-1) \, \ee^{-10 (x-1)^2}\,, \quad x \in [-a,a]\,,
\end{gathered}
\end{equation}
in the MKS system of units. 
For the time discretisation we apply the Lie--Trotter splitting method, based on Decomposition~I and the explicit Euler method, and for the space discretisation we use Fast Fourier techniques.  
In Figure~\ref{fig:Figure3}, we display the numerical result for the acoustic velocity potential~$\ps$ at time $t = 5 \cdot 10^{-3}$, obtained for $M = 6 \cdot 10^3$ space grid points and $N = 5 \cdot 10^4$ time steps. 
In order to reveal the effect of nonlinearity, we also display the profile of the solution to the linear wave equation $\patt \ps - c^2 \D \, \ps = 0$;
as typical in the nonlinear case, a slight one-sided steepening of the pulse over time while the wave travels to the left is observed, see also~\cite[Section 5.6.4]{Kaltenbacher2007}.
\section{Regularity results}
\label{sec:Regularity}
In this section, we deduce regularity results for the Westervelt equation and related evolution equations of hyperbolic and parabolic type. 
The stated results, ensuring that certain smoothness properties of the initial state are inherited by the solution to the original problem, the subproblems, and the associated variational equations, respectively, are essential ingredients in our convergence analysis of operator splitting methods for the time integration of the Westervelt equation.
We study the Westervelt equation subject to homogeneous Dirichlet boundary conditions on a regular domain $\Om \subset \RR^3$. 
The results can be extended to other boundary conditions, in the context of the Westervelt equation see for instance~\cite{ClasonKaltenbacherVeljovic2009};
they litterally carry over to the cases $\Om \subset \RR^d$ with $d = 1,2$, where in a single space dimension the continuous embedding $H^{1}(\Om) \hookrightarrow \nC(\Om)$ is utilised to deal with nonlinearities.  
Our approach for the Westervelt equation utilises a regularity result proven in~\cite{KaltenbacherLasiecka2009} for a reformulation in terms of the acoustic pressure. 
Although results on well posedness and regularity with respect to space and time can be found in the literature, see for instance~\cite[Thm~6, p.\,412]{Evans} for the hyperbolic case and~\cite[Thm~6, p.\,386]{Evans} for the parabolic case, we carry out the proofs of the stated results, since on one hand we only need local in time estimates as well as high spatial but low temporal regularity and on the other hand coefficients will depend on time in our setting. 
We note that estimates such as 
\begin{equation*}
\normbig{\nC([0,T],H^{k+5} \times H^{k+4})}{\big(\ps, \pat \ps\big)} \leq C\,, \qquad k \in \NN_{\geq 0}\,,
\end{equation*}
are at first deduced for even exponents $k+5 = 2 \ell$ with $\ell \in \NN_{\geq 0}$; 
the corresponding assertions for non-even and, more generally, for non-integer exponents then follow from the exact interpolation theorem, which states that boundedness of $F: X \to Y$ and $F: \Xtil \to \Ytil$ implies boundedness of $F: (X, \Xtil)_{\vartheta,q} \to (Y, \Ytil)_{\vartheta,q}$ for $0 < \vartheta < 1$ and $1 < q < \infty$, see~\cite[Thm~7.23]{AdamsFournier}.
\subsection{Westervelt equation} 
\paragraph{Initial-boundary value problem and reformulation.}
We consider~\eqref{eq:WesterveltOriginal} under homogeneous Dirichlet boundary conditions on a regular domain
\begin{equation}
\label{eq:Westervelt}
\begin{cases}
&\big(1 - \d \, \pat \ps(x,t)\big) \, \patt \ps(x,t) - \a \, \D \, \pat \ps(x,t) - \b \, \D \, \ps(x,t) \\
&\quad= 0\,, \quad (x,t) \in \Om \times (0,T]\,, \\
&\ps(x,t) = 0\,, \quad (x,t) \in \pa \Om \times [0,T]\,, \\
&\ps(x,0) = \ps^0(x)\,, \quad \pat \ps(x,0) = \ps^1(x)\,, \quad x \in \Om\,.
\end{cases} 
\end{equation}
As a first step, we formulate the Westervelt equation in terms of the acoustic pressure, related to the acoustic velocity potential via 
\begin{subequations}
\label{eq:WesterveltPressure}
\begin{equation}
p(x,t) = \r \, \pat \ps(x,t)\,, \quad \pat \ps(x,t) = \tfrac{1}{\r} \, p(x,t)\,, \quad (x,t) \in \overline{\Om} \times [0,T]\,, 
\end{equation}
for a positive constant $\r > 0$.
Differentation of the Westervelt equation yields
\begin{equation*}
- \, \d \, \big(\patt \ps(x,t)\big)^2 + \big(1 - \d \, \pat \ps(x,t)\big) \, \pattt \ps(x,t) - \a \, \D \, \patt \ps(x,t) - \b \, \D \, \pat \ps(x,t) = 0\,;
\end{equation*}
substituting $\pat \ps = \tfrac{1}{\r} \, p$ as well as $\patt \ps = \tfrac{1}{\r} \, \pat p$ and $\pattt \ps = \tfrac{1}{\r} \, \patt p$ implies 
\begin{equation*}
- \tfrac{\d}{\r^2} \, \big(\pat p(x,t)\big)^2 + \big(1 - \tfrac{\d}{\r} \, p(x,t)\big) \, \tfrac{1}{\r} \, \patt p(x,t) - \tfrac{\a}{\r} \, \D \, \pat p(x,t) - \tfrac{\b}{\r} \, \D p(x,t) = 0\,.
\end{equation*}
Altogether, we obtain the reformulation
\begin{equation}
\label{eq:Westervelt_p}
\begin{cases}
&\big(1 - \tfrac{\d}{\r} \, p(x,t)\big) \, \patt p(x,t) - \a \, \D \, \pat p(x,t) - \b \, \D p(x,t) \\
&\quad= \tfrac{\d}{\r} \big(\pat p(x,t)\big)^2\,, \quad (x,t) \in \Om \times (0,T]\,, \\
&p(x,t) = 0\,, \quad (x,t) \in \pa \Om \times [0,T]\,, \\
&p(x,0) = \r \, \ps^1(x)\,, \quad x \in \Om\,, \\
&\pat p(x,0) = \r \, \big(1 - \d \ps^1(x)\big)^{-1} \, \big(\a \, \D \ps^1(x) + \b \D\ps^0(x)\big)\,, \quad x \in \Om\,.
\end{cases} 
\end{equation}
\end{subequations}

\paragraph{Non-degeneracy and regularity result.}
Provided that the prescribed initial state satisfies the regularity assumptions 
\begin{subequations}
\label{eq:ThmRegularityp}
\begin{equation}
p(\cdot,0)\in H^2(\Om) \,, \quad \pat p(\cdot,0) \in H^1(\Om)\,, \quad \patt p(\cdot,0) \in L^2(\Om)\,, 
\end{equation}
with sufficiently small norms, the regularity result~\cite[Thm 1.1]{KaltenbacherLasiecka2009} ensures non-degeneracy as well as local existence and uniqueness of a solution to a weak form of~\eqref{eq:Westervelt_p}, that is, there exist constants $0 < \underline{\nu} < \overline{\nu} < \infty$ such that the relation
\begin{equation}
\label{eq:ThmRegularityp1}
\underline{\nu} \leq 1 - \tfrac{\d}{\r} \, p(x,t) \leq \overline{\nu}\,, \qquad (x,t) \in \overline{\Om} \times [0,T]\,, 
\end{equation}
holds and furthermore 
\begin{equation}
\begin{split}
&p \in \nC\big([0,T], H^2(\Om)\big) \cap \nC^1\big([0,T], H^1(\Om)\big) \cap \nC^2\big([0,T], L^2(\Om)\big) \\
&\qquad\qquad \cap H^2\big((0,T), H^1(\Om)\big)\,.
\end{split}
\end{equation}
\end{subequations}

\paragraph{Regularity result.}
The above statement permits to deduce the following regularity result for the Westervelt equation.
We note that the proof of Theorem~\ref{thm:TheoremRegularity} via ~\cite[Thms 1.2, 1.3]{KaltenbacherLasiecka2009} ensures well-posedness, globally in time, that is, there exist positive constants $\rho > 0$ and $M > 0$ such that the energy norm remains bounded 
\begin{equation*}
E_{\ps,1}(t) = \tfrac{1}{2} \Big(\normbig{L^2}{\pattt \ps(\cdot,t)}^2 + \normbig{L^2}{\nabla \patt \ps(\cdot,t)}^2 + \normbig{L^2}{\D \, \pat \ps(\cdot,t)}^2\Big) \leq M\,, \qquad t \geq 0\,, 
\end{equation*}
provided that $E_{\ps,1}(0) \leq \rho$;
moreover, due to the strong damping present in the equation, the energy decays exponentially.
However, for our purposes it suffices to prove local well-posedness and regularity results.

\begin{theorem}
\label{thm:TheoremRegularity}
Let $T > 0$ and $M > 0$ be arbitrary. 
\begin{enumerate}[(i)]
\item 
There exists a constant $\rho > 0$ such that if the initial data fulfill the regularity and compatibility assumptions 
\begin{equation*}
\begin{gathered}
\ps^0, \ps^1\in H_0^1(\Om) \cap H^4(\Om)\,, \qquad \normbig{H^4 \times H^4}{\big(\ps^0, \ps^1\big)} = \normbig{H^4}{\ps^0} + \normbig{H^4}{\ps^1} \leq \rho\,, \\
\patt \ps (\cdot, 0) = \big(1 - \d \ps^1\big)^{-1} \big(\a \, \D \ps^1 + \b \, \D \ps^0\big) \in H_0^1(\Om) \cap H^2(\Om)\,,
\end{gathered}
\end{equation*}
then the solution to the Westervelt equation~\eqref{eq:Westervelt} exists and satisfies the non-degeneracy condition  
\begin{equation*}
\underline{\nu} \leq 1 - \d \, \pat \ps(x,t) \leq \overline{\nu}\,, \qquad (x,t) \in \overline{\Om} \times [0,T]\,, 
\end{equation*}
with some constants~$0 < \underline{\nu} < \overline{\nu} < \infty$;
moreover, the relation 
\begin{equation*}
\begin{split}
&\ps \in 
\nC^1\big([0,T], H^4(\Om)\big) \cap \nC^2\big([0,T], H^2(\Om)\big) \cap \nC^3\big([0,T], L^2(\Om)\big) \\
&\qquad\qquad \cap H^3\big((0,T), H^1(\Om)\big) 
\end{split}
\end{equation*}
and the bound 
\begin{equation*}
\begin{split}
&\normbig{H^4 \times H^4\times H^2}{\big(\ps(\cdot,t),\pat \ps(\cdot,t),\patt \ps(\cdot,t)\big)} 
\leq \ee^{C t} \, \normbig{H^4 \times H^4}{\big(\ps^0,\ps^1\big)}\,, 
\\&\qquad t \in [0,T]\,,
\end{split}
\end{equation*}
are valid. 
\item 
For every $m \in \NN_{\geq 2}$ there exists a constant $\rho > 0$ such that if the initial data fulfull the regularity and compatibility assumptions 
\begin{equation*}
\begin{gathered}
\big(\ps^0, \ps^1\big) \in H^{2(m+1)}(\Om) \times H^{2m+1}(\Om)\,, \qquad \normbig{H^{2(m+1)} \times H^{2m+1}}{\big(\ps^0, \ps^1\big)} \leq \rho\,, \\
\patt \D^k \ps(\cdot, 0) = \big(1 - \d \ps^1\big)^{-1} \big(\a \D^{k+1} \ps^1 + \b \D^{k+1} \ps^0\big) + f_k(\cdot,0) \in H_0^1(\Om)\,,
\end{gathered}
\end{equation*}
for any integer $1 \leq k \leq m$, where 
\begin{equation*}
f_k = \D^k \Big(\big(1 - \d \pat \ps\big)^{-1} \D (\a \, \pat \ps + \b \, \ps)\Big) - \big(1 - \d \pat \ps\big)^{-1} \D^{k+1} (\a \, \pat \ps + \b \, \ps)\,, 
\end{equation*}
then the solution to~\eqref{eq:Westervelt} satisfies the relation 
\begin{equation*}
\big(\ps,\pat \ps\big) \in \nC\big([0,T],H^{2(m+1)}(\Om) \times H^{2m+1}(\Om)\big)
\end{equation*}
and the bound 
\begin{equation*}
\normbig{H^{2(m+1)} \times H^{2m+1}}{\big(\ps(\cdot,t), \pat \ps(\cdot,t)\big)} \leq \ee^{Ct} \, \normbig{H^{2(m+1)} \times H^{2m+1}}{\big(\ps^0,\ps^1\big)}\,.
\end{equation*}
\end{enumerate}
\end{theorem}
\emph{Proof.}
For notational simplicity, for a function $f: \overline{\Om} \times [0,T] \to \RR$ we henceforth write $f(t) = f(\cdot,t)$ for short. 
\begin{enumerate}[(i)] 
\item
(a) \emph{Non-degeneracy.} 
Inserting $p = \r \, \pat \ps$ into~\eqref{eq:ThmRegularityp1} immediately implies 
\begin{equation*}
0 < \underline{\nu} \leq \big(1 - \d \, \pat \ps(t)\big) \leq \overline{\nu} < \infty\,, \qquad t \in [0,T]\,.
\end{equation*}
(b) \emph{Basic regularity.}
For convenience, we recall the relations 
\begin{equation*}
\begin{split}
\pat p &= \r \, \patt \ps = \r \, \big(1 - \d \, \pat \ps\big)^{-1} \, \big(\a \, \D \, \pat \ps + \b \, \D \, \ps\big)\,, \\
\patt p &= \big(1 - \tfrac{\d}{\r} \, p\big)^{-1} \, \big(\a \, \D \, \pat p + \b \, \D p + \tfrac{\d}{\r} \big(\pat p\big)^2\big)\,, 
\end{split}
\end{equation*}
see also~\eqref{eq:Westervelt}--\eqref{eq:WesterveltPressure}.
The assumption $\ps^0, \ps^1 \in H^4(\Om)$ ensures that the initial data  
\begin{equation*}
\begin{split}
p(0) &= \r \, \ps^1\,, \\
\pat p(0) &= \r \, \big(1 - \d \, \ps^1\big)^{-1} \, \big(\a \, \D \ps^1 + \b \, \D \ps^0\big)\,, \\
\patt p(0) &= \big(1 - \tfrac{\d}{\r} \, p(0)\big)^{-1} \, \big(\a \, \D \, \pat p(0) + \b \, \D p(0) + \tfrac{\d}{\r} \, \big(\pat p(0)\big)^2\big) \\
&= \r \, \big(1 - \d \, \ps^1\big)^{-1} \, 
\Big(\a \, \D \, \big(\big(1 - \d \, \ps^1\big)^{-1} \, \big(\a \, \D \ps^1 + \b \, \D \ps^0\big)\big) 
+ \b \, \D \, \ps^1 \\
&\qquad+ \d \, \big(1 - \d \, \ps^1\big)^{-2} \, \big(\a \, \D \ps^1 + \b \, \D \ps^0\big)^2\Big)\,, 
\end{split}
\end{equation*}
chosen accordingly, fulfill the requirements
\begin{equation*}
p(0) \in H^2(\Om)\,, \qquad \pat p(0) \in H^1(\Om)\,, \qquad \patt p(0)  \in L^2(\Om)\,; 
\end{equation*}
an application of the regularity result~\eqref{eq:ThmRegularityp} thus implies 
\begin{equation*}
\begin{split}
&p \in \nC\big([0,T], H^2(\Om)\big) \cap \nC^1\big([0,T], H^1(\Om)\big) \cap \nC^2\big([0,T], L^2(\Om)\big) \\
&\qquad\qquad \cap H^2\big((0,T), H^1(\Om)\big)\,.
\end{split}
\end{equation*}
Together with the relation 
\begin{equation*}
\ps(t) = \ps^0 + \tfrac{1}{\rho} \int_{0}^{t} p(\t) \; \dd\t\,,
\end{equation*}
obtained by integrating $\pat \ps = \tfrac{1}{\r} \, p$ in time, this proves 
\begin{equation}
\label{eq:BasicRegularity}
\begin{split}
&\ps \in \nC^1\big([0,T], H^2(\Om)\big) \cap \nC^2\big([0,T], H^1(\Om)\big) \cap \nC^3\big([0,T], L^2(\Om)\big) \\
&\qquad\qquad\cap H^3\big((0,T), H^1(\Om)\big)\,.
\end{split}
\end{equation}
(c) \emph{Additional regularity.} 
Continuous Sobolev embeddings~\eqref{eq:AuxiliaryEstimates} imply 
\begin{equation*}
\normbig{H^1}{\big(1 - \d \, \pat \ps(t)\big) \, \patt \ps(t)} \leq C \, \big(1 + \norm{\nC^1([0,T],H^2)}{\ps}\big) \, \norm{\nC^2([0,T],H^1)}{\ps}\,, 
\end{equation*}
and thus by~\eqref{eq:Westervelt} we get 
\begin{equation*}
f = \a \, \D \, \pat \ps + \b \, \D \, \ps =\big(1 - \d \, \pat \ps\big) \, \patt \ps \in \nC\big([0,T], H^1(\Om)\big)\,.
\end{equation*}
Making use of the fact that~$\D \, \ps$ fulfills the differential equation 
\begin{subequations}
\label{eq:ProofVocf}
\begin{equation}
\pat \D \, \ps(t) = - \tfrac{\b}{\a} \, \D \, \ps(t) + \tfrac{1}{\a} \, f(t)\,, 
\end{equation}
employing the linear variation-of-constants formula 
\begin{equation}
\D \, \ps (t) = \ee^{-\frac{\b}{\a}t} \D \ps^0 + \tfrac{1}{\a} \int_{0}^{t} \ee^{-\frac{\b}{\a}(t-s)} f(s) \; \dd s\,, 
\end{equation}
\end{subequations}
and recalling that the initial state satisfies $\ps^0 \in H^4(\Om)$, we conclude that the relations $\D \, \ps(t) \in H^1(\Om)$ and by~\eqref{eq:ProofVocf}
$\D \, \pat \ps(t) \in H^1(\Om))$ are valid, that is, we obtain   
\begin{equation*}
\big(\ps, \pat \ps\big) \in \nC\big([0,T], H^3(\Om) \times H^3(\Om)\big)\,.
\end{equation*}
A reformulation of the Westervelt equation and differentiation in time leads to 
\begin{equation*}
\begin{split}
- \, \a \, \D \, \pat \ps &= \b \, \D \, \ps - \big(1- \d \, \pat \ps\big) \, \patt \ps\,, \\
- \, \a \, \D \, \patt \ps &= \b \, \D \, \pat \ps + \d \, \big(\patt \ps\big)^2 - \big(1- \d \, \pat \ps\big) \, \pattt \ps\,, \\
\a \, \patt \ps &= (-\D)^{-1} \, \big(- \a \, \D \, \patt \ps\big) \\
&= - \b \, \pat \ps + (-\D)^{-1} \Big(\d \, \big(\patt \ps\big)^2 - \big(1- \d \, \pat \ps\big) \, \pattt \ps\Big)\,;
\end{split} 
\end{equation*}
due to the the compatibility condition $\patt \ps(0) \in H^1_0(\Om) \cap H^2(\Om)$ an application of the Laplace operator under homogeneous Dirichlet boundary conditions to $\patt \ps$ is justified. 
Consequently, the relation 
\begin{equation*}
\begin{split}
f 
&=  \big(1 - \d \, \pat \ps\big) \, \patt \ps \\
&= \tfrac{1}{\a} \, \big(1 - \d \, \pat\ps\big) \, \Big(- \b \, \pat \ps \\
&\qquad+ (-\D)^{-1} \Big(\d \, \big(\patt \ps\big)^2 - \big(1- \d \, \pat \ps\big) \, \pattt \ps\Big)\Big)
\end{split}
\end{equation*}
follows. 
By the above considerations $\pat \ps\in \nC\big([0,T], H^3(\Om)\big)$ and by~\eqref{eq:BasicRegularity} the relations $\patt \ps \in \nC([0,T], H^1(\Om))$ and $\pattt \ps \in \nC\big([0,T], L^2(\Om)\big)$ hold;  
together with the regularisation $(-\D)^{-1}: L^2(\Om)\to H^2(\Om)$ and continuous Sobolev embeddings~\eqref{eq:AuxiliaryEstimates} we obtain 
\begin{equation*}
\begin{split}
\norm{H^2}{f(t)} &\leq C \, \normbig{H^2}{1 - \d \, \pat \ps(t)} \Big(\normbig{H^2}{\pat \ps(t)} + \normbig{H^1}{\patt \ps(t)}^2 \\
&\qquad+ \normbig{H^2}{1- \d \, \pat \ps(t)} \, \norm{L^2}{\pattt \ps(t)}\Big)\,, 
\end{split}
\end{equation*}
which implies 
\begin{equation*}
f \in \nC\big([0,T], H^2(\Om)\big)\,.
\end{equation*}
Proceeding as above, via~\eqref{eq:ProofVocf} we conclude that $\D \, \ps\,, \D \, \pat \ps \in \nC([0,T], H^2(\Om))$ holds and finally obtain 
\begin{equation*}
\big(\ps,\pat \ps\big) \in \nC\big([0,T], H^4(\Om) \times H^4(\Om)\big)\,.
\end{equation*}
(d) \emph{Additional regularity.} 
We utilise the relation 
\begin{equation*}
\begin{split}
\D \, \patt \ps 
&= \D \big(\big(1- \d \, \pat \ps\big)^{-1}\big) \, \big(\a \, \D \, \pat \ps + \b \, \D \, \ps\big) \\
&\qquad+ 2 \, \nabla \big(\big(1- \d \, \pat \ps\big)^{-1}\big) \cdot \big(\a \, \nabla \D \, \pat \ps + \b \, \nabla \D \, \ps\big) \\
&\qquad+ \big(1- \d \, \pat \ps\big)^{-1} \, \big(\a \, \D^2 \, \pat \ps + \b \, \D^2 \ps\big)\,,
\end{split}
\end{equation*}
which results from an application of the Laplace operator to the Westervelt equation and the identity $\D(f g) = \D f \, g + 2 \, \nabla f \cdot \nabla g + f \, \D g$.
Estimation by continuous Sobolev embeddings~\eqref{eq:AuxiliaryEstimates} implies 
\begin{equation*}
\begin{split}
\normbig{L^2}{\D \, \patt \ps(t)} 
&\leq 
C \, \normbig{L^{2}}{\D \big(\big(1- \d \, \pat \ps(t)\big)^{-1}\big)} \, \Big(\normbig{L^\infty}{\D \, \pat \ps(t)} + \normbig{L^\infty}{\D \, \ps(t)}\Big) \\
&\qquad+ C \, \normbig{L^4}{\nabla\big(\big(1- \d \, \pat \ps(t)\big)^{-1}\big)} \\
&\qquad\qquad\times \Big(\normbig{L^4}{\nabla \D \, \pat \ps(t)} + \normbig{L^4}{\nabla \D \, \ps(t)}\Big) \\
&\qquad+ C \, \normbig{L^{\infty}}{\big(1- \d \, \pat \ps(t)\big)^{-1}} \\
&\qquad\qquad\times \Big(\normbig{L^2}{\D^2 \, \pat \ps(t)} + \normbig{L^2}{\D^2 \ps(t)}\Big) \\
&\leq C \, \normbig{H^2}{\big(1- \d \, \pat \ps(t)\big)^{-1}} \, \Big(\normbig{H^4}{\ps(t)} + \normbig{H^4}{\pat \ps(t)}\Big) 
\end{split}
\end{equation*}
and due to $\big(\ps,\pat \ps\big) \in \nC\big([0,T], H^4(\Om) \times H^4(\Om)\big)$ proves 
\begin{equation*}
\patt \ps \in \nC\big([0,T];H^2(\Om)\big)\,.
\end{equation*}
\item 
\emph{Higher regularity.} 
In order to prove the stated higher regularity of the solution to the Westervelt equation, we proceed by induction, employing the regularity result deduced in the previous step. 
That is, we consider the partial differential equation for $\chi = \D^m \ps$ obtained from an $m-$fold application of the Laplace operator to the Westervelt equation under homogeneous Dirichlet boundary conditions 
\begin{equation}
\label{eq:PDEchi}
\begin{cases}
&\patt \chi(t) = \ahat(t) \, \D \, \pat \, \chi(t) + \bhat(t) \, \D \chi(t) + f_m(t)\,, \quad t \in (0,T]\,, \\
&\chi(t)\big\vert_{\pa \Om} = 0\,, \quad t \in [0,T]\,, \qquad \chi(0) = \D^m \ps^0\,, \quad \pat \chi(0) = \D^m \ps^1\,, 
\end{cases} 
\end{equation}
which is justified by the required compatibility conditions;  
here, we write 
\begin{equation*}
\ahat = \a \, \dhat\,, \qquad \bhat = \b \, \dhat\,, \qquad \dhat = (1 - \d \pat\ps)^{-1}\,,
\end{equation*}
for short and use the following reformulation of the right-hand side 
\begin{equation*}
\begin{split}
\D^m \big(\ahat \, \D \, \pat \ps + \bhat \, \D \, \ps\big) 
&= \D^m \big(\dhat \, \D (\a \, \pat \ps + \b \, \ps)\big) = \dhat \, \D^{m+1}(\a \, \pat \ps + \b \, \ps) + f_m \\
&= \ahat \, \D \, \pat \chi + \bhat \, \D \chi + f_m 
\end{split}
\end{equation*}
involving 
\begin{equation*}
f_m = \D^m (\dhat \, g) - \dhat \, \D^{m} g 
= \sum_{j=1}^{2m} \, \tbinom{2m}{j} \nabla^j \dhat \, \cdot \nabla^{2m-j} g\,, \qquad g = \D (\a \, \pat \ps + \b \, \ps)\big)\,. 
\end{equation*}
In order to show that 
\begin{equation*}
\big(\chi, \pat \chi) = \big(\D^m \ps, \D^m \pat \ps\big) \in \nC\big([0,T],H^2(\Om) \times H^1(\Om)\big)\,, \qquad m \in \NN_{\geq 2}\,,
\end{equation*} 
or, equivalently, 
\begin{equation*}
\big(\ps,\pat \ps\big) \in \nC\big([0,T],H^{2(m+1)}(\Om) \times H^{2m+1}(\Om)\big)\,, \qquad m \in \NN_{\geq 2}\,,
\end{equation*} 
holds, we carry out an induction proof with respect to integer $2 \leq k \leq m$ for the assertion
\begin{subequations}
\begin{equation}
\label{eq:IS}
\begin{gathered}
\big(\ps,\pat \ps\big) \in \nC\big([0,T],H^{2(k+1)}(\Om) \times H^{2k+1}(\Om)\big)\,, \\
\normbig{\nC(H^{2(k+1)} \times H^{2k+1})}{\big(\ps,\pat \ps\big)} \leq C \, \normbig{H^{2(k+1)} \times H^{2k+1}}{\big(\ps^0, \ps^1\big)}\,,
\end{gathered}
\end{equation}
based on the induction hypothesis
\begin{equation}
\label{eq:IV}
\begin{gathered}
\big(\ps,\pat \ps\big) \in \nC\big(H^{2k}(\Om) \times H^{2k-1}(\Om)\big)\,, \\
\normbig{\nC(H^{2k} \times H^{2k-1})}{\big(\ps,\pat \ps\big)} \leq C \, \normbig{H^{2k} \times H^{2k-1}}{\big(\ps^0,\ps^1\big)}\,;
\end{gathered}
\end{equation}
our previous considerations ensure that the induction beginning 
\begin{equation}
\label{eq:IB}
\begin{gathered}
\big(\ps,\pat \ps\big) \in \nC\big(H^4(\Om) \times H^4(\Om)\big)\,, \\
\normbig{\nC(H^4 \times H^4)}{\big(\ps,\pat \ps\big)} \leq C \, \normbig{H^4 \times H^4}{\big(\ps^0,\ps^1\big)}\,,
\end{gathered}
\end{equation}
\end{subequations}
is valid. 
A fundamental regularity result for linear evolution equations of hyperbolic type is given below in Proposition~\ref{prop:PropositionHyperbolic}; 
applying the first statement of Proposition~\ref{prop:PropositionHyperbolic} to \eqref{eq:PDEchi} with~$m$ replaced by~$k$, noting that by~\eqref{eq:IV}--\eqref{eq:IB} and the required smallness of the initial data the conditions 
\begin{equation*}
\begin{gathered}
a = \ahat\,, \qquad c = \bhat \geq \nu > 0\,, \\
b = 0\,, \qquad c = \bhat \in H^1\big((0,T),L^\infty(\Om)\big) \,, \qquad f = f_k \in L^2\big((0,T),L^2(\Om)\big)\,, 
\end{gathered}
\end{equation*}
are fulfilled, we obtain 
\begin{equation}
\label{eq:Est1}
\begin{split}
&\normbig{\nC([0,T],H^2 \times H^1)}{\big(\D^k \ps, \D^k \pat \ps\big)} \\
&\qquad \leq C \, \Big(\big(1 + \normbig{L^2((0,T),L^\infty)}{\pat \bhat}\big) \, \normbig{H^{2} \times H^{1}}{\big(\D^k \ps^0, \D^k \ps^1\big)} \\
&\qquad\qquad\quad + \normbig{L^2((0,T),L^2)}{f_k}\Big) \\
&\qquad \leq C \, \Big(\big(1 + \normbig{H^{4} \times H^{4}}{\big(\ps^0,\ps^1\big)}\big) \, \normbig{H^{2} \times H^{1}}{\big(\D^k \ps^0, \D^k \ps^1\big)} \\
&\qquad\qquad\quad + \normbig{L^2((0,T),L^2)}{f_k}\Big)\,,
\end{split}
\end{equation}
since the following bound holds 
\begin{equation*}
\normbig{L^2((0,T),L^\infty)}{\pat \bhat}
\leq \tfrac{\b \d}{\underline{\nu}^2} \, \normbig{L^2((0,T),L^\infty)}{\patt \ps} \leq C \, \normbig{H^{4} \times H^{4}}{\big(\ps^0,\ps^1\big)}\,.
\end{equation*}
It remains to estimate 
\begin{equation*}
\normbig{L^2((0,T),L^2)}{f_k} 
\leq \sum_{j=1}^{2k} \normbig{L^2((0,T),L^2)}{\nabla^j \dhat \cdot \nabla^{2k-j}\D(\a \pat \ps + \b \, \ps)}\,. 
\end{equation*}
For this purpose, we first consider the contribution $\dhat = \dtil(\pat\ps)$ involving the scalar function 
\begin{equation*}
\dtil: \RR \setminus \{\d^{-1}\} \longrightarrow \RR: \zeta \longmapsto \frac{1}{1 - \d \zeta}\,;
\end{equation*}
we note that non-degeneracy is essential for the regularity of~$\dtil$.
By~\eqref{eq:IV}--\eqref{eq:IB} we have $\pat \ps \in \nC([0,T],\nC^j(\overline{\Om}))$ for any $j \leq \max\{2,2k - 3\}$ and thus  
\begin{equation*}
\begin{split}
\normbig{\nC([0,T],W^{j,\infty})}{\dhat} 
&\leq \sup_{t \in [0,T]} \underset{|n| \leq j}{\sup_{\kappa \in \mathbb{N}_{\geq 0}^3}} \sup_{x \in \overline{\Om}} \absbig{\pax^{\kappa} \big(1 - \d \, \pat \ps(x,t)\big)^{-1}} \\
&\leq C \, \normbig{\nC^j([\overline{\nu}^{-1},\underline{\nu}^{-1}])}{\dtil} \, \sup_{t \in [0,T]} \normbig{\nC^{j}(\overline{\Om})}{\pat \ps(t)} \\
&\leq C \, \normbig{\nC^j([\overline{\nu}^{-1},\underline{\nu}^{-1}])}{\dtil} \, \normbig{\nC([0,T],H^{j+2})}{\pat\ps}\,;
\end{split}
\end{equation*}
the highest derivatives of $\dhat$ can be estimated utilising $\nabla \dhat =\d \, \dhat^{\,2} \nabla \pat \ps$ and 
\begin{equation*}
\begin{split}
\nabla^j \dhat &= \d \, \nabla^{j-1} \big(\dhat^{\,2} \nabla \pat \ps\big) 
= \d \sum_{i=0}^{j-1} \tbinom{j-1}{i} \, \nabla^i \dhat^{\,2} \, \nabla^{j-i} \, \pat \ps \\
&= \d \sum_{i=0}^{j-1} \tbinom{j-1}{i} \sum_{\ell=0}^{i} \tbinom{i}{\ell} \nabla^{\ell} \dhat \, \nabla^{i-\ell} \dhat \, \nabla^{j-i} \pat \ps\,.
\end{split}
\end{equation*}
For $j=2k-2$, $j=2k-1$, and $j=2k$ we get
\begin{equation*}
\begin{split}
\nabla^{2k-2} \dhat &= \underbrace{\d \sum_{i=1}^{2k-3} \tbinom{2k-3}{i} \, \nabla^i \dhat^{\,2} \, \nabla^{2k-2-i} \, \pat \ps}_{G_1} 
 + \underbrace{\d \, \dhat^{\,2}}_{G_1} \underbrace{\nabla^{2k-2} \, \pat \ps}_{G_2}\,, \\
\nabla^{2k-1} \dhat 
&= \underbrace{\d \sum_{i=2}^{2k-3} \tbinom{2k-2}{i} \, \nabla^i \dhat^{\,2} \, \nabla^{2k-1-i} \, \pat \ps}_{G_1} 
+ \underbrace{\d \, \dhat^{\,2}}_{G_1} \underbrace{\nabla^{2k-1} \, \pat \ps}_{G_3} \\
&\qquad + \underbrace{\d \, (2k-2) \nabla \, \dhat^{\,2}}_{G_1} \underbrace{\nabla^{2k-2} \, \pat \ps}_{G_2} 
+ \underbrace{\d \, \nabla^{2k-2} \, \dhat^{\,2}}_{G_2} \underbrace{\nabla \pat \ps}_{G_1}\,, \\
\nabla^{2k} \dhat  
&= \underbrace{\d \sum_{i=3}^{2k-3} \tbinom{2k-1}{i} \nabla^i \, \dhat^{\,2} \, \nabla^{2k-i} \, \pat \ps}_{G_1} 
+ \underbrace{\d \, \dhat^{\,2}}_{G_1} \underbrace{\nabla^{2k} \, \pat \ps}_{G_4} \\
&\qquad + \underbrace{\d \, (2k-1) \, \nabla \dhat^{\,2}}_{G_1} \underbrace{\nabla^{2k-1} \, \pat \ps}_{G_3}\, 
+ \underbrace{\d \, \tbinom{2k-1}{2} \nabla^{2} \, \dhat^{\,2}}_{G_1} \underbrace{\nabla^{2k-2} \, \pat \ps}_{G_2} \\
&\qquad+ \underbrace{\d \, (2k-1) \nabla^{2k-2} \, \dhat^{\,2}}_{G_2} \underbrace{\nabla^{2} \, \pat \ps}_{G_1}
\, + \underbrace{\d \, \nabla^{2k-1} \, \dhat^{\,2}}_{G_3} \underbrace{\nabla \pat \ps}_{G_1}\,, 
\end{split}
\end{equation*}
involving quantities that have in common that they satisfy 
\begin{equation*}
\begin{gathered}
G_1 \in \nC\big([0,T],L^\infty(\Om)\big)\,, \qquad G_2 \in \nC\big([0,T],L^6(\Om)\big)\,, \\ 
G_3 \in \nC\big([0,T],L^2(\Om)\big)\,, \qquad G_4 \in L^2\big((0,T),L^2(\Om)\big)\,, 
\end{gathered}
\end{equation*}
with $G_1,G_2,G_3$ bounded provided that~$\pat \ps$ is bounded in $\nC([0,T],H^{\max\{4,2k-1\}}(\Om))$ and~$G_4$ bounded  whenever~$\D^k \pat \ps$ is bounded in~$L^2((0,T),L^2(\Om))$. 
Using once more~\eqref{eq:IV}--\eqref{eq:IB} implies 
\begin{equation*}
\normbig{\nC([0,T],H^{\max\{4,2k-1\}})}{\pat\ps} \leq C \, \norm{H^{2k} \times H^{\max\{4,2k-1\}}}{\big(\ps^0,\ps^1\big)}\,, 
\end{equation*}
and thus we obtain the estimate 
\begin{equation*}
\begin{split}
\norm{L^2((0,T),L^2)}{f_k} 
&\leq C \, \normbig{H^{2k} \times H^{\max\{4,2k-1\}}}{\big(\ps^0,\ps^1\big)} \, \big(1+ \normbig{L^2((0,T),L^2}{\D^k \pat\ps}\big) \\
&\leq C \, \normbig{H^{2k} \times H^{\max\{4,2k-1\}}}{\big(\ps^0,\ps^1\big)} \, \big(1+ \normbig{\nC([0,T],H^1}{\D^k \pat \ps}\big)\,. 
\end{split}
\end{equation*}
Inserting into~\eqref{eq:Est1}, this finally implies the relation 
\begin{equation*}
\begin{split}
&\normbig{\nC([0,T],H^2 \times H^1)}{\big(\D^k \ps, \D^k \pat \ps\big)} \\
&\qquad \leq C \, \Big(\big(1 + \normbig{H^{4} \times H^{4}}{\big(\ps^0,\ps^1\big)}\big) \, \normbig{H^{2} \times H^{1}}{\big(\D^k \ps^0, \D^k \ps^1\big)} \\
&\qquad\qquad + \normbig{H^{2k} \times H^{\max\{4,2k-1\}}}{\big(\ps^0,\ps^1\big)} \, \big(1+ \normbig{\nC([0,T],H^1}{\D^k \pat \ps}\big)\Big)\,;
\end{split}
\end{equation*}
by restricting the size of the inital data the condition 
\begin{equation*}
\begin{split}
\normbig{H^{2k} \times H^{\max\{4,2k-1\}})}{\big(\ps^0,\ps^1\big)} < \tfrac{1}{2C}
\end{split}
\end{equation*}
can be fulfilled. 
Altogether, this implies~\eqref{eq:IS} and concludes the proof. 
$\hfill \diamondsuit$
\end{enumerate}
\subsection{Linear evolution equations of hyperbolic type}
\paragraph{Regularity result.} 
In the following, we state a regularity result for a linear wave equation subject to homogeneous Dirichlet boundary conditions on a regular domain  
\begin{equation} 
\label{eq:LinearHyperbolic}
\begin{cases}
&\patt \ps(x,t) = a(x,t) \, \D \, \pat \ps(x,t) + b(x,t) \, \pat \ps(x,t) \\
&\qquad\qquad\qquad +\; c(x,t) \, \D \, \ps(x,t) + f(x,t)\,, \quad (x,t) \in \Om \times (0,T]\,, \\
&\ps(x,t) = 0\,, \quad (x,t) \in \pa \Om \times [0,T]\,, \\
&\ps(x,0) = \ps^0(x)\,, \quad \pat \ps(x,0) = \ps^1(x)\,, \quad x \in \Om\,.
\end{cases} 
\end{equation}
We note that in the present context the coefficient~$b$ may have arbitrary sign and thus the term $b \, \pat \ps$ cannot be regarded as a damping term. 
\begin{proposition}
\label{prop:PropositionHyperbolic}
Assume that the coefficient functions satisfy the positivity and regularity conditions 
\begin{equation*} 
\begin{gathered}
a(x,t)\,, \, c(x,t)\geq \nu>0\,, \qquad (x,t) \in \Om \times [0,T]\,, \\
b \in L^2\big((0,T),L^3(\Om)\big)\,, \quad c \in W^{1,1}\big((0,T),L^\infty(\Om)\big)\,, \quad f \in L^2\big((0,T),L^2(\Om)\big)\,,
\end{gathered}
\end{equation*}
with $\norm{L^2((0,T),L^3)}{b}$ and $\norm{L^2((0,T),L^\infty)}{\pat c}$ sufficiently small.
\begin{enumerate}[(i)]
\item 
Then the solution to~\eqref{eq:LinearHyperbolic} satisfies the relation 
\begin{equation*} 
\big(\ps, \pat \ps\big) \in \nC\big([0,T],H^{2}(\Om) \times H^1(\Om)\big) 
\end{equation*}
and the bound 
\begin{equation*} 
\normbig{\nC([0,T],H^2 \times H^1)}{\big(\ps, \pat \ps\big)} 
\leq C \, \big(\normbig{H^{2} \times H^{1}}{\big(\ps^0,\ps^1\big)} + \normbig{L^2((0,T),L^2)}{f}\big) 
\end{equation*}
with constant depending on the respective norms of~$b$ and~$c$, but not on~$a$.
\item 
If in addition the relations 
\begin{equation*}
\begin{gathered}
a \in L^2\big((0,T),W^{1,\infty}(\Om) \cap W^{2,3}(\Om)\big)\,, \quad b \in L^2\big((0,T), H^2(\Om)\big)\,, \\
c \in L^2\big((0,T),H^2(\Om)\big) \cap H^1\big((0,T),L^\infty(\Om)\big)\,, \\
f \in \nC\big([0,T],H_0^1(\Om)\big) \cap L^2\big((0,T),H^2(\Om)\big)\,, \qquad \ps^1\in H_0^1(\Om)\,, \\
\patt \ps(\cdot,0) = a(\cdot,0) \, \D \ps^1 + b(\cdot,0) \, \ps^1 + c(\cdot,0) \, \D \ps^0 + f(\cdot,0) \in H_0^1(\Om)\,,
\end{gathered}
\end{equation*}
are satisfied and the norms 
\begin{equation*}
\normbig{L^2((0,T),W^{1,\infty})}{a}\,, \qquad 
\normbig{L^2((0,T),H^1)}{b}\,, \qquad 
\normbig{L^2((0,T),H^2)}{c}\,, 
\end{equation*}
are sufficiently small, then the solution to~\eqref{eq:LinearHyperbolic} satisfies the relation 
\begin{equation*} 
\big(\ps, \pat \ps\big) \in \nC\big([0,T],H^4(\Om) \times H^3(\Om)\big)
\end{equation*}
and the bound 
\begin{equation*} 
\normbig{\nC([0,T],H^4 \times H^3)}{\big(\ps, \pat \ps\big)} 
\leq C \, \big(\normbig{H^4 \times H^3}{\big(\ps^0,\ps^1\big)} + \normbig{L^2((0,T),H^2)}{f}\big)
\end{equation*}
with constant depending on the respective norms of $a,b,c$. 
If additionally 
\begin{equation*}
b \in L^2\big((0,T), W^{1,\infty}(\Om) \cap W^{2,3}(\Om)\big)\,, 
\end{equation*}
then smallness of $\norm{L^2((0,T),H^1)}{b}$ is not required. 
\end{enumerate}
\end{proposition} 
\paragraph{Proof.}
As before, for a function $f: \overline{\Om} \times [0,T] \to \RR$ we write $f(t) = f(\cdot,t)$ for short. 
\begin{enumerate}[(i)]
\item
Multiplication of the partial differential equation by~$\D \, \pat \ps$, integration with respect to space and time, and application of integration-by-parts yields
\begin{equation*}
\begin{split}
&\tfrac{1}{2} \, \normbig{L^2}{\nabla \pat \ps(t)}^2 + \int_0^t \normbig{L^2}{\sqrt{a(\t)} \, \D \, \pat \ps(\t)}^2 \; \dd\t 
+ \tfrac{1}{2} \, \normbig{L^2}{\sqrt{c} \, \D \, \ps(t)}^2 \\
&\qquad = \tfrac{1}{2} \, \normbig{L^2}{\nabla \ps^1}^2 + \tfrac{1}{2} \, \normbig{L^2}{\sqrt{c} \, \D \ps^0}^2 \\
&\qquad\qquad+ \int_0^t \int_\Om \Big(\tfrac{1}{2} \, \pat c(\t) \, \big(\D \, \ps(\t)\big)^2 - \big(f(\t) + b(\t) \, \pat \ps(\t)\big) \, \D \, \pat \ps(\t)\Big) \; \dd\t\,. 
\end{split}
\end{equation*}
Estimation by Hölder's and Young's inequality implies on the one hand 
\begin{equation*}
\begin{split}
&\int_0^t \int_\Om \tfrac{1}{2} \, \pat c(\t) \, \big(\D \, \ps(\t)\big)^2  \; \dd\t 
= \int_0^t \int_\Om \tfrac{1}{2} \, \frac{\pat c(\t)}{c(\t)} \, \big(\sqrt{c(\t)} \, \D \, \ps(\t)\big)^2  \; \dd\t \\
&\qquad\leq \nu^{-1} \, \normbig{L^1((0,T),L^\infty)}{\pat c} \, \sup_{\t \in [0,t]} \normbig{L^2}{\sqrt{c} \, \D \, \ps(\t)}^2 
\end{split}
\end{equation*}
and on the other hand
\begin{equation*}
\begin{split}
&-\int_0^t \int_\Om \big(f(\t) + b(\t) \, \pat \ps(\t)\big) \, \D \, \pat \ps(\t) \; \dd\t \\
&\qquad\leq 
\normBig{L^2((0,t),L^2)}{\frac{f + b \, \pat \ps}{\sqrt{a}}} \, \normbig{L^2((0,t),L^2)}{\sqrt{a} \, \D \, \pat \ps} \\
&\qquad\leq (2\nu)^{-1} \, \normbig{L^2((0,T),L^2)}{f} ^2 + (2\nu)^{-1} \, \normbig{L^2((0,T),L^3)}{b}^2 \, \sup_{\t \in [0,t]} \normbig{L^6}{\pat \ps(\t)}^2 \\
&\qquad\qquad+\normbig{L^2((0,t),L^2)}{\sqrt{a} \, \D \, \pat \ps}^2\,, 
\end{split}
\end{equation*}
since $(x_1 + x_2) \, x_3 = x_1 x_3 + x_2 x_3 \leq \frac{1}{2} x_1^2 + \frac{1}{2} x_2^2 + x_3^2$ for $x_1, x_2, x_3 \in \RR$, involving a constant that depends on~$\nu^{-1}$;
hence, using that the term involving~$\sqrt{a} \, \D \, \pat \ps$ cancels out, we obtain  
\begin{equation*}
\begin{split}
&\tfrac{1}{2} \, \normbig{L^2}{\nabla \pat \ps(t)}^2 + \tfrac{1}{2} \, \normbig{L^2}{\sqrt{c} \, \D \, \ps(t)}^2 \\
&\quad\leq \tfrac{1}{2} \, \normbig{L^2}{\nabla \ps^1}^2 + \tfrac{1}{2} \, \normbig{L^2}{\sqrt{c} \, \D \ps^0}^2 \\
&\quad\quad+ \nu^{-1} \, \normbig{L^1((0,T),L^\infty)}{\pat c} \, \sup_{\t \in [0,t]} \normbig{L^2}{\sqrt{c} \, \D \, \ps(\t)}^2 \\
&\quad\quad+ (2\nu)^{-1} \, \normbig{L^2((0,T),L^2)}{f}^2 
+ (2\nu)^{-1} \, \normbig{L^2((0,T),L^3)}{b}^2 \, \sup_{\t \in [0,t]} \, \normbig{L^2}{\nabla \pat \ps(\t)}^2\,, 
\end{split}
\end{equation*}
see also~\eqref{eq:AuxiliaryEstimates}.
Since $\norm{L^2((0,T),L^3)}{b}$ and $\norm{L^1((0,T),L^\infty)}{\pat c}$ are  assumed to be sufficiently small, taking the supremum over time, for a sufficiently small positive constant $\g > 0$ that is independent of the final time~$T$, 
we obtain 
\begin{equation*}
\begin{split}
&\g \, \Big(\normbig{\nC([0,T],H^2)}{\ps}^2 + \normbig{\nC([0,T],H^1)}{\pat \ps}^2\Big) \\
&\qquad\leq 
\big(1 - (2 \nu)^{-1} \, \normbig{L^2((0,T),L^3)}{b}^2\big) \, \sup_{\t \in [0,t]} \, \normbig{L^2}{\nabla \pat \ps(\t)}^2 \\
&\qquad\qquad+ \big(1 - \nu^{-1} \, \normbig{L^1((0,T),L^\infty)}{\pat c}\big) \, \sup_{\t \in [0,t]} \normbig{L^2}{\sqrt{c} \, \D \, \ps(\t)}^2 \\
&\qquad\leq \normbig{L^2}{\nabla \ps^1}^2 + \normbig{L^2}{\sqrt{c} \, \D \ps^0}^2  + (2 \nu)^{-1} \, \normbig{L^2((0,T),L^2)}{f}^2\,,
\end{split}
\end{equation*}
which proves the statement. 
\item 
Under the assumed additional regularity and compatibility conditions it is justified to apply the Laplace operator to both sides of~\eqref{eq:LinearHyperbolic}. 
The function $\pstil = \D \, \ps$ then solves an initial boundary value problem of a similar form 
\begin{equation*}
\begin{cases}
&\patt \pstil(t) = a(t) \, \D \, \pat \pstil(t) + \bbtil(t) \, \pat \pstil(t) + c(t) \, \D \pstil(x,t) + \ftil(t)\,, \quad t \in (0,T]\,, \\
&\pstil(t)\big\vert_{\pa \Om} = 0\,, \quad t \in [0,T]\,, \qquad \pstil(0) = \D \ps^0\,.
\end{cases} 
\end{equation*}
involving the coefficient functions $\bbtil = b + \D a$ and 
\begin{equation*}
\begin{split}
\ftil 
&= \D f + \D c \, \D \ps + \D b \, \pat \ps + 2 \, \big(\nabla c \cdot \nabla \D \ps + \nabla b \cdot \nabla \pat \ps + \nabla a \cdot \nabla \D \, \pat \ps\big) \\
&= \D f + \D c \, \pstil + \D b \, \pat \ps + 2 \, \big(\nabla c \cdot \nabla \pstil + \nabla b \cdot \nabla \pat \ps + \nabla a \cdot \nabla \pat \pstil\big)\,.
\end{split}
\end{equation*}
We utilise the bound
\begin{equation*}
\normbig{L^2((0,T),L^3)}{\tilde{b}} \leq \normbig{L^2((0,T),L^3)}{b} + \normbig{L^2((0,T),W^{2,3})}{a}
\end{equation*}
as well as 
\begin{equation}\label{eq:estftil}
\begin{split}
&\normbig{L^2((0,T),L^2)}{\tilde{f}}\\
&\qquad\leq \normbig{L^2((0,T),H^2)}{f} 
+ \normbig{L^2((0,T),L^2)}{\D c \, \pstil}
+ \normbig{L^2((0,T),L^2)}{\D b \, \pat \ps} \\
&\qquad\qquad+ 2 \, \normbig{L^2((0,T),L^2)}{\nabla c \cdot \nabla \pstil} + 2 \, \normbig{L^2((0,T),L^2)}{\nabla b \cdot \nabla \pat \ps} \\
&\qquad\qquad+ 2 \, \normbig{L^2((0,T),L^2)}{\nabla a \cdot \nabla \pat \pstil} \\
&\qquad\leq \normbig{L^2((0,T),H^2)}{f} + \normbig{L^2((0,T),L^2)}{\D c} \, \normbig{\nC((0,T),L^\infty)}{\pstil} \\
&\qquad\qquad+ \normbig{L^2((0,T),L^2)}{\D b} \, \normbig{\nC((0,T),L^\infty)}{\pat \ps} \\
&\qquad\qquad+ 2 \, \normbig{L^2((0,T),L^3)}{\nabla c} \, \normbig{\nC((0,T),L^6)}{\nabla \pstil} \\
&\qquad\qquad+ 2 \, \normbig{L^2((0,T),L^3)}{\nabla b} \, \normbig{\nC((0,T),L^6)}{\nabla \pat \ps} \\
&\qquad\qquad+ 2 \, \normbig{L^2((0,T),L^{\infty})}{\nabla a} \, \normbig{\nC((0,T),L^2)}{\nabla \pat \pstil} \\
&\qquad\leq \normbig{L^2((0,T),H^2)}{f} \\
&\qquad\qquad+ C\normbig{L^2((0,T),H^2)}{b}  \, \normbig{\nC((0,T),L^2)}{\pat \pstil} \\
&\qquad\qquad + C \, \normbig{L^2((0,T),H^2)}{c} \, \normbig{\nC((0,T),H^2)}{\pstil} \\
&\qquad\qquad+ 2 \, \normbig{L^2((0,T),W^{1,\infty})}{a} \, \normbig{\nC((0,T),H^1)}{\pat \pstil}\,, 
\end{split}
\end{equation}
obtained by means of Hölder's inequality and continuous Sobolev embeddings, see also~\eqref{eq:AuxiliaryEstimates}.   
By assumption, the coefficients $a,b,c$ are sufficiently small; 
together with the first assertion this proves the stated result.  
If additionally $b \in L^2((0,T), W^{1,\infty}(\Om) \cap W^{2,3}(\Om))$, in~\eqref{eq:estftil} we instead estimate the terms involving~$b$ by 
\begin{equation*}
\begin{split}
&\normbig{L^2((0,T),L^2)}{\D b \, \pat \ps} 
+ 2 \, \normbig{L^2((0,T),L^2)}{\nabla b \cdot \nabla \pat \ps} \\
&\qquad\leq 
\normbig{L^2((0,T),L^3)}{\D b} \, \normbig{\nC((0,T),L^6)}{\pat \ps} \\
&\qquad\qquad+ 2 \, \normbig{L^2((0,T),L^\infty)}{\nabla b} \, \normbig{\nC((0,T),L^2)}{\nabla \pat \ps} \\
&\qquad\leq C \, \big(\normbig{L^2((0,T),W^{2,3})}{b} 
+ \normbig{L^2((0,T),W^{1,\infty})}{b}\big) \, \normbig{\nC((0,T),H^1)}{\pat \ps}\,;
\end{split}
\end{equation*}
applying statement~(i) to $\norm{\nC((0,T),H^1)}{\pat \ps}$ shows that the smallness assumption on~$b$ is not needed.
$\hfill \diamondsuit$
\end{enumerate}
\subsection{Linear evolution equations of parabolic type}
\paragraph{Regularity result.} 
In the following, we state a regularity result for a linear reaction-diffusion equation subject to homogeneous Dirichlet boundary conditions on a regular domain  
\begin{equation} 
\label{eq:LinearParabolic}
\begin{cases}
&\pat \ps(x,t) \\
&\quad= a(x,t) \, \D \ps(x,t) + b(x,t) \, \ps(x,t) + f(x,t)\,, \quad (x,t) \in \Om \times (0,T]\,, \\
&\ps(x,t) = 0\,, \quad (x,t) \in \pa \Om \times [0,T]\,, \\
&\ps(x,0) = \ps^0(x)\,, \quad x \in \Om\,.
\end{cases} 
\end{equation}

\begin{proposition}
\label{prop:PropositionParabolic}
Assume that the coefficient functions satisfy the positivity and regularity conditions 
\begin{equation*} 
\begin{gathered}
0 < \nu \leq a(x,t) \leq \mu < \infty\,, \quad (x,t) \in \Om \times [0,T]\,, \\
a \in \nC\big([0,T],W^{1,3}(\Om)\big)\,, \quad b \in L^2\big((0,T),H^1(\Om)\big)\,, \quad f \in L^2((0,T),H^1(\Om))\,,
\end{gathered}
\end{equation*}
and that the norms of~$a$ and~$b$ are sufficiently small. 
\begin{enumerate}[(i)]
\item
Then the solution to~\eqref{eq:LinearParabolic} satisfies the relation 
\begin{equation*} 
\ps \in \nC\big([0,T],H^{2}(\Om)\big)
\end{equation*}
and the bound 
\begin{equation*} 
\normbig{\nC([0,T],H^2)}{\ps} + \normbig{L^2((0,T),H^1)}{\pat \ps} 
\leq C \, \big(\normbig{H^{2}}{\ps^0} + \normbig{L^2((0,T),H^1)}{f}\big)
\end{equation*}
with constant depending on the respective norms of~$a$ and~$b$.
\item 
If in addition the relations 
\begin{equation*}
\begin{gathered}
a \in L^2\big((0,T),H^3(\Om)\big)\,, \quad b \in L^2\big((0,T),H^3(\Om)\big)\,, \\
f \in \nC\big([0,T],H_0^1(\Om)\big) \cap L^2\big((0,T),H^3(\Om)\big)\,, \\
\pat \ps(\cdot,0) = a(\cdot,0) \, \D \ps^0 + b(\cdot,0) \ps^0 + f^0(\cdot,0) \in H_0^1(\Om)\,,
\end{gathered}
\end{equation*}
hold and the norms 
\begin{equation*}
\normbig{L^2((0,T),W^{1,\infty})}{a}\,, \quad \normbig{L^2((0,T),W^{2,3})}{a}\,, \quad \normbig{L^2((0,T),H^2)}{b}\,,
\end{equation*}
are sufficiently small, the solution to~\eqref{eq:LinearParabolic} satisfies the relation 
\begin{equation*} 
\ps \in \nC\big([0,T],H^4(\Om)\big)
\end{equation*}
and the bound 
\begin{equation*} 
\normbig{\nC([0,T],H^4)}{\ps} + \normbig{L^2((0,T),H^3)}{\pat \ps} \leq C \, \big(\normbig{H^{4}}{\ps^0} + \normbig{L^2((0,T),H^3)}{f}\big)
\end{equation*}
with constant depending on the respective norms of~$a$ and~$b$.
\item 
If for some $m \in \NN_{\geq 2}$ the regularity and compatibility conditions 
\begin{equation*}
\begin{gathered}
a \in L^2\big((0,T),H^5(\Om)\big)\,, \\
\quad f^{m-1} \in \nC\big([0,T],H_0^1(\Om)\big) \cap L^2\big((0,T),H^{2m+1}(\Om)\big)\,, \\
\D^j\ps(\cdot,0)\in H_0^1(\Om)\,, \\
\pat \D^j \ps(\cdot,0) = a(\cdot,0) \, \D^{j+1} \ps^0 + b^j(\cdot,0) \D^j \ps^0 + f^j(\cdot,0) \in H_0^1(\Om)\,, \\
\end{gathered}
\end{equation*}
are fulfilled for any integer $1 \leq j \leq m$, where $b^j = b^{j-1} + \D a$ and 
\begin{equation*}
f^j=\D f^{j-1} + \D b^{j-1}\D^{j-1}\ps + 2\nabla b^{j-1} \nabla \D^{j-1} \ps + 2 \nabla a\nabla \D^j \ps
\end{equation*}
for $1 \leq j \leq m-1$, and if the norm of~$a$ is sufficiently small, then the solution to~\eqref{eq:LinearParabolic} satisfies the relation 
\begin{equation*}
\ps \in \nC\big([0,T],H^{2(m+1)}(\Om)\big)
\end{equation*}
and the bound 
\begin{equation*} 
\begin{split}
&\normbig{\nC([0,T],H^{2(m+1)})}{\ps} + \normbig{L^2((0,T),H^{2m+1})}{\pat \ps} \\
&\quad \leq C \, \big(\normbig{H^{2(m+1)}}{\ps^0} + \normbig{L^2((0,T),H^{2m+1})}{f}\big)
\end{split}
\end{equation*}
with constant depending on the respective norms of~$a$ and~$b$.
\end{enumerate}
\end{proposition} 
\paragraph{Proof.}
The derivation of the first statement is in the lines of the proof of Proposition~\ref{prop:PropositionHyperbolic}, testing the partial differential equation with~$\frac{1}{a} \, \D \, \pat \ps$. 
Applying the obtained regularity result to the function $\D^k \ps$, which solves a related initial boundary value problem, then yields the remaining assertions. 
We omit the technical details.  
$\hfill \diamondsuit$
\section{Conclusions}
In this work, we have introduced and investigated operator splitting methods for the efficient time integration of the Westervelt equation modelling the propagation of high intensity ultrasound in nonlinear acoustics.
We have provided numerical comparisons for the first-order Lie-Trotter and second-order Strang splitting methods based on four decompositions, using explicit and implicit solvers for the numerical solution of the subproblems. 
The numerical examples confirm that time-splitting methods remain stable and retain their nonstiff orders of convergence for sufficiently regular problem data. 
For the the Lie-Trotter splitting method based on the computationally most favourable decomposition, we have carried out a rigorous stability and error analysis. 

Future work shall be concerned with an extension of the error analysis to the second-order Strang splitting method, justifying the use of an adaptive time stepsize control combining the first-order Lie--Trotter and second-order Strang splitting methods. 
Also, it remains to analyse the effect of additional errors caused by the numerical solution of the subproblems. 
Furthermore, it is of interest to study absorbing boundary conditions used for tackling unbounded domains or the excitation by Neumann boudary conditions and more advanced models of nonlinear acoustics such as Kuznetsov's equation. 
\appendix 
\section{Background on a compact local error expansion} 
\label{sec:LocalError}
A basic ingredient in the convergence analysis of time-splitting methods~\eqref{eq:Splitting} applied to nonlinear evolution equations of the form~\eqref{eq:ProblemAB} are suitable expansions for the local error reflecting the expected dependence on the time stepsize;  
due to the presence of unbounded nonlinear operators it is thereby essential to specify the required regularity assumptions on the initial state. 
For the sake of completeness, we detail the approach for the first-order Lie--Trotter splitting method exploited in~\cite{AuzingerHofstaetterKochThalhammer2013Nonlinear,DescombesThalhammer2012} in the context of nonlinear Schrödinger equations.  
For the convenience of the reader, we first collect auxiliary definitions and results, see also~\cite{AuzingerHofstaetterKochThalhammer2013Nonlinear,DescombesThalhammer2012,Lunardi}. 
\subsection{Prerequisites} 
We denote by $G: D(G) \subset X \to X$ as well as $H: D(H) \subset X \to X$ unbounded nonlinear operators.

\paragraph{Evolution operators.}
Accordingly to~\eqref{eq:EF} we denote by~$\nEG$ the evolution operator associated with the initial value problem 
\begin{equation}
\label{eq:EG}
\begin{cases}
&\tfrac{\dd}{\dd t} \, \nE_G(t,v) = G\big(\nE_G(t,v)\big)\,, \quad t \in (0,T]\,,\\ 
&\nE_G(t,v) = v\,.
\end{cases}
\end{equation}
More generally, in order to capture the additional dependence of the solution to a non-autonomous evolution equation on the initial time, we employ the notation 
\begin{equation*}
\begin{cases}
&\tfrac{\dd}{\dd t} \, \nE_G(t,t_0,v) = G\big(t,\nE_G(t,t_0,v)\big)\,, \quad t \in (t_0,T]\,, \\
&\nE_G(t_0,t_0,v) = v\,.
\end{cases}
\end{equation*}

\paragraph{Derivative with respect to initial value.}
The derivative of the evolution operator~$\nEG$ with respect to the initial value~$v$, denoted by~$\pa_{2} \nEG$, satisfies the non-autonomous linear initial value problem (variational  equation) 
\begin{equation}
\label{eq:dEG}
\begin{cases}
&\tfrac{\dd}{\dd t} \, \pa_{2} \nE_G(t,v) = G'\big(\nE_G(t,v)\big) \, \pa_{2} \nE_G(t,v)\,, \quad t \in (0,T]\,, \\
&\pa_{2} \nE_G(0,v) = I\,,
\end{cases}
\end{equation}
with right-hand side involving the Fréchet derivative of~$G$.

\paragraph{Inverse operator.}
Provided that the evolution operator associated with~\eqref{eq:EG} is well-defined for negative times, the identities 
\begin{equation*}
\nEG\big(- t,\nEG(t,v)\big) = v\,, \qquad \nEG\big(t,\nEG(- t,w)\big) = w\,, 
\end{equation*}
are valid. 
Differentiation with respect to~$v$ or~$w$, respectively, thus yields 
\begin{equation*}
\pa_{2} \nEG\big(- t,\nEG(t,v)\big) \, \pa_{2} \nEG(t,v) = I\,, \qquad \pa_{2} \nEG\big(t,\nEG(- t,w)\big) \, \pa_{2} \nEG(- t,w) = I\,. 
\end{equation*}
Setting $w = \nEG(t,v)$ the relation
\begin{equation*}
\pa_{2} \nEG\big(t,\nEG\big(- t,\nEG(t,v)\big)\big) \, \pa_{2} \nEG\big(- t,\nEG(t,v)\big) 
= \pa_{2} \nEG(t,v) \, \pa_{2} \nEG\big(- t,\nEG(t,v)\big) = I
\end{equation*}
follows. Altogether, this implies 
\begin{equation}
\label{eq:Inverse}
\big(\pa_{2} \nEG(t,v)\big)^{-1} = \pa_{2} \nEG\big(- t,\nEG(t,v)\big)\,.
\end{equation}

\paragraph{Linear variation-of-constants formula.}
A useful representation for the solution to an inital value problem involving a time-dependent linear operator and an additional inhomogeneity
\begin{subequations}
\label{eq:Cr}
\begin{equation}
\label{eq:CrEquation}
\begin{cases}
&\ddt \, \nE_{C+r}(t,t_0,v) = C(t) \, \nE_{C+r}(t,t_0,v) + r(t)\,, \quad t \in (t_0,T]\,, \\
&\nE_{C+r}(t_0,t_0,v) = v\,, 
\end{cases}
\end{equation}
relies on the linear variation-of-constants formula
\begin{equation}
\label{eq:CrSolution}
\nE_{C+r}(t,t_0,v) = \nE_{C}(t,t_0) \, v + \nE_{C}(t,t_0) \int_{t_0}^{t} \big(\nE_{C}(\t,t_0)\big)^{-1} r(\t) \; \dd\t\,; 
\end{equation}
here, we denote by $\nE_{C}(t,t_0)$ the resolvent (fundamental system) with the properties 
\begin{equation*}
\ddt \, \nE_{C}(t,t_0) = C(t) \, \nE_{C}(t,t_0)\,, \qquad \nE_{C}(t_0,t_0) = I\,.
\end{equation*}
We note that due to the linearity of the problem the evolution operator associated with the homogeneous problem is independent of the initial value. 
The above solution representation is verified by a brief calculation 
\begin{equation*}
\begin{split}
\ddt \, \nE_{C+r}(t,t_0,v) 
&= \ddt \, \nE_{C}(t,t_0) \, \bigg(v + \int_{t_0}^{t} \big(\nE_{C}(\t,t_0)\big)^{-1} r(\t) \; \dd\t\bigg) \\
&\qquad+ \nE_{C}(t,t_0) \, \big(\nE_{C}(t,t_0)\big)^{-1} r(t) \\
&= C(t) \, \bigg(\nE_{C}(t,t_0) \, v + \nE_{C}(t,t_0) \int_{t_0}^{t} \big(\nE_{C}(\t,t_0)\big)^{-1} r(\t) \; \dd\t\bigg) +  r(t) \\
&= C(t) \, \nE_{C+r}(t,t_0,v) +  r(t)\,.
\end{split}
\end{equation*}

\paragraph{Extension to parabolic equations.}
In situations where the linear non-autonomous evolution equation~\eqref{eq:CrEquation} is non-reversible in time, the inverse operator $(\nE_{C}(\t,t_0))^{-1}$ and thus the solution representation~\eqref{eq:CrSolution} is not defined, in general. 
However, for evolution equations of parabolic type the results deduced in~\cite[Sec.\,6]{Lunardi} ensure the existence of a family of bounded linear operators such that 
\begin{equation*}
\begin{gathered}
R(t,s) \, R(s,t_0) = R(t,t_0)\,, \quad R(s,s) = I\,, \\ 
\ddt R(t,s) = C(t) \, R(t,s)\,, \qquad 0 \leq t_0 \leq s \leq t \leq T\,, 
\end{gathered}
\end{equation*}
which leads to solution representation 
\begin{equation}
\nE_{C+r}(t,t_0,v) = R(t,t_0) \, v + \int_{t_0}^{t} R(t,\t) \, r(\t) \; \dd\t\,.
\end{equation}
\end{subequations}

\paragraph{Fundamental identity.}
A fundamental identity relates the defining function, the associated evolution operator and its derivative with respect to the initial value 
\begin{equation}
\label{eq:FundamentalIdentity}
\pa_{2} \nE_G(t,v) \, G(v) = G\big(\nE_G(t,v)\big)\,. 
\end{equation}
The identity is obtained by means of relation~\eqref{eq:EG} and the chain rule
\begin{equation*}
\begin{split}
\pa_{2} \nE_G(t,v) \, G(v) 
&= \pa_{2} \nE_G\big(t,\nE_G(0,v)\big) \, G\big(\nE_G(0,v)\big) = \tfrac{\dd}{\dd s}\Big\vert_{s=0} \, \nE_G\big(t,\nE_G(s,v)\big) \\
&= \tfrac{\dd}{\dd s}\Big\vert_{s=0} \, \nE_G(t+s,v) = G\big(\nE_G(t,v)\big)\,.
\end{split}
\end{equation*}

\paragraph{Lie-commutator.}
In accordance with the commutator of linear operators, the first Lie-commutator is defined by 
\begin{equation}
\label{eq:Commutator}
\big[G,H\big] \, v = G'(v) \, H(v) - H'(v) \, G(v)\,.
\end{equation}

\paragraph{A generalisation of the fundamental identity.}
An essential tool in the derivation of suitable local error expansions for time-splitting methods applied to nonlinear evolution equations is the relation  
\begin{subequations}
\label{eq:FundamentalIdentityGeneralised}
\begin{equation}
\begin{split}
\pa_{2} \nE_G(t,v) \, H(v) 
&= H\big(\nE_G(t,v)\big) \\
&\qquad+ \pa_{2} \nE_G(t,v) \int_{0}^{t} \big(\pa_{2} \nE_G(\t,v)\big)^{-1} \, \big[G,H\big]\big(\nE_G(\t,v)\big) \; \dd \t
\end{split}
\end{equation}
generalising~\eqref{eq:FundamentalIdentity}. 
In order to deduce the above integral representation it is useful to introduce the abbreviation 
\begin{equation}
\label{eq:FundamentalIdentityGeneralisedD}
D_{G,H}(t,v) = \pa_{2} \nE_G(t,v) \, H(v) - H\big(\nE_G(t,v)\big)\,. 
\end{equation}
Differentiation yields 
\begin{equation*}
\begin{split}
\ddt \, D_{G,H}(t,v) 
&= \ddt \, \pa_{2} \nE_G(t,v) \, H(v) - \ddt \, H\big(\nE_G(t,v)\big) \\
&= G'\big(\nE_G(t,v)\big) \, \pa_{2} \nE_G(t,v) \, H(v) - H'\big(\nE_G(t,v)\big) \, G\big(\nE_G(t,v)\big)\,,
\end{split}
\end{equation*}
see~\eqref{eq:EG}--\eqref{eq:dEG}, which further implies 
\begin{equation*}
\begin{split}
&\ddt \, D_{G,H}(t,v) - G'\big(\nE_G(t,v)\big) \, D_{G,H}(t,v) \\
&\qquad= G'\big(\nE_G(t,v)\big) \, \pa_{2} \nE_G(t,v) \, H(v) - H'\big(\nE_G(t,v)\big) \, G\big(\nE_G(t,v)\big) \\
&\qquad\qquad- G'\big(\nE_G(t,v)\big) \, \pa_{2} \nE_G(t,v) \, H(v) + G'\big(\nE_G(t,v)\big) \, H\big(\nE_G(t,v)\big) \\
&\qquad= G'\big(\nE_G(t,v)\big) \, H\big(\nE_G(t,v)\big) - H'\big(\nE_G(t,v)\big) \, G\big(\nE_G(t,v)\big) \\
&\qquad= \big[G,H\big]\big(\nE_G(t,v)\big)\,. 
\end{split}
\end{equation*}
The resulting initial value problem 
\begin{equation}
\label{eq:DIvp}
\begin{cases}
&\ddt \, D_{G,H}(t,v) = G'\big(\nE_G(t,v)\big) \, D_{G,H}(t,v) + \big[G,H\big]\big(\nE_G(t,v)\big)\,, \\
&D_{G,H}(0,v) = 0\,, 
\end{cases}
\end{equation}
\end{subequations}
can be cast into the form~\eqref{eq:Cr} with time-dependent linear operator $C(t) = G'\big(\nE_G(t,v)\big)$, associated resolvent $\nE_{C}(t,0) = \pa_{2} \nEG(t,v)$, 
remainder $r(t) = \big[G,H\big]\big(\nE_G(t,v)\big)$, and vanishing initial value. 
By the linear variation-of-constants formula the stated relation follows. 

\paragraph{Nonlinear variation-of-constants formula.}
A basic tool in the context of nonlinear evolution equations is the Gröbner--Alekseev formula;
as a generalisation of the linear variation-of-constants formula (Duhamel principle) it relates the solutions to two nonlinear evolution equations through an integral representation.  
For the convenience of the reader we recall the Gröbner--Alekseev formula adapted to the present situation. 

\begin{theorem}[Gröbner--Alekseev formula]
\label{thm:Groebner}
The solutions to the autonomous problem 
\begin{equation*}
\begin{cases}
&\tfrac{\dd}{\dd t} \, \nE_G(t - t_0,v) = G\big(\nE_G(t - t_0,v)\big)\,, \quad t_0 \leq t \leq T\,, \\
&\nE_G(0,v) = v\,, 
\end{cases}
\end{equation*}
and the related non-autonomous problem 
\begin{equation*}
\begin{cases}
&\tfrac{\dd}{\dd t} \, \nE_{G+R}(t,t_0,v) = G\big(\nE_{G+R}(t,t_0,v)\big) + R(t)\,, \quad t_0 \leq t \leq T\,, \\
&\nE_{G+R}(t_0,t_0,v) = v\,,
\end{cases}
\end{equation*}
satisfy the integral relation 
\begin{equation*}
\nE_{G+R}(t,t_0,v) - \nE_G(t-t_0,v) = \int_{t_0}^{t} \pa_{2} \nE_G\big(t-\t,\nE_{G+R}(\t,t_0,v)\big) \, R(\t) \; \dd\t\,.
\end{equation*}
\end{theorem}
\paragraph{Proof}
We recall the fundamental identity~\eqref{eq:FundamentalIdentity} which implies 
\begin{equation*}
\pa_{2} \nE_G(t-\t,w) \, G(w) = G(\nE_G(t-\t,w))\,, \qquad w = \nE_{G+R}(\t,t_0,v)\,.
\end{equation*}  
A brief calculation shows 
\begin{equation*}
\begin{split}
\nE_{G+R}(t,t_0,v) - \nE_G(t-t_0,v) 
&= \nE_G\big(t-\t,\nE_{G+R}(\t,t_0,v)\big)\Big\vert_{\t=t_0}^{t} \\
&= \int_{t_0}^{t} \tfrac{\dd}{\dd \t} \, \nE_G\big(t-\t,\nE_{G+R}(\t,t_0,v)\big) \; \dd\t \\
&= \int_{t_0}^{t} \Big(\pa_{2} \nE_G\big(t-\t,\nE_{G+R}(\t,t_0,v)\big) \, \big(G\big(\nE_{G+R}(\t,t_0,v)\big) + R(\t)\big) \\
&\qquad\qquad- G\big(\nE_G\big(t-\t,\nE_{G+R}(\t,t_0,v)\big)\big)\Big) \; \dd\t \\
&= \int_{t_0}^{t} \pa_{2} \nE_G\big(t-\t,\nE_{G+R}(\t,t_0,v)\big) \, R(\t) \; \dd\t\,,
\end{split}
\end{equation*}
which proves the statement. 
$\hfill \diamondsuit$ 
\subsection{Local error expansion} 
\paragraph{Approach.}
In the following, we consider the splitting operator associated with the first-order Lie--Trotter splitting method~\eqref{eq:Lie} as time-continuous function 
\begin{equation}
\label{eq:SLie}
\nSF(t,v) = \nEB\big(t,\nEA(t,v)\big)\,.
\end{equation}
We aim for the derivation of suitable evolution equations such that appropriate integral representations for their solutions finally yield the desired local error representation
\begin{equation*}
\nLF(t,v) = \nSF(t,v) - \nEF(t,v) = \nO\big(t^2\big).
\end{equation*}
Evidently, the alternative case $\nSF(t,v) = \nEA\big(t,\nEB(t,v)\big)$ is covered by exchanging the roles of~$A$ and~$B$, see also~\eqref{eq:Lie}.  
\subsubsection{Definition and reformulation of defect} 
\paragraph{Definition of defect.} 
The defect associated with the Lie--Trotter splitting method~\eqref{eq:SLie} is defined by 
\begin{equation*}
R_1(t,v) = \ddt \, \nSF(t,v) - F\big(\nSF(t,v)\big)\,. 
\end{equation*}

\paragraph{Reformulation of defect.} 
In order to obtain a first reformulation of the defect, we determine the derivative of the splitting operator. 
By means of the chain rule and the identity $\ddt \, \nEG(t,w) = G(\nEG(t,w))$ the relation 
\begin{equation*}
\begin{split}
\ddt \, \nSF(t,v)
&= \ddt \, \nEB\big(t,\nEA(t,v)\big) \\
&= B\big(\nEB\big(t,\nEA(t,v)\big)\big) + \pa_{2} \nEB\big(t,\nEA(t,v)\big) \, A\big(\nEA(t,v)\big) \\
&= B\big(\nSF(t,v)\big) + \pa_{2} \nEB\big(t,\nEA(t,v)\big) \, A\big(\nEA(t,v)\big) 
\end{split}
\end{equation*}
follows, which further implies 
\begin{equation}
\label{eq:R1Lie}
R_1(t,v) = D_{B,A}\big(t,\nEA(t,v)\big)\,,
\end{equation}
as shown by a brief calculation
\begin{equation*}
\begin{split}
R_1(t,v) 
&= \ddt \, \nSF(t,v) - F\big(\nSF(t,v)\big) \\
&= B\big(\nSF(t,v)\big) + \pa_{2} \nEB\big(t,\nEA(t,v)\big) \, A\big(\nEA(t,v)\big) - A\big(\nSF(t,v)\big) - B\big(\nSF(t,v)\big) \\ 
&= \pa_{2} \nEB\big(t,\nEA(t,v)\big) \, A\big(\nEA(t,v)\big) - A\big(\nSF(t,v)\big) \\
&= D_{B,A}\big(t,\nEA(t,v)\big)\,; \\
\end{split}
\end{equation*}
recall the decomposition $F = A + B$ and the abbreviation~\eqref{eq:FundamentalIdentityGeneralisedD}.
\subsubsection{Integral representation ensuring $\nLF(t,v) = \nO(t)$} 
\paragraph{Initial value problems for evolution and splitting operator.} 
Evidently, the evolution operator fulfills the initial value problem 
\begin{equation*}
\begin{cases}
&\ddt \, \nEF(t,v) = F\big(\nEF(t,v)\big)\,, \\
&\nEF(0,v) = v\,.
\end{cases}
\end{equation*}
The notion of the defect $R_1 = \ddt \nSF - F(\nSF)$ allows to consider the splitting operator as solution to the related initial value problem 
\begin{equation*}
\begin{cases}
&\ddt \, \nSF(t,v) = F\big(\nSF(t,v)\big) + R_1(t,v)\,, \\
&\nSF(0,v) = v\,.
\end{cases} 
\end{equation*}

\paragraph{Integral representation for local error.} 
The application of the nonlinear variation-of-constants formula (Gröbner--Alekseev formula) yields a first integral representation for the local error 
$\nLF(t,v) = \nSF(t,v) - \nEF(t,v) = \nE_{F+R_1}(t,0,v) - \nEF(t,v)$, namely 
\begin{equation}
\label{eq:LOtLie}
\nLF(t,v) = \int_{0}^{t} \pa_{2} \nEF\big(t-\t_1,\nSF(\t_1,v)\big) \, R_1(\t_1,v) \; \dd\t_1\,,
\end{equation}
see also Theorem~\ref{thm:Groebner}.
Provided that the integrand remains bounded on the underlying Banach space, this ensures $\nLF(t,v) = \nO(t)$.
\subsubsection{Integral representation ensuring $\nLF(t,v) = \nO(t^2)$} 
\paragraph{Integral representation for defect.} 
The representation~\eqref{eq:FundamentalIdentityGeneralised} implies the following integral representation for the defect
\begin{equation}
\label{eq:R1OtLie}
\begin{split}
R_1(\t_1,v) 
&= D_{B,A}\big(\t_1,\nEA(\t_1,v)\big) \\
&= \pa_{2} \nEB(\t_1,w) \int_{0}^{\t_1} \big(\pa_{2} \nEB(\t_2,w)\big)^{-1} \big[B,A\big]\big(\nEB(\t_2,w)\big) \; \dd\t_2 \, \Big\vert_{w = \nEA(\t_1,v)}\,,
\end{split}
\end{equation}
see also~\eqref{eq:R1Lie}.

\paragraph{Local error expansion.} 
Altogether, inserting the integral representation~\eqref{eq:R1OtLie} into~\eqref{eq:LOtLie} we obtain the local error expansion 
\begin{equation}
\label{eq:LLie}
\begin{split}
\nLF(t,v) 
&= \int_{0}^{t} \int_{0}^{\t_1} \pa_{2} \nEF\big(t-\t_1,\nSF(\t_1,v)\big) \, \pa_{2} \nEB(\t_1,w) \, \big(\pa_{2} \nEB(\t_2,w)\big)^{-1} \\
&\qquad\qquad\quad \times \big[B,A\big]\big(\nEB(\t_2,w)\big) \Big\vert_{w = \nEA(\t_1,v)} \; \dd\t_2 \, \dd\t_1\,,
\end{split}
\end{equation}
which ensures $\nLF(t,v) = \nO(t^2)$ provided that the integrand remains bounded in the underlying function space. 
\\[0.2cm]

\end{document}